\documentclass[
]{amsart}

 \usepackage{tikz}
 \usetikzlibrary{arrows,calc}
 \usepackage[export]{adjustbox}
 \usepackage{float}
 \usetikzlibrary{backgrounds,arrows.meta}

 \usepackage{chngcntr}
\usepackage[belowskip=-8pt,aboveskip=6pt]{caption}

\usepackage{pgfplots}
\pgfplotsset{compat=newest}
\usepgfplotslibrary{fillbetween}

\usepackage[utf8]{inputenc} 
\usepackage[T1]{fontenc}

\usepackage{amsmath}
\usepackage{amsfonts}
\usepackage{amssymb}
\usepackage{amsthm}
\usepackage{setspace}
\usepackage{amsrefs}

\usepackage{graphicx}

\usepackage{hyperref}
\hypersetup{hidelinks}

\title[]
{On regularity of $\overline\partial$-solutions on  $a_q$ domains with $C^2$ boundary in complex manifolds
}
\author[]{Xianghong Gong$^{\dag}$}

\date{\today}

 \address{Department of Mathematics,
 University of Wisconsin-Madison, Madison, WI 53706}
 \email{gong@math.wisc.edu}

\thanks{$^{\dag}$Partially supported by Simons
Foundation grant (award number:~505027) and NSF grant DMS-2054989}


\makeatletter
\newcommand*\bigcdot{\mathpalette\bigcdot@{.75}}
\newcommand*\bigcdot@[2]{\mathbin{\vcenter{\hbox{\scalebox{#2}{$\m@th#1\bullet$}}}}}
\makeatother

\newcommand{\dist}{\operatorname{dist}}

\newtheorem{thm}{Theorem}[section]
\newtheorem{cor}[thm]{Corollary}
\newtheorem{prop}[thm]{Proposition}
\newtheorem{lemma}[thm]{Lemma}
\newtheorem{remark}[thm]{Remark}

\newcommand{\dv}{d\upsilon}
\newcommand{\ds}{d\sigma}

\theoremstyle{definition}

\newtheorem{defn}[thm]{Definition}

\newtheorem{rem}[thm]{Remark}

\newtheorem{notation}[thm]{Notation}

\renewcommand{\th}[1]{\begin{thm}\label{#1}}
\renewcommand{\eth}{\end{thm}}
\newcommand{\co}[1]{\begin{cor}\label{#1}}
\newcommand{\eco}{\end{cor}}
\renewcommand{\le}[1]{\begin{lemma}\label{#1}}
\newcommand{\ele}{\end{lemma}}
\newcommand{\pr}[1]{\begin{prop}\label{#1}}
\newcommand{\epr}{\end{prop}}

\newcommand{\AL}[1]{\begin{align}
#1
\end{align}
}

\newcommand{\ga}{\begin{gather}}
\newcommand{\ega}{\end{gather}}
\newcommand{\gan}{\begin{gather*}}
\newcommand{\egan}{\end{gather*}}
\newcommand{\al}{\begin{align}}
\newcommand{\eal}{\end{align}}
\newcommand{\aln}{\begin{align*}}
\newcommand{\ealn}{\end{align*}}
\newcommand{\eq}[1]{\begin{equation}\label{#1}}
\newcommand{\eeq}{\end{equation}}

\newcommand{\DD}[2]{\frac{\partial #1}{\partial #2}}
\newcommand{\f}[2]{\frac{#1}{#2}}

\newcommand{\ip}[1]{\left<#1\right>}

\newcommand{\ci}{~\cite}

\newcommand{\cc}{{\bf C}}
\newcommand{\nn}{{\bf N}}

\newcommand{\rr}{{\bf R}}

\newcommand{\pp}{{\bf P}}

\newcommand{\ov}{\overline}

\newcommand{\RE}{\operatorname{Re}}

\renewcommand{\dbar}{\overline\partial}

\newcommand{\cL}{\mathcal}

\newcommand{\all}{\alpha}

\newcommand{\gaa}{\gamma}

\newcommand{\del}{\delta}
\newcommand{\Del}{\Delta}
\newcommand{\var}{\varphi}
\newcommand{\e}{\epsilon}
\newcommand{\om}{\omega}
\newcommand{\Om}{\Omega}

\newcommand{\la}{\lambda}

\newcommand{\pd}{\partial}

\newcommand{\jq}[1]{\langle #1\rangle}

\newcommand{\re}[1]{(\ref{#1})}
\newcommand{\rea}[1]{$(\ref{#1})$}
\newcommand{\rl}[1]{Lemma~\ref{#1}}
\newcommand{\nrc}[1]{Corollary~\ref{#1}}
\newcommand{\rp}[1]{Proposition~\ref{#1}}
\newcommand{\rt}[1]{Theorem~\ref{#1}}
\newcommand{\rd}[1]{Definition~\ref{#1}}

\newcommand{\rf}[1]{Figure~\ref{#1}}

\newcommand{\rla}[1]{Lemma~$\ref{#1}$}
\newcommand{\nrca}[1]{Corollary~$\ref{#1}$}
\newcommand{\rpa}[1]{Proposition~$\ref{#1}$}
\newcommand{\rta}[1]{Theorem~$\ref{#1}$}

\newcommand{\supp}{\operatorname{supp}}

\newcommand{\db}{\dbar}

\newcounter{pp}
\newcommand{\bpp}{\begin{list}{$\hspace{-1em}(\alph{pp})$}{\usecounter{pp}}}
\newcommand{\epp}{\end{list}}

\newcounter{ppp}
\newcommand{\bppp}{\begin{list}{$\hspace{-1em}(\roman{ppp})$}{\usecounter{ppp}}}
\newcommand{\eppp}{\end{list}}

\def\beq{\begin{equation}}
\def\eeq{\end{equation}}

\begin{document}

\begin{abstract}
We study regularity of solutions $u$ to   $\overline\partial u=f$ on a relatively compact $C^2$ domain $D$ in a complex manifold of dimension $n$, where $f$ is a $(0,q)$ form.  Assume  that there are  either $(q+1)$ negative or $(n-q)$ positive Levi eigenvalues at each  point  of boundary $\partial D$.
Under the necessary condition that a locally $L^2$ solution exists  on the domain,
we show  the existence of the solutions on the closure of the domain that gain $1/2$ derivative when  $q=1$ and $f$ is in the H\"older-Zygmund space $\Lambda^r( D)$ with $r>1$. For $q>1$, the same regularity  for the solutions is achieved when   $\partial D$ is either sufficiently smooth  or  of $(n-q)$ positive Levi eigenvalues everywhere on $\partial D$.
\end{abstract}

 \maketitle


\setcounter{thm}{0}\setcounter{equation}{0}

\section{Introduction}\label{sec1}

Let $D$ be a relatively compact  domain in a complex manifold $X$ of dimension $n$. We say that $D$ satisfies the {\it condition} $a_q$   if boundary  $\pd D\in C^2$ and  its Levi-form  on complex tangent space $T_\zeta^{(1,0)}(\partial D)$   has   either $(q+1)$ negative or $(n-q)$ positive eigenvalues at each   $\zeta\in\pd D$. We are interested in the regularity of solutions $u$ to the $\db$-equation, $\db u=f$, on $\ov D$.
 We will study the case where $f$ is a $V$-valued $(0,q)$ form for a holomorphic vector bundle $V$ on $X$. Denote by $\Lambda_{(0,q)}^r({D}, V)$ the space of $V$-valued $(0,q)$ forms whose coefficients are in H\"older--Zygmund space $\Lambda^r({D})$ defined in Sect.~\ref{h-space}; note that $\Lambda^r(D)$ is the  H\"older space $C^r(\ov D)$ when $r\in(0,\infty)\setminus\nn$.
To ensure the existence of solutions,  we impose a minimum requirement that there is   an $L_{loc}^2$ solution $u_0$ on $D$  and seek a possibly different  solution   of better regularity.

Our main results are the following.
\th{regsol} Let $r\in(1,\infty)$ and $q\geq 1$.
Let $D $ be a relatively compact domain with $C^2$ boundary in a complex manifold $X$ satisfying the condition $a_q$. Let $V$ be a holomorphic vector bundle on $X$.
 Then there exists a bounded
linear $\db$-solution operator $H_q\colon \Lambda_{(0,q)}^r({D}, V)\cap\db L^2_{loc}(D)\to \Lambda_{(0,q-1)}^{r+1/2}({D}, V)$, 
 provided $(a)$
$q=1$ or $\pd D$ has $(n-q)$ positive Levi eigenvalues at each point on $\pd D$;  or $(b)$
$\pd D\in\Lambda^{r+\f{5}{2}}$.
\eth
Note that $H_q$ is independent of $r$ and it provides a smooth ($C^\infty$) linear $\db$-solution operator for smooth forms in both cases.
When $\pd D\in C^2$,   case (a) provides a satisfactory regularity result for $\db$-solutions in the H\"older--Zygmund spaces for $q=1$. For $q>1$, we have the following.
\th{regsol+} Let   $q\geq 2$ and keep notations in \rta{regsol} with $\pd D\in C^2$. Then there exists a finite $r_0$ such that if $r>r_0$, there exists a bounded
linear $\db$-solution operator $$H_q^{r}\colon \Lambda_{(0,q)}^r( D, V)\cap \db L^2_{loc}(D)\to \Lambda_{(0,q-1)}^{r-5/2}( D).$$ Further, $H_q^{r}$ maps $C_{(0,q)}^\infty(\ov D, V)\cap \db L^2_{loc}(D) $ into $ C_{(0,q-1)}^{\infty}(\ov D)$.
\eth
The value $r_0$ is stated in a detailed version of \rt{regsol+} in \rt{nash-moser-c12}.

\medskip

We  first state  some closely related  results on $\db$-solutions $u$ to $\db u=f$ on strictly pseudoconvex domains $D$ in $\cc^n$:  After  work of  Lieb--Grauert~\cite{MR273057} and Kerzman~\cite{MR0281944}, Henkin--Romanov~\cite{MR0293121} achieved the sharp $C^{1/2}$ solutions for $f\in L^\infty$  by integral formulas. The $C^{k+1/2}$ solutions for  $f\in C^k$ $(k\in\nn)$ were obtained by Siu~\cite{MR330515} for $(0,1)$ forms   and by Lieb--Range~\cite{MR597825}  for forms with $q\geq1$ when $\pd D$ is sufficiently smooth.
 For  $\pd D\in C^2$, Theorem~\ref{regsol} and analogous results for a homotopy formula were proved in~\cite{MR3961327}  through the construction of a homotopy formula. These results were extended by Shi~\cite{MR4244873} to a weighted Sobolev spaces with a gain less than $1/2$ derivative and by  Shi--Yao~\cites{MR4688544,shi-yao-ajm} to  $H^{s+1/2,p}$ space gaining $1/2$ derivative  for $s>1/p$ when $\pd D\in C^2$  and  for $s\in\rr$ when $\pd D$ is sufficiently smooth. It is worthy to point out that Shi--Yao achieved the first regularity result for negative order $s$, and Yao~\cite{MR4739361} further showed that a similar result holds for convex finite multitype smooth domains in $\cc^n$. Also, Gong--Lanzani~\cite{MR4289246} obtained $\Lambda^{r+1/2}$ (with $r>1$) regularity gaining $1/2$ derivative on strongly $\cc$-linear convex $C^{1,1}$ domains $D$.

Next, we mention results on $\db$ solutions on $a_q$ domains in complex manifolds, which are also called $Z(q)$ domains in ~\cite{MR0461588}*{p.~57}).
 The basic   estimate for $\db$ was proved by Morrey~\cite{MR0099060} for $(0,1)$-forms
and by Kohn~\cites{MR0153030,MR0208200} for forms of any type on strongly pseudoconvex manifolds with smooth boundary. These results lead to  regularity of the $\db$-Neumann operator for strictly pseudoconvex manifolds and more general for compact manifolds whose boundary satisfies  property $Z(q)$ (see~\cites{MR0153030,MR0208200} and \cite{MR0177135}*{Thm.~3.9}).  By work of H\"ormander~\cite{MR0179443}, the basis estimate is equivalent to the condition $Z(q)$ on $D$ (see~\cite{MR0177135}*{Prop.~3.12} and \cite{MR0461588}*{Thm.~3.2.2} for details).  When $\pd D\in C^\infty$, sharp regularity results for $\db$ solutions were obtained by Greiner-Stein~\cite{MR0499319} for $(0,1)$ forms and Beals--Greiner--Stanton~\cite{MR886418} for $(0,q)$ forms under condition $Z(q)$ through the study of the regularity of $\db$-Neumann operator in $L^{k,p}$ and $\Lambda^r$ spaces.
The condition $a_q$ also ensures the stability of the solvability of the $\db$-equation on $(0,q)$ forms; namely if  $f=\db u_0+\tilde f$ on $D$ while $\tilde f$ is a $\db$-closed form on a larger domain and $u_0$ has the regularity in the sought-after class, then $\tilde f=\db\tilde u$ for some $L^2$ form $\tilde u$. This stability is useful to obtain regularity for $\db$ solutions as shown by Kerzman~\cite{MR0281944}; first one seeks regularity for $u_0$ without solving the $\db$-equation. Then $u_0+\tilde u$ provides a desired solution based on regularity of $u_0$ and  the interior regularity of $\tilde u$ from the elliptic theory on systems of partial differential equations.

\medskip

To prove our results, we will use integral formulas to obtain local solutions near each boundary point of $D$. We then use the Grauert bumping method as in~\cite{MR0281944} to construct $\tilde f$.    To provide background for our results, let us mention regularity results for $\db$-solutions on the transversal intersection of domains in $\cc^n$.   Range--Siu~\cite{MR338450} obtained $C^{1/2-\epsilon}$ estimate with any positive $\e$ for a real transversal intersection of sufficiently smooth strictly pseudoconvex domains. For higher order derivatives,
J. Michel~\cite{MR928297}  obtained $C^{k+1/2-\e}$ estimate for $\db$ solutions on a certain intersection of smooth strictly pseudoconvex domains. J.~Michel and Perotti~\cite{MR1038709} extended the result to real transversal intersection of  strictly pseudoconvex domains with sufficiently smooth boundary.
We should also mention that the local version of \rt{regsol} was proved by Laurent-Thi\'ebaut and Leiterer \cite{MR1207871} when $\pd D\in C^\infty$ and $k\in \nn$.
  Ricard~\cite{MR1992543} obtained regularity for concave wedge with $C^{k+2}$ boundary and convex wedges with $C^2$ boundary. The reader is referred to
 Barkatou~\cite{MR1888228}  and
Barkatou-Khidr~\cite{MR2844676}  for further results in this direction.
 However, all existing integral formulas for $\db$ solutions, including ours for $q>1$, require boundary to be sufficiently smooth when concavity is present.

\smallskip

 It is well-known that concavity of the domains is useful; the classical Hartogs theorem says that a holomorphic function on a $1$-concave domain extends to a holomorphic function across the boundary. Therefore, on a $1$-concave domain in a complex manifold, all $\db$-solutions for $(0,1)$ forms must have the same regularity regardless the smoothness of the boundary of the domain. In Section 7 we will successfully implement this idea to prove \rt{regsol} $(i)$ with $q=1$.

 When $q>1$, not all $\db$ solutions have the same regularity.
 To prove \rt{regsol+}, we first derive an estimate for a solution operator when $\pd D$ is sufficiently smooth.  Then we apply a Nash-Moser iteration method by solving the $\db$-equation  on the subdomains $D_k\subset D_{k+1}$ that have smooth boundary and in the limit, we obtain a desired solution on the closure of $D=\cup D_k$.

\smallskip

We organize the paper as follows.

In Section 2, we formulate an approximate local homotopy formula on a suitable neighborhood of a boundary point $\zeta_0\in\pd \Om$. In Sections 3 and 4. we derive (genuine) local homotopy formulas for $\db$-closed $(0,q)$ forms near $(n-q)$ convex and $(q+1)$ concave boundary points of an $a_q$ domain. There we follow approaches developed in Lieb--Range~\cite{MR597825} and Henkin--Leiterer~\cite{MR986248}.
While a homotopy formula for forms that are not necessary $\db$-closed can be derived for $(n-q)$ convex points without extra conditions,  such a formula on the concave side of the boundary turns out to be subtle. In fact, Laurent-Thi\'ebaut and Leiterer~\cite{MR1621967}*{Prop.~0.7}  proved that there is no {\it local} homotopy formula near a $(q+1)$ concave boundary for $(0,q)$ forms that has {\it good} estimates, and see also Nagel--Rosay~\cite{MR1016447} on $\db_b$ for domains in sphere in $\mathbb C^3$.  Such phenomena leads to difficulties for strictly pseudoconvex hypersurfaces for local CR embedding problem in Webster~\cite{MR995504} and Polykov~\cite{MR2088929} on global CR embeddings of concave compact CR manifolds in critical dimensions.

Section 5 contains some elementary facts on H\"older--Zygmund spaces, where we derive an equivalence characterization on the H\"older--Zygmund norms solely relying on a version of Hardy-Littlewood lemma.

Section 6 contains main local estimates for homotopy operators that appear in the local homotopy formula. One of main purposes of the section is to derive precise estimates that reflex the convexity of the H\"older-Zygmund norms; see \re{hqfr12-c} for strictly $(n-q)$ convex  $C^2$ domains  and \re{hqfr12} for strictly $(q+1)$-concave domains. We emphasize that the estimates do not require the forms to be $\db$-closed; in fact, the estimates hold for $(q+1)$ concave boundary points. These estimates immediately give us the desired regularity stated in \rt{regsol} for local solutions.

In Section 7   we show as a novelty how the Hartogs extension theorem can be used to study the regularity of $\db$ solutions for $(0,1)$ forms.
A local version of \rt{regsol} (a) for $q=1$ is proved in this section.

In Section 8 we show the existence of global solutions with the desired regularity by using local solutions in Sections 6 and 7 and the interior regularity of elliptic systems. We also derive a global estimate for $\db$-solutions.
Using this global estimate, we employ the Nash--Moser smoothing operator to prove a detailed version of \rt{regsol+} in Section 9.

The paper has two appendices. In Appendix A, we recall the existing regularity on the signed distance function near a $C^2$ hypersurface in a Riemannian manifold. In Appendix B, we describe a stability result of solvability of $\db$-equation on $a_q$ domains with $C^2$ boundary using results in H\"ormander~\cite{MR0179443}.

\medskip
{\bf Acknowledgments.} The author is grateful to the anonymous referee for many  helpful suggestions and a proof of \re{eq:abcd}.

 \setcounter{equation}{0}

\section{A local approximate
 homotopy formula}

 Let $X$ be a complex manifold of dimension $n$.
Let $D$ be a relatively compact  domain in $X$ defined by $\rho<0$, where $\rho$ is a $C^2$ defining function  with $\nabla\rho(\zeta)\neq0$ when $\rho(\zeta)=0$. Throughout the paper,  $\nabla^k\rho$ is the set of $k$-th partial derivatives in local coordinates. We recall the following definitions  in Henkin--Leiterer~\cite{MR986248} and H\"ormander~\cite{MR0179443}.
\begin{defn} Let $D\subset X$ be a $C^2$ domain defined by $\rho<0$  as above.
\bpp
\item For $\zeta\in\pd D$, the {\it  Levi-form}
$L_\zeta\rho$ in local holomorphic coordinates $z$ on $X$ is   the complex Hessian
$$H_\zeta\rho(w)=\sum_{j,k=1}^n\DD{^2H}{ z_j\pd \bar z_k}(\zeta)w_j\bar w_k,\quad w\in\cc^n
$$
restricted  on the complex tangent space  $T_\zeta^{1,0}(\pd D)$, where the latter is identified with $\{w\in\cc^n\colon\sum\DD{\rho}{z_j}(\zeta) w_j=0\}$.
\item $\pd D$ is {\it strictly $q$-convex} (resp. {\it $q$-concave})  at $\zeta\in\pd D$,  if $L_\zeta\rho$  has  at least $q$ positive (resp. negative) eigenvalues.
    \item   $D$ satisfies the {\it condition $a_q$}  if   $L_\zeta\rho$ has at least either $(q+1)$ negative eigenvalues or $(n-q)$ positive eigenvalues  for every $\zeta\in\pd D$.
  \epp
  \end{defn}
Therefore,  a  domain is strictly pseudoconvex if and only if it is strictly $(n-1)$-convex.
 To see a domain satisfying the condition $a_q$,
let $\pp^n$ be the   complex project space. Let $B_{q}^r\subset\pp^n$ be defined by
$$
B_{q}^r\colon |z^{q}|^2+\cdots+ |z^n|^2< r|z^0|^2+\cdots+ r|z^{q-1}|^2
$$
with $1\leq q\leq n$ and $r>0$. Then $B^r_{q}$ is  both strictly $(n-q)$-convex and  $(q-1)$-concave.
 When $\pd D$ is Levi non-degenerate, the condition $a_q$ for $\pd D$ is equivalent to the number of negative Levi eigenvalue not being $q$ at every point $\zeta\in\pd D$. Thus  $B^{r_2}_{q_1}  \setminus\ov{  B^{r_1}_{q_2}}$   satisfies the condition $a_q$ if $r_2>r_1$, $q_2\geq q_1$, and
$$ q\in\{0,\dots, n-1\}\setminus\{ q_1-1, n-q_2-1\}.
$$

We first fix  notation. If $E^1,\dots, E^k$ are sets, we write $E^{1\dots k}=E^1\cap\cdots\cap E^k$.

 The main purpose of this section is to construct an approximate homotopy formula on a subdomain $D'$ of $D$, where $\pd D'$ and $\pd D$ share a piece of boundary containing a given point $\zeta_0\in\pd D$. The $D'$ will have the form $D^{12}=D^1\cap D^2$ where $D^1=D$ and $D^2$ is a ball in suitable coordinates for $\pd D$.
Since the construction of $D'$ is local, we   assume that $D$ is  contained in $\cc^n$.

\begin{defn}\label{pck} Let $k\geq1$ and $r>1$.
 Let $k\geq1$ and $r>1$. Let $D^1,\dots, D^\ell$ be open sets in a complex manifold $X$.
We say that  $\om:=D^{1\dots\ell}$ is a {\it  $C^k$
 $($resp. $\Lambda^r$ with $ r>1)$   transversal intersection} (of $D^1,\dots, D^\ell$), if  $\ov\om$ is compact, there are $C^k$ (resp. $\Lambda^r$) real functions $\rho_1,\dots,\rho_\ell$ on a neighborhood $U$ of  $\ov\om$ such that $D^j\subset U\colon\rho_j<0$, and for any $1\leq j_1<\cdots<j_i\leq\ell$,
$$
d\rho_{j_1}(\zeta)\wedge\cdots\wedge d\rho_{j_i}(\zeta) \neq0, \quad\text{when $\rho_{j_1}(\zeta)=\cdots=\rho_{j_i}(\zeta)=0$}.
$$
\end{defn}
We will   need Leray maps, following notation in \cite{MR986248}.
\begin{defn} Let $D$ be a domain in $\cc^n$. Let $S\subset \cc^n\setminus D$ be a $C^1$ submanifold in $\cc^n$. We say that $g\colon D\times S\to \cc^n$ is a {\it Leray map}   if $g\in C^1(D\times S)$ and
$$
g( z,\zeta )\cdot(\zeta-z)\neq 0, \quad \zeta\in S, \quad z\in D.
$$
 Throughout the paper, we use
\eq{g0}
g_0( z,\zeta )=\ov\zeta-\ov z.
\end{equation}
\end{defn}

   Let $g^j \colon D\times S^j\to\cc^n$ be   $C^1$ Leray   mappings for $j=1,\dots,\ell$.
Let $w=\zeta-z$.      Define  \gan
\omega^i=\f{1}{2\pi i}\f{g^i\cdot dw}{g^i\cdot w},
\quad
\Omega^i=\omega^i\wedge(\ov\pd\omega^i)^{n-1},\\
\Omega^{01}=\omega^0\wedge\omega^1\wedge\sum_{\alpha+\beta=n-2}
(\ov\pd\omega^0)^{\alpha}\wedge(\ov\pd\omega^1)^{\beta}.
\end{gather*}
Here both differentials $d$ and $\db$ are in $z,\zeta$ variables.
In general, define
\gan
\Omega^{1\cdots \ell}=\omega^{g_1}\wedge\cdots\wedge\omega^{g_\ell}\wedge\sum_{\alpha_1+\cdots+\all_\ell=n-\ell}
(\ov\pd\omega^{g_1})^{\alpha_1}\wedge\cdots(\ov\pd\omega^{g_\ell})^{\all_\ell}.
\end{gather*}
Decompose $\Om^{\bigcdot}=\sum\Om_{(0,q)}^{\bigcdot}$,
where $\Om_{(0,q)}^{\bigcdot}$
  has type $(0,q)$  in $z$. Hence $\Om^{i_1,\dots, i_\ell}_{(0,q)}$ has type $(n,n-\ell-q)$    in $\zeta$. Set
   $\Om^{\bigcdot}_{0,-1}=0$ and $\Omega_{0,n+1}^{\bigcdot}=0$. The Koppelman lemma says that
   $$
   \db_z\Om^{1\cdots \ell}_{(0,q-1)}+\db_\zeta\Om^{1\cdots \ell}_{(0,q)}=\sum(-1)^{j}\Om_{(0,q)}^{1\cdots\hat j\cdots\ell}.
   $$
See Chen-Shaw~\cite{MR1800297}*{p.~263} for a proof. We will use  special cases:
\ga
\label{kop1}\db_\zeta\Om_{(0,q)}^1+\db_z\Om_{(0,q-1)}^1=0, 
\\
\label{kop2}\db_\zeta\Om^{01}_{(0,q)}+\db_z\Om_{(0,q-1)}^{01}=-\Om_{(0,q)}^1+\Om_{(0,q)}^0, 
\\
\label{kop3}\db_\zeta\Om^{012}_{(0,q)}+\db_z\Om_{(0,q-1)}^{012}=-\Om_{(0,q)}^{12}+\Om_{(0,q)}^{02}-\Om^{01}_{(0,q)}.
\end{gather}
Here each identity holds in the sense of distributions on the set  where the forms are non-singular.

To integrate on   submanifolds of   $\cc^n$, let us see how a sign changes when an exterior differentiation interchanges with integration. Following~\cite{MR1800297}*{p.~263},   define
$$
  \int_{y\in M} u(x,y)dy^J\wedge dx^I =  \left \{\int_{y\in M}u (x,y)dy^J\right\}dx^I
$$
for a   function $u$ on a manifold $M$ with boundary. For  the exterior differential $d_x$, we have
$$ 
d_x\int_{y\in M}\phi(x,y) =(-1)^{\dim M}\int_{y\in M}d_x\phi(x,y).
$$ 
Then,  Stokes' formula has the form
\ga
\int_{y\in\pd M}\phi(x,y)\wedge \psi(y)=\int_{y\in M}\Bigl\{ d_y \phi(x,y)\wedge \psi(y)+
(-1)^{\deg \phi} \phi(x,y)\wedge d\psi(y)\Bigr\},
\label{checksign+}
\end{gather}
where $\deg\phi$ is the total degree of $\phi$ in $(x,y)$.
 When $D^{1\dots\ell}$ is a   $C^1$ transversal intersection, we choose orientations so that
$$
\int_{D^{1\dots\ell}}\, df=\sum_{i=1}^\ell\int_{\ov D^{1\dots\ell}\cap\pd D^i}f.
$$
Suppose $D^{12}$ is a  $C^1$  transversal intersection. Then we define
\eq{d12s}
S^i:=\ov {D^{12}}\cap \pd D^i, \quad \pd{D^{12}}=S^1\cup S^2;
\quad S^{12}:=\pd S^1, \quad S^{21}:=\pd S^2.
\end{equation}
Thus    Stokes' formula has the following special cases
$$
\int_{D^1\cap D^2}df=\int_{S^1}f+\int_{S^2}f,\quad
\int_{S^1}df=\int_{S^{12}}f, \quad
\int_{S^2}df=\int_{S^{21}}f.
$$
Next, we introduce integrals on  domains and lower-dimensional sets:
\ga \label{defnLR}
R_{D; q}^{i_1\dots i_\ell}f(z):=\int_{D}\Om^{i_1\dots i_\ell}_{(0,q)}(z,\zeta)\wedge f(\zeta), \quad
L_{i_1\cdots i_\mu; q}^{j_1\dots j_\nu}f:=\int_{S^{i_1\cdots i_\mu}}\Omega_{(0,q)}^{j_1\cdots j_\nu}\wedge f.
\end{gather}

Let $E_D\colon C^0(\ov D)\to C^0_0(\cc^n)$ be the Stein extension operator such that
$$ 
\|E_Du\|_{\Lambda^r(\cc^n)}\leq C_r(D)\|u\|_{\Lambda^r(D)}
$$ 
where the H\"older--Zygmund norm and its equivalent norms are defined in Section~\ref{h-space}. Note that the extension exists for any bounded Lipschitz domain $D$. See~\cite{MR3961327}*{Prop.~3.11} for a proof of the extension property and references therein.

\vspace{3ex}
\begin{figure}[!ht!]\counterwithin{figure}{section}
\centering
\begin{tikzpicture}[even odd rule]

 \begin{scope}
\clip  (0,0)  circle(2.5);
 \fill[gray!40] (0,0)  ellipse (4 and 1.75)  (0,-3)  ellipse (4 and 3);  


 \end{scope}

  \begin{scope}
\clip  (0,0)  circle(2.5);
 \fill[white]
 (0,-3)  ellipse (4 and 3); 
\end{scope}

  \draw[thick,domain=(1/8)*pi:(7/8)*pi,color=black, name path=D1]
 plot[smooth]({0+4*cos (\x r)},{-3
 +3*sin(\x r)});

  \draw[domain=(2.7/8)*pi:(5.3/8)*pi,color=black, name path=U]
 plot[smooth]({0+4*cos (\x r)},{0
 +1.75*sin(\x r)});  

  \draw[dashed, thick] (0,0) circle (2.5); 

\fill (0,0)  circle (.125);


   \node at (1.75,1.25) {$U^1$};
     \node at (-4.95,1.5) {$D^1\colon \rho^1<0,\  D^2\colon \rho^2<0$};
      \node at (-4.75,.8) {$\rho^2(z)=|z|^2-r^2_2$};

   \node at (0,-.75) {$S^1=\ov{D^{12}}\cap \pd D^{1}$};
      \draw[->] (-1.12,-.55)--(-1.12,-.175);

     \node at (0,-1.35) {$S^2=\ov{D^{12}}\cap \pd D^{2}$};
      \draw[->] (-1.325,-1.35)--(-2.0,-1.35);
        \node at (0,-2.25) {$D^2$};
      \node at (3.1,-1.5) {$D^{1}$};
       \node at (-.195,.65) {$S^1_+=\pd U^1\setminus S^1$};
      \draw[->] (-1.35,.75)--(-2.125,1.25);
 \draw[->] (-1.15,.95)--(-1.15,1.6);
\end{tikzpicture}
 \vspace{3ex}
    \caption{ $\pd D^{12}=S^1\cup S^2$ and $\pd U^1=S^1\cup S^1_+$}
    \vspace{2ex}
  \label{fig:U}
\end{figure}

The main purpose of this section is to
  derive the following approximate homotopy formula on a bounded domain  that is a  $C^1$ transversal intersection.

\begin{notation}
In all figures, when a domain is denoted by a letter such as $D^1, U^1$ in \rf{fig:U}, the letter is always placed inside the domain and next to the boundary of the domain. In all figures,  $D^2$ is   the ball $B_{r_2}:=\{z\in\cc^n\colon|z|^2<r_2^2\}$.
\end{notation}
\begin{prop}[see \rf{fig:U}] \label{hf} Let $D^{12}\subset\cc^n$ be a bounded domain of    $C^1$ transversal intersection. Let $S^1,S^2$ be given by \rea{d12s}. Let $U^1\subset\cc^n\setminus D^{1}$ be a bounded domain with   piecewise $C^1$ boundary such that $\pd U^1=S^1\cup S^1_+$ with $S^1_+=\pd U^1\setminus S^1$.
Suppose that $g^1(\cdot,\cdot)$ is a $C^1$ Leray map on $D^{12}\times \ov{U^1}$ and $g^2$ is a $C^1$ Leray map on $D^{12}\times S^2$.
 Let $ f$ be  a $(0,q)$-form  such that $  f$ and $\db f$ are in $C^1(\ov{ D^{12}})$. Then on
 $D^{12}$
  we have
\al{}\label{tsqf}
 f&=L_{1;q}^1f+L_{2;q}^2f+L_{12; q}^{12}f+\db  H_q f+H_{q+1}\db f,\quad \text{if $q\geq1$},\\
\label{tsqf+}
 f&=L_{1;0}^1f+L_{2;0}^2f+L_{12;0}^{12}f+H_1\db  f, \hspace{5.5em} \text{if $q=0$}
\end{align}
 where for  $s=q,q+1$ with $q\geq1$ and $E=E_{D^{12}}$
\al{}
\label{hqf}
 H_s f&:=H^{(1)}_s f+ H_s^{(2)}f,
\\
\label{hq1} H^{(1)}_s f&:=R_{U^1\cup D^{12}; s-1}^0 E f+R_{U^1;s-1 }^{01}[\db,E] f,
 \\
 \label{hq2}
 H^{(2)}_sf&:=-R^1_{U^1;s-1}Ef+L_{1^+;s-1}^{01} Ef +L_{2;s-1}^{02} f+L_{12;s-1}^{012}f,\\
L^{01}_{1^+; s-1}Ef&:=\int_{S^1_+}\Omega^{01}_{0,s-1}\wedge Ef,
\\
H_0 f&:=\int_{\pd D^{12}}\Om_{(0,0)}^1 f-\int_{U^1}\Om_{(0,0)}^1\wedge E\db  f=\int_{U^1}\Om_{(0,0)}^1\wedge [\db, E]  f.
 \label{H0f}
\end{align}
\epr
\begin{rem}Formula \re{tsqf} is called   an {\it approximate  homotopy formula} due to the presence of the boundary integrals $L^1_{1;q},L^{2}_{2;q},L^{12}_{12;q}$. These boundary integrals will be transformed under further conditions on the Levi-form of $\pd D^1$.
\end{rem}
\begin{proof}
Consider first $q\geq1$. We  recall the Bochner--Martinelli--Koppelman formula~\cite{MR1800297}*{p.~273} and a version for domains with piecewise $C^1$ boundary~\cite{MR986248}*{Def.~3.1, p.~46; Thm.~3.12, p.~53}:
 \aln 
  f(z)&=\db_z\int_{D^{12}}\Om_{(0,q-1)}^0(z,\zeta)\wedge f(\zeta)+\int_{D^{12}}\Om_{(0,q)}^0(z,\zeta)
 \wedge\db f(\zeta)\\
 \nonumber &\quad +\int_{\pd{D^{12}}}\Om_{(0,q)}^0(z,\zeta)\wedge f(\zeta).
 \end{align*}
Here and in what follows, we assume $z\in D^{12}$.   To ease notation, we drop the degree indicator in $L^1_{1;q}$ and write it as $L^1_{1}$ and do the same for other operators. 
To apply Stokes' formula on $\pd D^{12}$, we use notation in  \re{kop1}-\re{kop3} and rewrite the last term as
\aln
\int_{\pd {D^{12}}}\Om_{(0,q)}^0(z,\zeta)\wedge f(\zeta)
&=\int_{S^1}\Om_{(0,q)}^0(z,\zeta)\wedge f(\zeta)+\int_{S^2}\Om_{(0,q)}^0(z,\zeta)\wedge f(\zeta).
\end{align*}
Let us transform two boundary integrals via   Koppelman's lemma. By  \re{kop2}, we get
\aln
&\int_{S^1}\Om_{(0,q)}^0(z,\zeta)\wedge f(\zeta)=
L_1^1f(z)
+\int_{S^1}\left(\db_\zeta\Om_{(0,q)}^{01}(z,\zeta)+\db_z\Omega^{01}_{(0,q-1)}\right)\wedge f(\zeta)
\\
&\quad =L_1^1f(z)
+L^{01}_{12}f(z)-\int_{S^1}\Om_{(0,q-1)}^{01}(z,\zeta)\wedge\db_\zeta f(\zeta)-
\db_z\int_{S^1}\Om_{(0,q-1)}^{01}(z,\zeta)\wedge f(\zeta),
\end{align*}
where the last identity is obtained by Stokes' formula for $S^1$ with $\pd S^1=S^{12}$.
Analogously, we get
\al
L_2^0f=L_2^2f+L_{21}^{02}f
-L_{2}^{02}\db f  -
\db L_2^{02} f.
\nonumber
\end{align}
Using $L^{02}_{21}f=-L^{02}_{12}f$, we  get
$$
L^{01}_{12}f(z)+L_{21}^{02}f(z)=-L_{12}^{12}f(z)+\int_{S^{12}}\db_z\Om^{012}_{(0,q-1)}(z,\zeta)\wedge f(\zeta)+\int_{S^{12}}\db_\zeta\Om^{012}_{(0,q)}(z,\zeta)\wedge f(\zeta).
$$
Using Stokes' theorem to the last term and  $\pd(S^1\cap S^2)=\emptyset$, we obtain from \re{checksign+}
$$
L^{01}_{12}f+L_{21}^{02}f=-L_{12}^{12}f+\db L_{12}^{012}f
+L_{12}^{012}\db f.
$$
This shows that
\aln 
  f(z)&=-
\db_z\int_{S^1}\Om_{(0,q-1)}^{01}(z,\zeta)\wedge f(\zeta)+\db_z\int_{D^{12}}\Om_{(0,q-1)}^0(z,\zeta)\wedge f(\zeta)\\ \nonumber&\quad-\int_{S^1}\Om_{(0,q)}^{01}(z,\zeta)\wedge\db  f(\zeta)
+\int_{D^{12}}\Om_{(0,q)}^0(z,\zeta)
 \wedge\db f(\zeta)\\
\nonumber &\quad-L_{12}^{12}f(z) + \db L_{12}^{012}f(z)+L_{12}^{012}\db f(z)
  -L_2^{02}\db f(z) -\db L_2^{02}f(z).
 \end{align*}
Next, we transform both   integrals on $S^1$ into volume integrals using  Stokes' formula. Here we need to modify the methods in Lieb--Range~\cite{MR597825} and~\cite{MR3961327}, since $S^1$ has boundary.

For the rest of the proof let $E$ be the Stein extension operator $E_{D^{12}}$.

With orientations, we have $\pd U^1= S^1_+-S^1$.
 By Stokes' formula and
 \cite{MR3961327}*{(2.12)},  we have 
\aln 
&-\int_{\zeta\in S^1}\Om_{(0,q-1)}^{01}( z,\zeta )\wedge  f(\zeta)+ L^{01}_{1^+}E f(z)=\\
 &
 \nonumber\qquad \qquad\int_{U^1 }
\Om_{(0,q-1)}^{01}(z,\zeta)\wedge\db E f(\zeta) +
\int_{U^1 }
\Om_{(0,q-1)}^{0}(z,\zeta)\wedge E f(\zeta)\\
&\nonumber\qquad
\qquad-\int_{U^1 }
\Om_{(0,q-1)}^{1}(z,\zeta)\wedge E f(\zeta)+
\int_{U^1 }
\db_z\Om_{(0,q-2)}^{01}(z,\zeta)\wedge E f(\zeta).
\end{align*}
This shows that on $D^{12}$
\al\label{keyidsim}
R^0_{D^{12}}f-L_{1}^{01} f&=-  L^{01}_{1^+}E f + R_{U^1 }
^{01}\db E f\\
&\quad +
R_{U^1\cup D^{12} }
^{0} E f-R_{U^1 }^
{1} E f+\db
R_{U^1 }
^{01} E f.\nonumber
\end{align}
After applying $\db$, the last term will be dropped.
This shows that
\al\label{keyidsim+}
\db R^0_{D^{12}}f- \db L_1^{01}f &=-\db  L^{01}_{1^+}Ef-
\db
R_{U^1}^{1} E f\\
&\quad+\db  R_{U^1 }
^{01}\db E f+
\db R_{U^1\cup D^{12}  }
^{0} E f.
\nonumber\end{align}
We apply \re{keyidsim}  in which $f$ is replaced by $\db f=E\db f$ to obtain
\aln
R_{D^{12}}^0\db f-L_{1}^{01}\db f&=-L_{1+}^{01}E\db f
 +R_{U^1 }
^{01}\db E\db f +
R_{U^1\cup D^{12} }
^{0} E\db f\\
&\nonumber \quad-R_{U^1 }
^{1} E\db f  +
\db R_{U^1 }
^{01} E\db f.
\end{align*}
We can pair the last term with the second last term in \re{keyidsim+} to form the desired commutator $[E,\db] f$. Finally, we have $\db E\db f=[\db, E]\db f$.
This completes the proof of \eqref{tsqf}. The above proof is still valid for \eqref{tsqf+} (case $q=0$), as the Koppelman lemma holds with $\Om_{0,-1}=0$.

Strictly speaking, the above computation is only valid when $\pd D\in C^3$, since the Koppelman lemma can be verified easily when all Leray maps $W_j\in C^2$. When $\pd D^i\in C^2$, one can still verify the integral formula on the domain $D^1\cap D^2$ by smoothing $g^j$. For instance, see \cite{MR4289246}*{p.~6808} for details.
\end{proof}

\setcounter{equation}{0}
\section{A local homotopy formula for
 $(n-q)$ convex configuration}

The main purpose of this section is to construct a local homotopy formula near a boundary point of strictly $(n-q)$ convex.

In \rp{hf}, we have derived a  local approximate homotopy formula \re{tsqf}:
$$ f=L_{1;q}^1f+L_{2;q}^2f+L_{12;q}^{12}f+\db  H_q f+H_{q+1}\db f.
 $$
 To obtain a genuine local homotopy formula, we will show that the boundary integrals   $L^1_{1;q}f, L^2_{2;q}f, L^{12}_{12;q}f$
  vanish when the boundary is   $(n-q)$ convex and the Leray mappings $g^1,g^2$
are chosen appropriately.

The constructions in this section and the next are inspired by  Henkin--Leiterer~\cite{MR986248}.

We first transform a   $(n-q)$-convex domain $D$ into a new form $D^1$.

\le{convex-rho}Let $D\subset  U_0$ be defined by $\rho^0<0$ with $\rho^0\in C^2(U_0)$.
Suppose that $\pd D$ is $(n-q)$-convex at $\zeta\in\pd D$  and $ \nabla\rho^0(\zeta)\neq0$.
\bpp
\item There are an open set $U_1\Subset U_0$ containing $\zeta$ and  a  biholomorphic mapping $\psi$ defined on a neighborhood of $\ov{U_1}$  such that $\psi(\zeta)=0$, $U:=\psi(U_1)$ is a polydisc, and  $D^1:=\psi(U_1\cap D)$ is defined by
\eq{qconv-nf}
\rho^1(z)=-y_{n}+\la_1|z_1|^2+\cdots+\la_{q-1}|z_{q-1}|^2+|z_{q}|^2+
\cdots+|z_{n}|^2+R(z)<0,
\end{equation}
where $|\la_j|<1/4$ and $R(z)=o(|z|^2)$.  There exists $r_1=r_1(\nabla\rho^0,\nabla^2\rho^0)>0$  such that the boundary  $\pd D^1$ intersects the sphere $\pd B_r$ transversally when $0<r<r_1$. Furthermore, the function $R$ in \rea{qconv-nf} is in $C^a(B_{r_1})$ $($resp. $\Lambda^a(B_{r_1}))$,
 when $\rho^0\in C^a(U_0)$ with $a\geq2$ $($resp. $\Lambda^a(U_0)$ with $a>2)$.
\item Let $\psi$ be as above. There exists $\delta(D)>0$ such that if $\tilde D$ is defined by $\tilde\rho^0<0$ and $\|\tilde\rho^0-\rho^0\|_{C^2(U_0)}<\delta(D)$, then  $\psi(U_1\cap\tilde D)$ is given by
\eq{qconv-nf-t}
\tilde\rho^1(z)=-y_{n}+\la_1|z_1|^2+\cdots+\la_{q-1}|z_{q-1}|^2+|z_{q}|^2+
\cdots+|z_{n}|^2+\tilde R(z)<0
\end{equation}
with  $\|\tilde R-R\|_{\Lambda^a(B_{r_1})}\leq C_a\|\tilde\rho^0-\rho^0\|_{\Lambda^a(U_0)}$ for $a>2$
and  $\|\tilde R-R\|_{C^a(B_{r_1})}\leq C_a\|\tilde\rho^0-\rho^0\|_{C^a(U_0)}$ for $a\geq2$.
 There exists $r_1>0$ such that  the boundary  $\pd(\psi(U_1\cap \tilde D))$ intersects the sphere $\pd B_{r_2}$ transversally when $r_1/2<r_2<r_1$.
\epp
Here $\delta(D)$ depends on the modulus of continuity of $\nabla^2\rho^0$.
\ele
\begin{rem}\bpp\item  Throughout the paper, we denote by $C(\nabla\rho^0,\nabla^2\rho^0)$, such as the above $r_1(\nabla\rho^0,\nabla^2\rho^0)$, a constant depending on $\nabla\rho^0,\nabla^2\rho^0$ in local coordinates.  In particular, when $\rho_0$ is a defining function of an $a_q$ domain $D$, the constant $C(\nabla\rho^0,\nabla^2\rho^0)$ depends on the Levi-form $L_\zeta \rho^0$ and the low and upper bounds of $|\nabla_\zeta\rho^0|$  for $\zeta\in \pd D$.
\item
We will   refer to $(\tilde D^1,\tilde\rho^1)$ as $(D^1, \rho^1)$ indicating the various constants in estimates are {\it uniform} in  $\tilde\rho^0$
 or {\it stable} under small $C^2$ perturbations of $\pd D$ via $\rho^0$.  Set
 $$D^2\colon \rho^2(z)=|z|^2-r_2^2<0$$  with restriction $r_1>r_2>r_{1}/2$. We will  assume that   $D^2$ is contained in the polydisc $ U$ in \rl{convex-rho} $(a)$.
 \epp
\end{rem}
\begin{proof}
$(a)$      We may assume that $\zeta=0$. Permuting coordinates yields $\rho^0_{z_n}\neq0$. Let $\tilde z_n=2\rho^0_\zeta\cdot(\zeta-z)-\sum \rho^0_{\zeta_j\zeta_k}(\zeta_j-z_j)(z_k-\zeta_k)$ and $\tilde z'=z'$. Then $\rho_1(\tilde z):=\rho^0(z)$, where the new domain has a defining function
$$
\rho_1(z)=-y_n+\sum a_{j\ov k}z_j\ov z_k+o(|z|^2).
$$
 Choose a nonsingular matrix $A$ for a linear change of coordinates
 \eq{Utz}
 A\tilde z=z \quad \text{ with $\tilde z_n=z_n$}.
  \end{equation}
 Set $\rho_2(\tilde z):=\rho_1(z)$. The new domain has the definition function
  $$\rho_2(z)=-y_n+\sum_{j<q}\la_j|z_j|^2+\sum_{j=q}^{n-1}|z_{j}|^2 +\sum_{j=1}^{n}\RE\{a_jz_j\ov z_n\}+o(|z|^2),
$$
where $a_1,\dots, a_{n-1}\in\cc$ and $a_n\in\rr$.  Set $\rho_3:=\rho_2+\mu\rho_2^2$ with $\mu=a_{n}+1$. We get
$$
\rho_3(z)=-y_n+\sum_{j<q}\la_j|z_j|^2+\sum_{j=q}^{n}|z_{j}|^2+\RE\Bigl\{\mu z_n^2+\sum_{j=1}^{n-1}a_jz_j\ov z_n\Bigr\}+o(|z|^2).
$$
Using new coordinates $\tilde z_n=z_n-i\mu z_n^2+iz_n\sum_{j=1}^{n-1}a_jz_j$ and $\tilde z'=z'$,
we get the defining function of a new domain:
$$
\rho_4(z)=-y_n+\la_1|z_1|^2+\cdots+\la_{q-1}|z_{q-1}|^2+|z_{q}|^2+
\cdots+|z_{n}|^2+y_n\RE\{\sum_{j=1}^{n}b_jz_j\}+o(|z|^2).
$$
 We get $|\la_j|<1/4$ after a   dilation. Then $\rho_4+\rho_4\RE\{\sum b_jz_j\}$, renamed as $\rho^1$,  has the form \re{qconv-nf}. As usual, the implicit function theorem can be proved by using the inverse mapping theorem~\cite{MR0385023}*{pp.~224-225}. Then the inverse mapping theorem for Zygmund spaces~\cite{GG} yields the desired smoothness of $R$.
 The details are left to the reader. The transversality of $\pd D^1$ and $\pd D^2$ also follows from the computation below.

$(b)$ The above construction of $\psi$ is explicit in $\rho^0$ with the exception of the linear change of coordinates $
\tilde z= A^{-1}z$ in \re{Utz} that is fixed for all small perturbations $\tilde\rho^0$ of $\rho^0$. Thus, it is easily to check that  $\|\tilde R-R\|_{C^a(R_{r_1})}\leq C_a\|\tilde\rho^0-\rho^0\|_{C^a(U_0)}$ for $a\geq2$ and $\|\tilde R-R\|_{\Lambda^a(B_{r_1})}\leq C_a\|\tilde\rho^0-\rho^0\|_{\Lambda^a(U_0)}$ for $a>2$.

 We want to show that $\nabla \tilde\rho^1$ is not proportional to $\nabla\rho^2$ on the common zero set of $\tilde\rho^1,\rho^2$. Suppose that
$\nabla\rho^2=\mu\nabla\tilde\rho^1$
when $\tilde\rho^1(z)=\rho^2(z)=0$. We get $2y_n=\mu (-1+2 y_n+\tilde R_{y_n})$. When $r<1/4$ and $\delta(D)$ are sufficient small, by $|z|=r$ we obtain $|-1+2 y_n+\tilde R_{y_n}|<1/2$. Hence $-\mu^{-1}y_n\in(1/4,3/4)$ as
$$
\|\tilde R\|_{C^2(B_{r_1})}<1/C.
$$  For $j<n$, we have
$2y_j=2\mu\tilde\rho_{y_j}$. This shows that $|y_j|\leq C|y_n|$. Also, $|x_k|\leq C|y_n|$. Thus $\tilde\rho^1(z)=0$ implies $|y_n|\leq C'|y_n|^2+|\tilde R(z)|$. In view of $C'|y_n|<1/2$, we get
$$
|z|\leq C\|R\|_{C^2(B_{r_1})}|z|^2+C\del(D).
$$
By choosing $\del(D)$ depending on $r_1$, we get $|z|<r_1^2/C$. The latter contradicts the vanishing of $\rho^2(z)=|z|^2-r^2_2$ since $r>r_1/2$.
 \end{proof}

Recall that our original domain $D$ is normalized as $D^1$. We now fix notation. Let $(D^1,U,\phi,\rho^1)$ be as in \rl{convex-rho}. Thus $\rho^1$ is given by \re{qconv-nf} (or \re{qconv-nf-t}). Recall  that
\eq{rho2}
\rho^2(z)=|z|^2-r^2_2
\end{equation}
 where $0<r_2<r_1$ and $r_1/2<r_2<r_1$     for \rl{convex-rho} $(a), (b)$. Let us define
\ga{}\label{d12}
D^1\colon\rho^1<0, \quad D^2\colon\rho^2<0, \quad D^{12}=D^1\cap D^2,\\
\label{s12}
\pd D^{12}=S^1\cup S^2, \quad S^i\subset\pd D^i,\\
\label{g02}
g^2( z,\zeta )=\Bigl(\DD{\rho^2}{\zeta_1},\dots, \DD{\rho^2}{\zeta_n}\Bigr)=\ov\zeta.
\end{gather}
It is well-known that
\eq{W2}
| g^2(z,\zeta)\cdot(\zeta-z)|>0, \quad (z,\zeta)\in D^2\times\pd D^2.
\end{equation}

\le{}Let $(D^1,U,\phi,\rho^1)$ be as in \rla{convex-rho}. Define
\eq{}\label{g1}
g^{1}_j( z,\zeta )=\begin{cases}
\DD{\rho^1}{\zeta_j},&q\leq j\leq n,\\
\DD{\rho^1}{\zeta_j}+(\ov\zeta_j-\ov z_j),& 1\leq j<q.
\end{cases}
\end{equation}
Then for $\zeta,z\in U$ and by shrinking $U$ if necessary, we have
\al{}\label{W1-dist}
2\RE\{ g^1( z,\zeta )\cdot(\zeta-z)\}&\geq \rho^1(\zeta)-\rho^1(z)+\f{1}{2}|\zeta-z|^2, \\
|g^1(z,\zeta)\cdot(\zeta-z)|&\geq 1/C_*, \quad \forall\zeta\in S^{12}, z\in B_{\tilde r_2}\label{g1z}\end{align}
where $\tilde r_2<r_2/{C}$. Here  $C,C_*$ depend $\nabla\rho^1,\nabla^2\rho^1$ $($and hence only on $\nabla\rho^0,\nabla^2\rho^0)$,  and $r_2$; $C_*$ is independent of $\tilde r_2$.
\ele
\begin{proof} We have
$
\RE\{ g^1( z,\zeta )\cdot(\zeta-z)\}=\RE\{ \rho^1( z,\zeta )\cdot(\zeta-z)\}+\sum_{j<q}|\zeta_j-z_j|^2.
$
 By Taylor's theorem, we have
\al\label{taylor}
\rho^1(z)-\rho^1(\zeta)&=2\RE\{\rho^1_{\zeta}\cdot(z-\zeta)\} +\RE\{\rho^1_{\zeta_j\zeta_k}(z_j-\zeta_j)(z_k-\zeta_k)\}\\
&\quad+
H_\zeta\rho^1(z-\zeta)+R( z,\zeta ).\nonumber
\end{align}
Note that $H_\zeta\rho^1$, restricted to the $(z_q,\cdots, z_n)$ subspace, is a positive definite quadratic form. Also, $\rho^1_{\zeta_j\zeta_k}$ and second-order derivatives of $R$ are small.  We can show that
\aln\rho^1(z)-\rho^1(\zeta)
&\geq2\RE\{\rho^1_{\zeta}\cdot(z-\zeta)\}+\sum_{j\geq q}|\zeta_j-z_j|^2-c|\zeta-z|^2
\end{align*}
where $c<1/2$. Now \re{W1-dist} follows  from \re{qconv-nf}. Note that \re{g1z} follows from \re{W1-dist} as $\rho^1(\zeta)=0$, $|\zeta-z|^2\geq r_2^2/4$  for $\zeta\in S^{12}\subset (\pd D^{1})\cap\pd D^2_{r_2}$, and $|\rho^1(z)|<C\tilde r_2$  for $z\in B_{\tilde r_2}$.  By \re{taylor} and \re{W1-dist} for case \re{qconv-nf},  we also get \re{W1-dist} for the  $\tilde\rho^1$ in
\re{qconv-nf-t} when $\|\tilde\rho^0-\rho^0\|_2\leq C\del(D)$ is  small.
\end{proof}

\begin{figure}[!ht!]\counterwithin{figure}{section}
\centering
\begin{tikzpicture}[scale=.75,even odd rule]
 \begin{scope}
\clip  (0,-3)  ellipse (4 and 3); 
 \fill[gray!40]
  (0,0)  circle (2.5);  
 \end{scope}

  \draw[thick,domain=(1/8)*pi:(7/8)*pi,color=black, name path=D1]
 plot[smooth]({0+4*cos (\x r)},{-3
 +3*sin(\x r)});

  \draw[dashed, thick] (0,0) circle (2.5); 

  \fill (0,0)  circle (.125);   \node at (0,.4) {$0$};

\node at (0,-1) {$D^{12}$};
   \node at (0,2) {$D^{2}$};
    \node at (3.0,-1.5) {$D^{1}$};
\end{tikzpicture}
 \vspace{-9ex}
    \caption{Convex configuration: $D^{12}=D^1\cap D^2$}
\vspace{3ex}
 \label{fig:convex}
\end{figure}

\begin{defn}\label{def:3.4}
\bpp \item As in~\cite{MR986248}*{p.~80}, the $(U,D^1, D^2,\psi,\rho^1, \rho^2)$  in \rl{convex-rho}~$(a)$ and
\re{rho2}-\re{d12} is called an {\it  $(n-q)$-convex configuration} (of type I).  \rl{convex-rho} $(b)$ is referred to as the {\it stability}  of the configuration. In brevity,  we call $(D^1,D^2)$ an $(n-q)$-convex configuration. See \rf{fig:convex}, i.e. the figure for type I in~\cite{MR986248}*{p.~80}.
\item The $g^1,g^2$ given by \re{g02} and \re{g1} are called the {\it canonical Leray maps} for the $(n-q)$-convex configuration $(D^1,D^2)$.
    \epp
\end{defn}
Note that $g^2(z,\zeta)$ is holomorphic in $z$ and  $g^1(z,\zeta)$ is anti-holomorphic in merely $q-1$ variables of $z$.  Checking the types, we have
\al{}
\label{type-2}
 &\Omega_{(0,k)}^{2}( z,\zeta )=0, \  k\geq1;\quad \db_z\Om^2_{(0,0)}(z,\zeta)=0;\\
\label{type-1}
& \Omega_{(0,k)}^{1}( z,\zeta )=0, \  k\geq q;\quad \db_z\Om_{(0,q-1)}^1( z,\zeta )=0;
\\
\label{type-12}
&\Omega_{(0,k)}^{12}( z,\zeta )=0,\  k\geq q;\quad \db_z\Om_{(0,q-1)}^{12}( z,\zeta )=0.
\end{align}
By  \re{type-2}-\re{type-12}, we have for $z\in D^{12}$
\ga\label{L110}
L_{i;q}^if(z)=\int_{S^i}\Om_{(0,q)}^i(z,\zeta)\wedge f(\zeta)=0, \  i=1,2; \quad L_{12;q}^{12}f(z)=0;\\
 \db_z\int_{U_1}\Om_{(0,q-1)}^1(z,\zeta)\wedge Ef(\zeta)=0, \quad \int_{U_1}\Om_{(0,q)}^1(z,\zeta)\wedge E\db f(\zeta)=0.
\end{gather}

This shows that for $s=q,q+1$,  the $H^{(2)}_s$ in \re{hq2} can be replaced by
\eq{nhq2}
H^{(2)}_sf=L_{1^+;{s-1}
}^{01} Ef +L_{2;s-1}^{02} f+L_{12;s-1}^{012}f.
\end{equation}
Therefore, we have obtained the following local homotopy formula.
\begin{thm}\label{hf-c} Let $ 0<q\leq n$.
Let $(D^1,D^2)$ be a $(n-q)$-convex configuration with   Leray maps  \rea{g02}-\rea{g1}.
Let $U^1=D^2\setminus\ov {D^1}$ and  $S_+^1=\pd D^2\setminus D^1$.
 Suppose that $  f$ is a $(0,q)$ form such that $f$ and $\db f$ are in $C^1(\ov {D^{12}})$. Then on $D^{12}$
\gan
 f= \db  H_q f+  H_{q+1}\db f
\end{gather*}
with  $H_s=H_s^{(1)}+H_s^{(2)}$, for $H^{(1)}_s$   defined by \rea{hq1} and $H_s^{(2)}$   defined by \rea{nhq2}.
\end{thm}
We remark that the   kernel  in $H^{(2)}_s$ is smooth, since $S^1_{+}, S^2$ and $ S^{12}$ do not intersect small neighborhoods of the origin in $\pd D^1$. Therefore, terms in $H^{(2)}_s$ can be estimated easily, while the main term $H_s^{(1)}$     will be estimated in Section~\ref{sec:.5est}.

\setcounter{equation}{0}
\section{A local $\db$ solution operator for
 $(q+1)$-concave configuration}

 We recall   from  \rp{hf} the approximate local homotopy formula
\eq{f=l11} f=L_{1;q}^1f+L_{2;q}^2f+L_{12;q}^{12}f+\db  H_q f+H_{q+1}\db f.
 \end{equation}
 As the   strictly $(n-q)$-convex case, we will show that  $L^1_{1;q}f, L^2_{2;q}f$ vanish when the boundary is   $(q+1)$-concave and the Leray mappings $g_1,g_2$
are chosen appropriately. However, the boundary integral $L^{12}_{12;q}f$ may not vanish. We will show that this term is $\db$-closed for a $\db$-closed $f$, and this allows us to use a third Leray mapping to transform it into a genuine $\db$ solution operator
$f=\db H_qf$ for possibly different $H_q$.
Thus, the $(q+1)$ concavity is sufficiently to construct a $\db$ solution operator.

The presence of $L_{12;q}^{12}f$ will lead to a subtlety.
For a local homotopy formula for forms which are not necessarily $\db$-closed, we  however  need an {\it extra} negative Levi eigenvalue, which will be assumed at the end of the section.

The following is a restatement of \rl{convex-rho}, by considering the complement $(D^1)^c$ and $-\rho_0$ where $\rho_0$ defines $D^1$.
 \le{concave-rho}
  Let $D\subset  U_0$ be defined by $\rho^0<0$ with $\rho^0\in C^2(U_0)$.
Suppose that $\pd D$ is $(q+1)$-concave at $\zeta\in\pd D$ and  $ \nabla\rho^0(\zeta)\neq0$.
Then all assertions, including dependence of various
 constants on $\nabla\rho^0,\nabla^2\rho^0$,  in \rla{convex-rho} on $\rho^1,\psi$ are valid, provided \rea{qconv-nf}-\rea{qconv-nf-t} are replaced by
\begin{gather}\label{rho1-v}
\rho^1(z)=-y_{q+2}-|z_1|^2-\cdots-|z_{q+2}|^2+\la_{q+3}|z_{q+3}|^2+
\cdots+\la_n|z_n|^2+R(z),\\
\label{rho1-t-v}
\tilde\rho^1(z)=-y_{q+2}-|z_1|^2-\cdots-|z_{q+2}|^2+\la_{q+3}|z_{q+3}|^2+
\cdots+\la_n|z_n|^2+\tilde R(z)
\end{gather}
with $|\la_j|<1/4$ for $j>q+2$.
\ele
As in \rl{convex-rho}, we assume $\psi(U_1)$ is a polydisc $U$.

When $\rho^1$ has the form \re{rho1-v}, as in~\cite{MR986248}*{pp.~118-120} define
\eq{}\label{HL2pg81}
g^1_{j}( z,\zeta )=\begin{cases}
\DD{\rho^1}{z_j},&1\leq j\leq q+2,\\
\DD{\rho^1}{z_j}+\ov z_j-\ov \zeta_j,& q+3\leq j\leq n.
\end{cases}
\end{equation}
Note that this kind of Leray mappings was   used by Hortmann~\cites{MR422688, MR627759} for strictly concave domains.
Then we have
\eq{W1s-dist}
-2\RE\{g^1( z,\zeta )\cdot(\zeta-z)\}\geq \rho(\zeta)-\rho(z)+|\zeta-z|^2/C.
\end{equation}
An essential  difference between the $q$-convex and $(q+1)$ concave cases is that $g^1( z,\zeta )$ is no longer $C^\infty$ in $z$ when $\pd D$ is only finitely smooth. Nevertheless, a useful feature is that $g^1(z,\zeta)$ is holomorphic in $\zeta_1,\dots, \zeta_{q+2}$.

As $(n-q)$-convex case, we   use
$
D^2_{r_2}\colon
\rho^2(z):=|z|^2-r_2^2<0.
$
We still take Leray maps $g^0(z,\zeta)=\ov\zeta-\ov z$ and
$
g^2(z,\zeta)=(\f{\pd\rho^2}{\pd\zeta_1},\dots, \f{\pd\rho^2}{\pd\zeta_n})=\ov \zeta.
$
Denote by $\text{deg}_\zeta$  the degree of a form in $\zeta$. We get
\ga
\label{Ca-type-1}
\text{deg}_\zeta\,  \Omega_{(0,*)}^{1}(z,\zeta)\leq 2n-q-2.
\end{gather}
Therefore, we still have \re{L110}. Thus \re{f=l11} becomes
\eq{f=l11+} f= L_{12; q}^{12}f+\db  H_q f+H_{q+1}\db f.
 \end{equation}

However, unlike the $(n-q)$ convex case,
$
L_{12;q}^{12}f=\int_{S^{12}} \Om^{12}_{(0,q)}\wedge f
$
may not be identically zero. Let us  transform this integral on $S^{12}$ via  Stokes' formula.  We intersect $D^1\cap D^2$ with a third domain
\eq{defD3}
D^3\colon \rho^3<0, \quad 0\in D^3
\end{equation}
where
\eq{rho3}
\rho^3(z):=-y_{q+2}+\sum_{j=q+3}^n3|z_j|^2-r^2_3
\end{equation}
with $0<r_3<r_2/{C_n}$; see \rl{y123} below for further restrictions.
As in~\cite{MR986248}*{p.~120}   define
\ga
\label{defW-3}
g^{3}_j( z,\zeta )=\begin{cases}
0,&1\leq j<q+2,\\
i,&j= q+2,\\
3(\ov\zeta_j+\ov z_j),& q+3\leq j\leq n.
\end{cases}
\end{gather}

We can verify
\eq{HH}
\RE\{g^3(z,\zeta)\cdot(\zeta-z)\}=\rho^3(\zeta)-\rho^3(z).
\end{equation}
Note that $0\in D^{23}$.
See \rf{fig:concave} for relations of $D^1, D^2, D^3$.

 Then
we have the following.
\le{y123} Let $\rho^i, g^i$ be defined by \rea{rho2}, \rea{g02}, \rea{rho1-v}, \rea{HL2pg81}, \rea{rho3}-\rea{defW-3} for $i=1,2,3$.
\bpp
\item
There exists $r_1$, depending on $\nabla\rho^0,\nabla^2\rho^0$, such that  $\pd D^1, \pd D^2,\pd D^3
$  pairwise intersect transversally when $0<C_nr_3<r_2<r_1$ and $r_1<2r_2$.
\item Let $\tilde\rho^0,\tilde\rho^1$ be as in \rla{concave-rho} and $\tilde D^1$ be defined by $\tilde\rho^1<0$.
If $\del(D)$ is sufficiently small, $\|\tilde \rho^0-\rho^0\|_2<\del(D)$, and $1/{C_n'}<C_nr_3<r_2<r_1$ and $r_2<r_1/2$, then $\pd\tilde D^1,\pd D^2,\pd D^3$  pairwise intersect transversally.
 \item Let $r_1,r_2,r_3$ be as in $(b)$. Here $C_n'$ depends only on $n$.   Then
 \ga\label{gizC}\pd\tilde D^1\cap\pd D^2\cap D^3=\emptyset,\quad
 \pd\tilde D^1\cap\pd D^3\cap\pd D^2=\emptyset,\\
 |g^i(z,\zeta)\cdot(\zeta-z)|\geq 1/C_*, \quad \forall(z,\zeta)\in B_{r_4} \times \ov{D^2\setminus(D^1\cup D^3)}, \ i=0,2,3,
 \label{gicZ+}\\
S^{12}\subset \ov{D^2\setminus(D^1\cup D^3)}.
\label{gicZ=}
 \end{gather}
 \epp
Here  $C_n,C_n'$ depend only on $\nabla\rho^1,\nabla^2\rho^1$ $($and hence only on $\nabla\rho^0,\nabla^2\rho^0)$.
\ele

\begin{figure}
\vspace{-40ex}
\centering
\begin{tikzpicture}[scale=.75, even odd rule]
 \begin{scope}
\clip (0,5.5) ellipse (8 and 6);
 \fill[gray!40] (0,-3)  ellipse (4 and 3) (0,0)  circle (3);
 \end{scope}

 \draw[dashed,thick, domain=(0.5/8)*pi:(7.5/8)*pi,color=black, name path=D1]
 plot[smooth]({0+4*cos (\x r)},{-3
 +3*sin(\x r)});

  \draw[dashed, thick] (0,0) circle (3);

   \draw[dashed,thick, domain=(10.65/8)*pi:(13.45/8)*pi, name path=D3]
 plot[smooth]({0+8*cos (\x r)},{5.5+6*sin(\x r)});

 \begin{scope}
\clip (0,-3) ellipse (4 and 3);
 \fill[gray!40,]     (0,6)  ellipse (0 and 6) (0,0)  circle (3);
 \end{scope}

 \draw[dashed, thick,domain=(2.8/8)*pi:(5.2/8)*pi,color=black, name path=D1p]
 plot[smooth]({0+4*cos (\x r)},{-3
 +3*sin(\x r)});
   \draw[domain=(11.5/8)*pi:(12.5/8)*pi, name path=D3p]
 plot[smooth]({0+8*cos (\x r)},{5.5+6*sin(\x r)});

 \tikzfillbetween[of=D1p and D3p]{white};

\fill (0,0)  circle (.125);

\draw[thick, domain=1.1*pi:1.9*pi,] plot[smooth] ({3*cos(\x r)},{3*sin(\x r)});

\draw[thick,domain=(1.95/8)*pi:(2.75/8)*pi,color=black, name path=D1]
 plot[smooth]({0+4*cos (\x r)},{-3 +3*sin(\x r)});
 \draw[ thick,domain=(5.27/8)*pi:(6.05/8)*pi,color=black, name path=D1]
 plot[smooth]({0+4*cos (\x r)},{-3
 +3*sin(\x r)});
\draw[thick, domain=(11.4/8)*pi:(12.6/8)*pi, name path=D3]
 plot[smooth]({0+8*cos (\x r)},{5.5+6*sin(\x r)});

   \node at (0,.4) {$0$};
\node at (0,1.65) {$D^{123}$};
   \node at (0,-1.5) {$\overline{D^{2}}\setminus(D^1\cup D^3)$};
    \node at (3.75,-1.125) {$D^{1}$};
     \node at (3.75,.5) {$D^{3}$};
\end{tikzpicture}
 \vspace{-8ex}
    \caption{Concave configuration: $D^{123}=D^1\cap D^1\cap D^3$}
      \vspace{3ex}
   \label{fig:concave}
\end{figure}

\begin{proof}
$(a)$ Suppose $\nabla\rho^1(z)=\mu\nabla\rho^3(z)$, and $\rho^1(z)=\rho^3(z)=0$. We have $-1+2y_{q+2}+R_{y_{q+2}}=\mu$. This shows that $-3/2<\mu<-1/2$. We also have $z_j+o(|z|)=0$ for $j=1,\dots q+1$,
$x_{q+2}+o(|z|)=0$,  and
$\la_jz_j+o(|z|)=3\mu z_j$
for $j>q+2$. The latter implies that $|z_j|=o(|z|)$ since $|3\mu|-|\la_j|>1/4$. Hence, $\rho^1(z)=0$ yields $y_{q+2}=o(|z|)$. This shows that $z=0$, which contradicts $\rho^3(0)<0$.

To show   $\pd D^2$ and $\pd D^3$ intersect transversally, suppose that at an intersection point $z$ we have $\nabla\rho^2=\mu\nabla\rho^3$. We first get $2y_{q+2}=-\mu$.    Then $\nabla\rho^2=\mu\nabla\rho^3$ implies that $|z|\leq C_n|y_{q+2}|$. We get $|y_{q+2}|\geq r_2/{C_n}$ and $|y_{q+2}|\leq r_3+|z|^2\leq r_3+C_n^2|y_{q+2}|^2$.  Then $2r_3\geq|y_{q+2}|\geq r_2/{C_n}$, a contradiction to the assumption $r_3<r_2/C$.
  Note that $\pd D^1,
 \pd D^2$
 intersect transversally was proved in \rl{convex-rho}.

 $(b)$ We leave the details to the reader.

 $(c)$ Suppose $\tilde\rho^1(z)=0$ and $\rho^3(z)<0$. Then we have
 $$
 |z_1|^2+\cdots+|z_{q+2}|^2+\sum_{j>q+2}(\la_j+3)|z_j|^2-\tilde R(z)<r_3^2.
 $$
This shows that $|z|<2r_3$.  We obtain the first identity in \re{gizC}. The proof of the second identity in \re{gizC} is similar.

We now verify  \re{gicZ+}. Cases for $i=0,2$ are trivial.  Case $i=3$ follows from \re{defD3}.
Finally, \re{gicZ=} follows from \re{gizC}.
\end{proof}

\begin{defn}\label{ccav}\bpp \item As in~\cite{MR986248}*{p.~119}, the $(U,D^1, D^2,D^3,\psi,\rho^1,\rho^2,\rho^3)$  in \rla{y123}
is called a {\it $(q+1)$-concave configuration},
while the stability in \rl{convex-rho} for the corresponding $(q+1)$-concave case will be called   as the {\it stability}  of the configuration. In short, $(D^1,D^2,D^3)$, shown in  \rf{fig:concave} or figure for type I in~\cite{MR986248}*{p.~119}, is called a {\it $(q+1)$-concave configuration}.
\item The $g^1,g^2,g^3$ in \rea{HL2pg81}, \rea{g02} and \rea{defW-3}  are called the {\it canonical Leray maps} of the configuration.
    \epp
\end{defn}

Note that the anti-holomorphic differentials $d\ov\zeta_j,d\ov z_k$  appear in $\Om^3$ as a wedge product in some of
$$
d(\ov\zeta_j-\ov z_j), \quad q+3\leq j\leq n.
$$
Consequently,  $d\ov z_j$ and $d\ov\zeta_j$, having the {\it same} index,  cannot appear in $\Om_{(0,q)}^3$ simultaneously. Therefore
\ga\label{type-3}
\text{deg}_\zeta\, \Om^3_{(0,\ell)}( z,\zeta )\leq n+([n-(q+3)+1]-\ell)= 2n-q-\ell-2, \quad\forall\ell.
\end{gather}

We now derive a result analogous to \cite{MR986248}*{Lem. 13.6 $(iii)$, p.~122} but for different boundary integrals $L^{\bigcdot}_{12}$.

\begin{lemma} Let $0<q\leq n-2$. Let $(D^1,D^2,D^3)$ be   a  $(q+1)$ concave configuration with Leray maps $(g^1,g^2,g^3)$. Let $U^1=D^2\setminus\ov {D^1}$ and  $S_+^1=\pd D^2\setminus D^1$. Then
\begin{gather}
\label{db13=0}
\Om^{13}_{(0,\ell)}( z,\zeta )=0,\quad\ell<q;\quad \db_\zeta\Om^{13}_{(0,q)}( z,\zeta )=0.
\end{gather}
\bpp\item
Suppose that $f \in C^1_{(0,q)}(\ov{ D^{1}})$ is $\db$-closed   on $D^{1}$. Then on   $D^{123}$,
\begin{align}
\label{dL2312}L_{12;q}^{13}f&=0,\\
L_{12;q}^{12}f&=-L_{12;q}^{23}f+\db L_{12;q-1}^{123}f+L_{12;q}^{123}\db f.
\label{L1212}
\end{align}
\item If $(D^1,D^2,D^3)$ is a $(q+2)$-concave configuration and $0<q\leq n-3$, then \rea{dL2312} and \rea{L1212} are   valid on   $D^{123}$
for any   $f\in C_{(0,q)}^1(\ov{ D^{1}})$.
\epp
\ele
\begin{proof}To verify
\re{db13=0}, we note that $\Om^{13}_{(0,q)}( z,\zeta )$ has type $(0,q)$ in $z$. It has type $(n,n-2-q)$ in $\zeta$ and it is holomorphic in $\zeta_1,\dots, \zeta_{q+2}$. After taking $\db_\zeta$, it has type $(n,n-q-1)$ in $\zeta$ for anti-holomorphic variables $\zeta_{q+3}, \dots, \zeta_n$. However, the number of these anti-holomorphic variables is $<n-(q+3)+2=n-q-1$. We have verified \re{db13=0}.

$(a)$
To verify \re{dL2312}, we need an approximation theorem of Henkin--Leiterer~\cite{MR986248}*{Lemma 13.5 (iii), p. 122} in which $r=n-q-2$ for our $(q+1)$ concave domain $D^1$.

 Fix $z\in D^{123}$. Let  $K=\ov{D^2\setminus(D^1\cup D^3)}$. See \rf{fig:concave} for $K$ and $D^{123}$.
 Note that $D^{123}$ is open, $K$   is compact, and they are disjoint.

 By \re{W1s-dist}, $g^1(z,\zeta)\cdot(\zeta-z)\neq0$ for   $\zeta\in \ov{D^2}\setminus D^1$.  By \re{gicZ+}, we know that $\Om^{13}$ is a continuous $(n,n-q-2)$ in $\zeta\in K$.  Consequently, we can find a sequence $\om^\nu_z$ of   continuous $(n,n-q-2)$-forms in $\zeta\in U\supset \ov{D^2}$ such that  $\om^\nu_z$ are $\db$-closed   on $U$ and  converge to $\Om^{13}(z,\cdot)$ uniformly on $K$ as $\nu\to\infty$. Using a standard smoothing, we may assume that $\om^\nu_z$ are $C^1$ in $\zeta\in \ov{D^2}$, $\db$-closed and approximate $\Om^{13}(z,\cdot)$ uniformly on $K$.
   By \re{gicZ=}, $S^{12}$ is contained in $K$.  We obtain for each fixed $z\in D^{123}$
$$
L_{12;q}^{13}f(z)=\lim_{\nu\to\infty}\int_{S^{12}}\om^\nu_z(\zeta)\wedge f(\zeta).
$$
Since
   $S^1\subset \pd D^{12}\subset U$ and $\db_\zeta \om^\nu_z=0$ on $U$,  Stokes' formula implies
   $$
   \int_{S^{12}}\om^\nu_z(\zeta)\wedge f(\zeta)= \int_{S^1}\om^\nu_z(\zeta )\wedge \db_\zeta f=0.
   $$
Hence $L_{12;q}^{13}f(z)=0$. By
$
\Omega_{(0,q)}^{12}=\Omega_{(0,q)}^{13}-\Omega_{(0,q)}^{23}+\db_\zeta \Omega_{(0,q)}^{123}+\db_z \Omega_{(0,q-1)}^{123},
$
we obtain \re{L1212} through integration on $S^{12}$.

$(b)$
Note that when $(D^1,D^2,D^3)$ is a $(q+2)$ concave configuration. We have $L^{13}_{12;q}f=0$ for any $(0,q)$ forms $f$ that are not necessarily $\db$-closed by \re{db13=0} that now holds for $\ell<q+1$. We get the desired conclusion immediately.
 \end{proof}
\begin{remark}\label{no-control}
In the proof for case $(i)$, we do not have any control on $\om^\nu$ outside $\ov {D^2\setminus(D^1\cup D^3)}$ other than the uniform convergence. Therefore, in $(a)$ it is crucial that   $f$ is $\db$-closed.\end{remark}

So far, we have been following Henkin--Leiterer~\cite{MR986248}. We could derive a homotopy formula on $D^1\cap D^2\cap D^3$ as in~\cite{MR986248}. However, since we only need a local homotopy formula near a boundary point,
we now departure from the approach in~\cite{MR986248}. Let us still use the approximate homotopy formula on $D^{12}$. Modify it only for $z\in D^1\cap D^2\cap B_{r_4}$ using mainly \re{gizC}-\re{gicZ=} to construct a $\db$ solution operator on this smaller domain. To this end, we fix $0<r_4<r_3/{c_n}$ so that
\eq{defr4}
 B_{r_4}\Subset D^2\cap D^3.
\end{equation}
 Thus our starting point is still the approximate homotopy formula \re{f=l11}. We will however use Koppelman's lemma for $g^1,g^2,g^3$ on the set $S^{12}$.

We recall the Koppelman homotopy formula for the ball
$$
g=\db T_{B_{r_4};q}g+T_{B_{r_4};q+1}\db g
$$
for $g\in C^1_{(0,q)}(\ov B_{r_4})$; see~\cite{MR0774049}*{Cor.~1.12.2, p.~60; Cor.~2.1.4, p.~68}.
  We now transform $L_{12;q}^{23}f$ in \re{L1212}.

The following is analogous to \cite{MR986248}*{Lem.~13.7, p.~125} for $L_{23;q}^{23}$.
\begin{lemma} Let $1\leq q\leq n-2$. Let $(D^1,D^2_{r_2},D^3_{r_3})$ be  a  $(q+1)$-concave configuration.
 Let $L_{12,\bigcdot}^{23}$ be defined by \rea{defnLR}. Assume $r_4$ satisfies \rea{defr4}.
For  $f\in C_{(0,q)}^1(\ov{ D^{1}})$, we have for $z\in\ov{B_{r_4}}$
\gan{}\db L^{23}_{12;q}f(z)=\int_{S^{12}}\Om^{23}_{(0,q+1)}(z,\zeta)\wedge \db f(\zeta),\\
L^{23}_{12;q}f=\db T_{B_{r_4};q}L^{23}_{12;q}f+T_{B_{r_4}, q+1}L_{12;q+1}^{23}\db f, \\
L_{12;q}^{23}\colon   C^1_{(0,q)}(
\ov {D^1}
)\cap\ker\dbar\to   C^\infty_{(0,q)}(\ov B_{r_4} )\cap\ker\dbar.
\end{gather*}
\ele
\begin{proof}
 By \re{gicZ+}-\re{gicZ=}, we have $g^i(z,\zeta)\cdot(\zeta-z)\neq0$ for $i=2,3$,  $z\in\ov B_{r_4}$ and $\zeta\in S^{12}$.  Thus, the form $\Om^{23}(z,\zeta)$ is  smooth in
 $z\in B_{r_4}$ and $\zeta\in S^{12}\subset\pd D^2$.

    We have
$$
\db_\zeta\Om_{(0,q+1)}^{23}+\db_z\Om_{(0,q)}^{23}=\Om_{(0,q+1)}^2-\Om_{(0,q+1)}^3.
$$
By \re{type-2}, $\Om_{(0,q+1)}^2=0$. Thus $\int_{S^{12}}\Omega_{(0,q+1)}^2( z,\zeta )\wedge f(\zeta)=0$.
  By  \re{type-3}, the $\zeta$-degree   of $f(\zeta)\wedge \Om^3_{(0,q+1)}( z,\zeta )$ is less than $ 2n-3$, which is less than $\dim (S^1\cap S^2)$.
This shows that
$$
\int_{S^{12}}\Om_{(0,q+1)}^3( z,\zeta )\wedge f(\zeta)=0, \quad q>0.
$$
By Stokes' formula and $\pd(S^{12})=\emptyset$, we obtain
$$
\db L^{23}_{12;q}f(z)=-\int_{S^{12}}\db_\zeta\Om_{(0,q+1)}^{23}( z,\zeta )\wedge f( \zeta )=\int_{S^{12}} \Om_{(0,q+1)}^{23}( z,\zeta )\wedge\db f(\zeta ).
\qedhere
$$
\end{proof}
In summary, we have the following local $\db$-solution operator for the concave case.
\th{cchf-closed}Let $1\leq q\leq n-2$. Let $(D^1,D^2 ,D^3 )$ be a  $(q+1)$-concave configuration. Let $U^1=D^2\setminus\ov {D^1}$ and  $S_+^1=\pd D^2\setminus D^1$. Assume $0<r_4<r_3/C_n$. Let $f\in
C^1_{(0,q)}(\ov{D^{1}})
$ be  $\db$-closed.
 On $D^{1}\cap B_{r_4}$, we have $$
f=\db H_qf
$$
with $H_q=H_q^{(1)}+H_q^{(2)}+H_q^{(3)}$. Here $ H^{(1)}_q$ and $H^{(2)}_q$  are given by \rea{hq1}, \rea{nhq2},    and
$$ H_{q}^{(3)}f=
L_{12;q-1}^{123}f-T_{B_{r_4}; q}L^{23}_{12;q} f.
$$
\eth

Although it is not used in this paper, for potential applications, it is worthy to state the following local homotopy formula if we have an extra negative Levi eigenvalue: if $\pd D^1$ is strictly $(q+2)$ concave, then
$
 L_{12;q}^{23}f=0
$
by \re{db13=0} for a $(0,q)$-form $f$.
Therefore, we have the following.
\th{cchf} Let $1\leq q\leq n-3$. Let $(D^1,D^2 ,D^3 )$ be a  $(q+2)$-concave configuration. Assume $r_4$ satisfies \rea{defr4}. Let $U^1=D^2\setminus\ov {D^1}$ and  $S_+^1=\pd D^2\setminus D^1$.  Let $f\in
C^1_{(0,q)}(\ov{D^{1}})
$ with $\db f\in C^1(\ov{D^1})$. On $D^{1}\cap B_{r_4}$, we have \gan{}
f=\db H_qf+  H_{q+1}\db f  
\end{gather*}
with $H_q=H_q^{(1)}+H_q^{(2)}+H_q^{(3)}$. Here  $ H^{(1)}_q$ and $H^{(2)}_q$  are given by \rea{hq1}, \rea{nhq2}   and
$$
H_{q}^{(3)}f=
L_{12;q-1}^{123}f-T_{B_{r_4};q}L^{23}_{12;q} f.
$$
\eth

As in the convex case, the integral kernels  in $H^{(2)}_q, H^{(3)}_q$  have smooth kernels, since $S^1_{+}, S^2, S^{12}$ do not intersect small neighborhood of the origin in $\pd D^1$. There is another main difference in kernels between the convex and concave cases. For the concave case, when the boundary  $\pd D$ is in $ \Lambda^{m}$,   the Leray functions \re{HL2pg81} is $\Lambda^{m-1}$ in $z$ when $\rho\in\Lambda^m$.  Anyway $H_q^{(1)}$ is  the main term to be estimated.

\setcounter{equation}{0}

\section{H\"{o}lder-Zygmund spaces and a Hardy-Littlewood lemma}\label{h-space}
In this section, we recall the H\"{o}lder--Zygmund norms and indicate how a Hardy-Littlewood lemma can be used to derive the estimates. We will also discuss various equivalent norms.

By
 a bounded Lipschitz  domain $D$ in $\rr^n$, we mean that  there are a finite open covering $\{U_i\}_{i=1}^N$ of $\ov D$, rigid affine transformations $A_i$, positive numbers $\del_0,\del_1,L$ such that    $A_i(\ov D\cap\ov{ U_i})$ is  defined $\del_1\geq x_n\geq R_i(x')$ with $x'\in[0,\del_0]^{n-1}$ and $|R_i(\tilde x')-R_i(x')|\leq L|\tilde x'-x'|$. We say that a constant $C(D)$ depending on $D$ is {\it stable} under small perturbations of $D$, if $C(\tilde D)$ can be chosen independent of $\tilde D$ if $A_i(\ov{\tilde D}\cap \ov{U_i})$ is defined by $\del_1\geq x_n\geq\tilde R_i(x')$, where  $|\tilde R_i(x')-R_i(x')|<\e$ and $|\tilde R_i(\tilde x')-\tilde R_i(x')|\leq (L+\e)|\tilde x'-x'|$ for some $\e>0$.

Let $\nn$ (resp. $\nn_+$) be the set of nonnegative (resp. positive) integers.
For $k\in\nn$, let $C^k(\ov D)$ be the space of $C^k$ functions on $\ov D$ with the standard norm $\|\cdot\|_{C^k(D)}$. When $a=k+\all$  with $0<\all<1$, let $C^a(\ov D)$ be the space of $C^k$ functions $f$ on $\ov D$ such that $$
\|f\|_{C^r(D)}:=\|f\|_{C^k(D)}+\sup_{
x,y\in D, x\neq y}\f{|\nabla^k  f(y)-\nabla^k  f(x)|}{|y-x|^{\all}} <\infty.
$$
  Define $\Del_hf(x)=f(x+h)-f(x)$ and
$\Del^2_h f(x)=f(x+2h)+f(x)-2f(x+h)$.
When $k\in\nn_+$, define $\Lambda^{k}(\rr^n)$ to be the space of all functions $f\in C^{k-1}(\rr^n)$  such that
\gan 
\| f\|_{\Lambda^{k}(\rr^n)}:=\|f\|_{C^{k-1}(\rr^n)} +\sup_{
h\neq0,x\in\rr^n}+ \f{|\Del^2_h\nabla^k  f(x)|}{|h|} <\infty.
\end{gather*}
We then define for $k\in\nn_+$
\gan
\|f\|_{\Lambda^{k}(D)}:=\inf_{\tilde f|_D=f}\|\tilde f\|_{\Lambda^{k}(\rr^n)},\\
{\Lambda}^{k}(D):=\{f\in C^k(\ov D)\colon \|f\|_{\Lambda^{k}(D)}<\infty\}.
\end{gather*}
For convenience, set $\Lambda^r(D):=C^r(\ov D)$ and $\|\cdot\|_{\Lambda^r(D)}:=\|\cdot\|_{C^r(D)}$ when $r\in(0,\infty)\setminus\nn$.

Following~\cite{MR2017700}*{Def.~3.15}, define an intrinsic norm in difference operators
$$
\|f\|^{{in}}_{\Lambda^{r}(D)}=\| f\|_{C^0( D)}+ \sup_{
x+jh\in D, j=0,\dots, ([r]+1); h\neq0}\left\{ \f{|\Del^{[r]+1}_hf(x)|}{|h|^{r}}\right\}, \quad \forall f\in\Lambda^r(D).
$$

It is known that when $D=\rr^n$, $|\cdot |^{in}_{\rr^n;r}$ and $|\cdot|_{\rr^n;r}$ are equivalent~\cite{MR0781540}*{Thm. 2.5.12, ~p.~110; Thm.~2.5.13 (i), p.~115}.
A classical theorem~\cite{MR0521808}*{Thm.~18.5, p.~63} says that given any function $f$ on a bounded Lipschitz domain satisfying $\|f\|^{in}_{\Lambda^r(D)}<\infty$ (without assuming $f\in\Lambda^r(D)$), there is an extension $\tilde f\in\Lambda^r(\rr^n)$ such that $\tilde f|_D=f$; consequently, $f\in\Lambda^r(D)$.

\begin{notation}To simplify notation, sometime we denote the norm $\|\cdot\|_{C^a(D)}$ by $\|\cdot\|_{D;a}$ or simply $\|\cdot\|_a$, and also we denote
$\|\cdot\|_{\Lambda^a(D)}$ by $|\cdot|_{D;a}$ or simply $|\cdot|_a$.
\end{notation}

\le{HL}Let $0<\beta\leq1$.
Let $ D\subset\rr^n$ be a bounded and connected Lipschitz domain. Suppose that   $f$ is in $C_{loc}^{[\beta]+1}( D)$ and
$$
|\nabla^{[\beta]+1} f(x)|\leq A\dist(x,\pd D)^{\beta-[1+\beta]}.
$$
Fix $x_0\in D$. Then
$
\|f\|^{in}_{ \Lambda^{\beta}(D)}\leq \|f\|_{C^0(D)}+C_\beta A,
$
where constants $C_0,C_\beta$ are stable under small perturbations of $D$.
\ele
\begin{proof}
Suppose $\beta=1$. We need  to estimate $\Del_h^2f(x)$ when $x$ is close to the boundary and $h\in\rr^n$ is small.
Take a boundary point $x_1$ of $\pd D$. We may assume that $x_1=0$. By definition, we may assume that $  D$ is defined by $x_n>g(x')$ where $g$ is a Lipschitz function satisfying
$|g(\tilde x')-g(x')|\leq L|\tilde x'-x'|$ with $L>1$. Then
$$
\dist(x,\pd D)\geq (x_n-g(x'))/{\sqrt{L^2+1}}.
$$
Suppose that $x,x-h, x+h$  are in $ D$ and  close to the origin.
 We first consider the special case  when $x,h$ satisfy
\eq{spcs}
x+th\in  D, \quad \dist(x+th,\pd D)\geq |h|, \quad \forall t\in[0,2].
\end{equation}
 Set $u(s)=f(x+h+sh)+f(x+h-sh)-2f(x+h)$. Thus $u(0)=u'(0)=0$ and $|u''(s)|\leq C_nLA|h|$.  This shows that $|u(1)|\leq C_nLA|h|$.

\begin{figure}
\centering
\begin{tikzpicture}[scale=.75]
\draw[very thick](pi-5,0.15)--(pi-3.5,0.15-.415)--(pi-1.5,0.15+.25) -- (pi-.3+.5,.25-.4)--(pi+.75+2,0.15)--(pi+3.5+2,-.2);


   \draw[dashed, 
   ] (pi-3,0.15-.2)-- +(35:5)
arc (35:145:5) -- cycle;

   \draw[thick, 
    ] (pi,-.1)-- +(35:5)
arc (35:145:5) -- cycle;


   \draw[dashed, 
   ] (pi+3,.125)-- +(35:5)
arc (35:145:5) -- cycle;

    \filldraw[black] (pi-3,.5+.1-.35) circle (2pt) (pi,.5+.1-.35) circle (2pt) (pi+3,.5+.1-.35) circle (2pt);

    \filldraw[black] (pi-3,.25+.75*4+1) circle (2pt) (pi,.25+.75*4+1) circle (2pt) (pi+3,.25+.75*4+1) circle (2pt);


  \draw[->]  (pi-3,.5+.1-.35)--(pi-3,.25+.75*4+.98) node[midway,right] {$\tilde h$};
    \draw[->,dashed] (pi-3,.5+.1-.35) --  (pi-.1,.5+.1-.35)node[midway, below] {$h$};

  \node at (pi-5.5,3-.75) {$V^*_{x-h}$};

    \node at (pi+2.625,3-.75) {$V^*_{x}$};

       \node at (pi,.5) {$x$};

   \node at (pi+5.125,3-.75) {$V^*_{x+h}$};

    \node at (pi+5,.125) {$D$};

\end{tikzpicture}
     \vspace{3ex}
    \caption{Lift points $x\pm h\in V^*_{x\pm h}$ to $x\pm\tilde h\in V^*_x$}
 \vspace{3ex}
 \label{fig:lifting}
\end{figure}

 For the general case,  set $\tilde h=(0',4\sqrt{L^2+1}|h|)$.  When $y\in D$ is close to the origin and $|\tilde h|<1$, we have
\ga
y+t\tilde h\in D, \quad  \dist(y+t\tilde h,\pd D)\geq t|\tilde h|/{\sqrt{L^2+1}}, \quad t\in[0,2].\label{yavn+}
\end{gather}
 Rewrite  $\Del^2_hf(x)= f(x+2h)+f(x)-2f(x+h)$ as
\al\label{thedecom}
\Del^2_hf(x) &=2f(x+2h+\tilde h)+2f(x+\tilde h)-4f(x+h+\tilde h)\\
\nonumber
&\quad - f(x+2\tilde h)-f(x+2h+2\tilde h)+2f(x+h+2\tilde h)\\
\nonumber
&\quad + f(x)+f(x+2\tilde h)-2f(x+\tilde h)\\
\nonumber
&\quad+ f(x+2h)+f(x+2h+2\tilde h)-2f(x+2h+\tilde h)\\
\nonumber
&\quad - 2f(x+h)-2f(x+h+2\tilde h)+4f(x+h+\tilde h).
\end{align}
We estimate each row on the right-hand side of \re{thedecom}. Let us  denote by $[a,b]$ the line segment connecting two points $a,b$ in $\rr^n$.
By   \re{yavn+}, we have
$$
[x+\tilde h,x+\tilde h+2h]\subset D, \quad
\dist(x+\tilde h+t h,\pd D)>|h|,
$$
for $t\in[0,2]$.
Therefore, we can estimate the first row using the estimation for the special case  \re{spcs} in which $x$ is replaced by $x+\tilde h$. The second row is estimated similarly.   For the third row, by \re{yavn+} and $x\in D$, we get
$$
\dist(x+\tilde h+s\tilde h,\pd D)> (1+s)|h|> (1-|s|)| h|, \quad s\in[-1,1].
$$
Take $u(s)=f(x+\tilde h+s\tilde h)+f(x+\tilde h-s\tilde h)-2f(x+\tilde h)$.  This yields $|u''(s)|\leq C_n AL| h|/{(1-s)}$ for $s\in[0,1]$ and
$|u(1)|\leq C_nAL| h|$ by
 $
u(1)=\int_{0}^1(1-s)u''(s)\, ds$.
 This  gives us the desired estimate for the third row.  The last two rows in \re{thedecom} can be estimated similarly.

The proof is standard for $0<\beta<1$, by using the three-term decomposition  $\Del_hf(x)=\Del_{\tilde h}f(x)+\Del_hf(x+\tilde h)-\Del_{\tilde h}f(x+h)$ and \re{yavn+}.
\end{proof}

Very recently, Shi--Yao~\cite{Shi-Yao-3}*{Def.~3.3, Thm.~1.1} show among other results  that when $r=k+\beta$ with $k\in\nn$ and $0<\beta\leq1$,  $\|f\|_{\Lambda^r(D)}$ and $\sum_{|\all|\leq k} \|\nabla^\all f\|_{\Lambda^\beta(D)}^{in}$ are equivalent using Rychkov's universal extension for Besov spaces $B^r_{p,q}(D)$. They also show the stability of constants in their estimates under small deformation of Lipschitz domains~\cite{Shi-Yao-3}*{Rmk.~6.9}.

Let us also observe that  Dispa ~\cite{MR2017700}*{Thm.~3.18} shows that there exists a constant $C_r(D)$ such that
\eq{key-est}
 \|f\|_{\Lambda^r(D)}\leq C_r(D)\|f\|^{{in}}_{\Lambda^r(D)}, \quad r>0.
\end{equation}
 By definition and Taylor formulae, it is easy to see that $|f|^{{in}}_{\rr^n;r}\leq C_r|f|_{\rr^n;r}$ and hence
$\|f\|^{{in}}_{\Lambda^r(D)}\leq C_r \|f\|_{\Lambda^r(D)}.
$
Therefore, we have
\begin{cor}\label{all-equiv-excluding-int} The norms   $ \|\cdot \|^{in}_{\Lambda^r(D)}$ and $\|\cdot \|_{\Lambda^r(D)}$ are equivalent.
\end{cor}

The proof in~\cite{MR2017700} uses only the non-universal part of Rychkov extension in ~\cite{MR1721827}*{Thms.~2.2-2.3}. Thus one can see  that the constant $C_r(D)$ in \re{key-est} depends only on the Lipschitz norm of the graph function of $\pd D$ and hence it is stable under small perturbations of $D$.

We also need the following  result.
\begin{lemma}\label{inpd} Let $r=k+\beta$ with $k\in\nn$ and $0<\beta\leq1$. Let $D$ be a bounded Lipschitz domain.
 Then
\begin{align}\label{finL}
\|f\|^{in}_{\Lambda^{r}(D)}&\leq C_r(D)(\|f\|_{C^0(D)}+\|\nabla^kf\|^{in}_{\Lambda^\beta(D)}).
\end{align}
Further, $C_r(D)$ is stable under a small perturbation of Lipschitz domains.
\end{lemma}
\begin{proof}It suffices  to verify that if $x+jh$ are in $ D$ for $j=0,\dots, (k+1)$,
$$
|\Del^{k+[\beta]+1}_hf(x)|  \leq C_r|h|^{r}  \|\nabla^k  f\|^{in}_{\Lambda^\beta(D)}.
$$
The case $k=0$ follows from definition. Suppose $k>0$. By a simpler analogue of \re{thedecom}, we have a three-term decompose
$
\Del_hf(x)=\sum_{i=1}^3 \Del_{h^{i}}f(y^i)
$
where $y^i=x+\tilde h^i$, $\tilde h^i, h^i$ depend only on $h$ such that $[y^i,y^i+h^i]\subset D$ and $|h^i|\leq C|h|$.
  Thus
with $\pd_if(x)=\f{\pd}{\pd x_i}f(x)$ we have
$$
\Del_{h^i}f(y^i)=\sum_{j=1}^nh_j^i\int_0^1\pd_j f(y^i+sh^i)\, ds.
$$
Repeat this.   We can write $\Del^k_hf(x)$ as a linear combination of
\eq{01k}
\int_{[0,1]^k} \pd^\all f(\tilde y^\all+s_1h^{(1)}+\cdots +s_kh^{(k)})\, ds \times \prod_{i=1}^k h^{(i)}_{\all_i}
\eeq
where $|\all|=k$, and $\tilde y^\all=x+\tilde h^\all$, $\tilde h^\all,h^{(i)}$ depend only on $h$. Further, $|h^{(i)}|\leq C_k(D)|h|$.
Thus we obtain \re{finL}  applying a three-term decomposition to the integrand in \re{01k} when $0<\beta<1$ or applying \re{thedecom} when $\beta=1$.
\end{proof}

We will need interpolation.
\begin{prop}\label{first-equiv} Let $E\colon C^0(\ov D)\to C^0_0(\rr^n)$ be the Stein extension. Then
\gan{}
 \|f\|_{C^r(D)}\leq\|Ef\|_{C^r(\rr^n)}\leq C_r(D) \|f\|_{C^r(D)},\quad   \ r\in[0,\infty); \\
 \|f\|_{\Lambda^r(D)}\leq \|Ef\|_{\Lambda^r(\rr^n)}\leq C_r(D)  \|f\|_{\Lambda^r(D)}, \quad   \ r\in(0,\infty).
\end{gather*}
\end{prop}
\begin{proof} The first and the third inequalities follow from definitions.
The other two are properties of the Stein extension operator; see~\ci{MR3961327}*{Prop.~3.11} for details.
\end{proof}

\begin{defn}[\cite{MR0230022}*{p.~167}] Let $(X_0,\|\cdot\|_0), (X_1,\|\cdot\|_1)$ be two Banach spaces contained in a Banach space $\cL X$. Set
\aln{}
K(t,f;X_0,X_1)&=\inf_{f=f_0+f_1}\{\|f_0\|_0+t\|f_1\|_1\},\quad t>0,\\
 \|f\|_{\theta;X_0,
X_1}&=\inf_{t>0} t^{-\theta}K(t,f;X_0,X_1),\quad  0<\theta<1.
\end{align*}
We now specialize the intermediate spaces via H\"older spaces and set
$$
\|f\|^*_{\Lambda^{r_\theta}(D);r_0,r_1}:=\|f\|_{\theta;C^{r_0}(\ov D),C^{r_1}(\ov D)}, \quad r_\theta=(1-\theta)r_0+\theta r_1,
$$
where $r_0,r_1$ are given and $0<\theta<1$.
\end{defn}
\begin{lemma}\label{fsEf} $\|f\|^*_{\Lambda^{r_\theta}(D);r_0,r_1}$ and $\|Ef\|_{\Lambda^{r_\theta}(\rr^n)}$ are equivalent.
\end{lemma}
\begin{proof}  It is well-known that $\|Ef\|_{\Lambda^{r_\theta}(\rr^n)}$ is equivalent to $\|Ef\|_{\Lambda^{r_\theta}(\rr^n);r_0,r_1}^*$.  Let us write  $\|\cdot\|_{C^r(\rr^n)}$, $\|\cdot\|_{C^r(\ov D)}$ as $\|\cdot\|_r,\|\cdot\|_{D,r}$ respectively.
 Thus  it suffices to show
\eq{fseffs}
\|f\|^*_{\Lambda^{r_\theta}(D);r_0,r_1}\leq \|Ef\|^*_{\Lambda^{r_\theta}(\rr^n);r_0,r_1}\leq C_r(D)\|f\|^*_{\Lambda^{r_\theta}(D);r_0,r_1}.
\end{equation}
Set
$  K(t,f;D):=K(t,f; C^{r_0}(\ov D), C^{r_1}(\ov D))$ and $$
K(t,Ef):=K(t,Ef; C^{r_0}(\rr^n), C^{r_1}(\rr^n)).$$
Suppose $f=f_0+f_1$ on $\ov D$. Since $E$ is linear, then $Ef=Ef_0+Ef_1$. By definition,
$
K(t,Ef)\leq \|Ef_0\|_{r_0}+t\|Ef_1\|_{r_1}\leq
C_{r_0,r_1}(D)(\|f_0\|_{D;r_0}+t\|f_1\|_{D;r_1})$. Thus $$K(t,Ef)\leq C_{r_0,r_1}(D)K(t,f;D),\quad \|Ef\|^*_{\Lambda^{r_\theta}(\rr^n);r_0,r_1}\leq C_r(D)\|f\|^*_{\Lambda^{r_\theta}(D);r_0,r_1}.
$$

If $Ef=\tilde f_0+\tilde f_1$, we have $f=\tilde f_0|_{\ov D}+\tilde f_1|_{\ov D}$. By definition,
$
K(t,f;D)\leq \|\tilde f_0|_{\ov D}\|_{D;r_0}+t\|\tilde f_1|_{\ov D}\|_{D;r_1}\leq \|\tilde f_0\|_{r_0}+t\|\tilde f_1\|_{r_1}.
$
Thus $K(t,f;D)\leq K(t,Ef)$ and hence $\|f\|^*_{\Lambda^{r_\theta}(D);r_0,r_1}\leq \|Ef\|^*_{\Lambda^{r_\theta}(\rr^n);r_0,r_1}$. We have verified \re{fseffs} and hence the lemma.
\end{proof}

In summary, we have the following.
\begin{cor}\label{all-equiv} The norms $\|Ef\|_{\Lambda^r(\rr^n)}, \|f\|_{\Lambda^r(D)}, \|f\|_{\Lambda^r(D)}^{{in}},$
and $\|f\|_{\Lambda^r(D);r_0,r_1}^*$ are equivalent with constant factors that are stable under small perturbations of $D$ by Lipschitz domains.
\end{cor}

\begin{prop}[\cite{MR0230022}*{Thm.~3.2.23, p.~180}]\label{int-norm} Let $(X_0,\|\cdot\|_{X_0}), (X_1,\|\cdot\|_{X_1})$ $($resp. $(Y_0,\|\cdot\|_{Y_0}), (Y_1,\|\cdot\|_{Y_1}))$ be two Banach spaces continuously embedded in a Banach space $\cL X$ $($resp. $\cL Y)$.   If $L\colon \cL X\to\cL Y$ is   linear and
$$
\|Lf_i\|_{Y_i}\leq M_i\|f_i\|_{X_i}, \quad i=0,1
$$
then $\|Lf\|_{\theta;Y_0,Y_1}\leq M_0^{1-\theta}M_1^\theta\|f\|_{\theta;X_0,X_1 }$ for $0<\theta<1$.
\end{prop}

\begin{cor}[\cite{GG}*{Prop.~3.4}]
\label{int-dist} Let $D\subset\rr^n$ be a bounded Lipschitz domain.  Let $E\colon C^0(\ov D)\to C^0_0(\rr^n)$ be the Stein extension. Then
$$
|[\nabla, E]f(x)|\leq C_r(D)\|f\|_{\Lambda^r(D)}\dist(x,D)^{r-1}, \quad r>1.
$$
Further, $C_r(D)$ is stable under small perturbations of $\pd D$.
\end{cor}

\setcounter{thm}{0}\setcounter{equation}{0}
\section{$\f{1}{2}$-gain estimates for local homotopy operators}\label{sec:.5est}

In this section, we derive the estimates for homotopy operators. We will give   precise estimates which are potentially useful for applications. We remark that the local estimates do not require the forms to be $\db$ closed. Throughout the paper, we denote by $A(z,\zeta,\nabla\rho^1,\dots,\nabla^k\rho^1)$ a polynomial in $z,\zeta,\nabla\rho^1,\dots,\nabla^k\rho^1$.

We first consider the $(n-q)$ convex case. In this case the result is essentially in~\cite{MR3961327}. We simplified proof by using Lemmas~\ref{HL}, \ref{inpd} and equivalent norms,  which are applicable for $C^2$ domains.
\begin{thm}\label{conv-est} Let $r\in(1,\infty)$ and $1\leq q\leq n-1$.
Let $(D^1,D^2)$ be an $(n-q)$-convex configuration.
The homotopy operator $H_q$ in \rta{hf-c} satisfies
\al{}\label{hqfr12-c}
\|H_q\var\|_{\Lambda^{r+1/2}(D^{1}\cap D^2_{r_3})}&\leq C_r(\nabla\rho^1,\nabla^2\rho^1)\|\var\|_{\Lambda^r(D^{12})}, \quad r_1/2<r_3<3r_2/4,\\
\|H_q\var\|_{C^{3/2}(D^{1}\cap D^2_{r_3})}&\leq C(\nabla\rho^1,\nabla^2\rho^1)\|\var\|_{C^1(D^{12})}, \quad\  r_1/2<r_3<3r_2/4.
\label{C-hqfr12-c}
\end{align}
Moreover,  $C_r(\nabla\rho^1,\nabla^2\rho^1)$ is stable under small $C^2$ perturbations of $\rho^1$.
\end{thm}
\begin{proof}
To simplify  notation, write $\|\var\|_{\Lambda^r(D^{12})}, \|H_q\var\|_{\Lambda^{r+1/2}(D^{12}_{r_3})}$ as $|\var|_{r}, |H_q\var|_{ r+1/2}$.
   Recall the homotopy operator
$$
H_{q}\var=(H^{(1)}_q+H^{(2)}_q)\var
$$
with
$$
H^{(1)}_q\var=R_{U^1\cup D^{12}}^0 E \var+R_{U^1 }^{01}[\db,E] \var.
$$
Near the origin, each term in $H^{(2)}_q\var$ given by \re{nhq2} is a linear combination of integrals of the form
$$ Kf(z):=\int_{S^I}\f{A(\nabla_\zeta\rho^1,\nabla^2_\zeta\rho^1,\zeta, z) f(\zeta )}{(g^1\cdot(\zeta-z))^a(g^2\cdot(\zeta-z))^b |\zeta-z|^{2c}}\, dV
$$
where $A$ is a polynomial, $f$ is the coefficients of $\var$, and  $a,b,c$ are   integers. Note that $S^I$ is one of $S^1_+, S^{12}, S^2$. For $z$ close to the origin and $\zeta\in S^I$ we have $|g^i(\zeta,z)\cdot(\zeta-z)|>c$.  Thus
$
| Kf|_{r+1/2}\leq C_r \|f\|_{0}.
$
Here and in what follows $C_r$ denotes a constant depending on $\nabla\rho^1,\nabla^2\rho^1$.

We now estimate the main term $H^{(1)}\var$.
Decompose it as
 \eq{dbECV}
H^{(1)}\var(z)=\int_{D^{12}\cup U^1}\Om_{(0,q)}^{0}(z,\zeta)\wedge E\var(\zeta)+ \int_{U^1}\Om_{(0,q)}^{01}(z,\zeta)\wedge[\db,E]\var(\zeta).
 \end{equation}
Denote the first integral by $K_1\var$. Since $D^{12}_{r_3}$ is contained in a relatively compact ball in $D^{12}\cup U^1$, by estimates on the Newtonian potential~\cite{MR999729}*{p.~316}, we have
$$
|K_1\var|_{r+1/2}\leq C|\var|_{r-1/2}.
$$
Note that the above is proved in~\cite{MR999729} when $r$ is not an integer. When it is an integer, the estimate follows the interpolation for the Zygmund spaces.
The last integral in \re{dbECV} can be written as a linear combination of
 \gan 
  K_2f(z):=\int_{U^1} f(\zeta )\f{A(\nabla_\zeta\rho^1,\nabla^2_\zeta\rho^1,\zeta, z)N_{1}(\zeta-z )}
 {\Phi^{n-j}(z,\zeta)|\zeta -z |^{2j}}\, dV(\zeta), \quad 1\leq j<n,\\
 \Phi(z,\zeta)=g^1(z,\zeta)\cdot(\zeta-z),
\end{gather*}
where  $A$ is a polynomial, $N_m(\zeta)$ denotes a monomial in $\zeta,\ov \zeta$ of degree $m$, and $f$ is a coefficient of the form $[\db,E]\var$ and hence by \nrc{int-dist}
$$
|f(\zeta)|\leq C_r|\var|_r \dist(\zeta,D^1)^{r-1}.
$$
Here and in what follows, $N_m(\zeta)$ denotes a monomial of $\zeta$ with degree $m\geq0$
whose coefficient is a smooth function in $z,\zeta$.    

 Fix $\zeta_0\in\pd D^1$.  We first choose local coordinates such that $s_1(\zeta),s_2(\zeta)$ and $t(\zeta)=(t_3,\dots, t_{2n})(\zeta)$ vanish at $\zeta_0$, $D^1$ is defined by $s_1<0$,  and
 \ga \label{LbPhi}
  |\Phi(z,\zeta)|\geq  c_*(\dist(z,\pd D^1)+ s_1(\zeta)+|s_2(\zeta)|+|t(\zeta)|^2),\\  C|\zeta-z|\geq |\Phi(z,\zeta)|\geq c_*|\zeta-z|^2,\\
|\zeta-z|\geq c_*(\dist(z,\pd D^1)+ s_1(\zeta)+|s_2(\zeta)|+ |t(\zeta)|),
 \label{LbPhi+}
  \end{gather}
  for $z\in D^1,\zeta\notin D^1$.

Let $r=k+\all$ with integer $k\geq1$ and  $0<\all\leq1$. Assume  $z\in D^{12}_{r_3}$.

Consider first the case   $0<\all<1/2$. 
We have
$|\pd^{k+1}Kf(z)|\leq C_r|\var|_{r}I(z)$, where
\eq{dk+2}
I(z):=\int_{[0,1]\times[-1,1]^{2n-1}}\f{s_1^{r-1}(s_1+|s_2|+|t|)^{-(2j+b-1)}\, d s_1ds_2dt}{(\dist(z,\pd D^1)+ s_1+|s_2|+|t|^2)^{(n-j)+a}}
\end{equation}
with $a+b=k+1$ and $1\leq j<n$. The worst term occurs when $j=n-1$ and $a=k+1$. Therefore, using polar coordinates for $t(\zeta)$ we obtain $I(z)\leq C   I_1(z)$ for
\eq{I_1(z)}
  I_1(z):= \int_{[0,1]^3}\f{s_1^{r-1}t^{2n-3}\, d s_1ds_2dt}{(\dist(z,\pd D^1)+ s_1+s_2+t^2)^{k+2}(s_1+s_2+t)^{2n-3}}.
\eeq
 Using polar coordinates for $(s_1,s_2)$, we obtain $I_1(z)\leq CI_2(z)$ with
\eq{hatI}
  I_2(z)= \int_{[0,1]^2}\f{s^{\all+1} \, d sdt}{(\dist(z,\pd D^1)+ s+t^2)^{3}}
\leq C\dist(z,\pd D^1)^{(\all+1/2)-1}
\end{equation}
where the last inequality is proved in \cite{MR3961327}*{Lem.~4.3 (i)}.

Note that the estimate is valid for $r=1$ and $\all=0$, when we  replace $|\var|_1$ by $\|\var\|_1$. This gives us
\re{C-hqfr12-c}.

The case $1/2\leq \all\leq 1$ is obtained by  interpolation via \rp{int-norm} in which $X_i=\Lambda^{r_i}( D)$,  $Y_i=\Lambda^{r_i+1/2}( D)$ for $i=0,1$ with $r_0=k+\all/3$, $r_1=k+1+\all/3$, $r_\theta=k+\all$, and $\theta=2\all/3$.

To be used later, we can also give a direct estimate as follows. When $1/2\leq \all\leq 1$, we use $|\nabla^{k+2}Kf(z)|\leq C_r |\var|_{r} I(z)$ defined by \re{dk+2} in which $a+b=k+2$. Then $I(z)\leq C  I _3 (z)$ for
\eq{tildeI}
  I _3 (z):= \int_{[0,1]^3}\f{s_1^{r-1}t^{2n-3}\,d s_1ds_2dt}{(\dist(z,\pd D^1)+ s_1 +s_2 +t^2)^{k+3}(s_1+s_2+t)^{2n-3}},\end{equation}
which is less than $C \dist(z,\pd D^1)^{(\all+1/2)-2}$ because $\all+1/2<2$.
\end{proof}

\le{henint}Let $\beta,\mu_1\in[0,\infty)$, $0\leq\la\leq\mu_1$,  and $0<\del<1$. Set
$$ 
\beta':=\beta-\f{\mu_1+\la-3}{2}.
$$ 
 Then for $m\geq0$,
$$ 
\int_{[0,1]^{3}}\f{ s_1^\beta(\del+s_1+s_2+t^2)^{-1-\mu_1}t^{m}}{(\del+s_1+s_2+t)^{m+\la-\mu_1}
}\, ds_1ds_2dt
<
\begin{cases}
C\del^{\beta'},&\beta'<0;\\
C,&\beta'>0.
\end{cases}
$$ 
\ele
\begin{proof} In our application, we have $m=2(n-1)-1\geq1$. It suffices to consider $m=0$. We consider the integral in the following regions.

$(i)$ $s_2>\max\{s_1,\delta,t\}$. On this region the integral is less than
$$
\int_{s_2=\del}^1\int_{s_1=0}^{s_2}\int_{t=0}^{s_2}\frac{s_1^{\beta}}{s_2^{\la+1 }}\, dt ds_1ds_2\leq
\int_{\del}^1 s_2^{\beta+1-\la}\,  ds_2<\int_{\del}^1s_2^{\beta'-1}\, ds_2<C\del^{\beta'}
$$
if $\beta'<0$. Also the integral bounded by a constant if $\beta'>0$. The same bounds can be obtained for the integral on regions $(ii)$ $s_1>\max\{\delta,s_2,t\}$  and $(iii)$ $\del>\max\{s_1,s_2,t\}$.

$(iv)$ $t^2>\max\{\del,s_1,s_2\}$. On this region, the integral is less than
$$
\int_{t=\sqrt\del}^1\int_{s_1={t}^2}^1\int_{s_2={t}^2}^1\frac{s_1^{\beta}\, ds_1ds_2dt}{t^{2(1+\mu_1)+\la -\mu_1}}\leq  \int_{t=\sqrt\del}^1t^{2\beta+2-\mu_1-\la }\, dt<C\del^{\beta-(\mu_1+\la -3)/2}
$$
when $\beta'
<0$.
When $\beta'>0$, the integral is bounded by a constant.

 $(v) \  t^2<\del+s_1+s_2<t$. On this region, the integral is less than
\eq{triple}
\int_{(s_1,s_2)\in[0,1]^2}\int_{t= \del+s_1+s_2}^{\sqrt{\del+s_1+s_2}} \frac{s_1^{\beta}\, dt ds_1ds_2}{(\del+s_1+s_2)^{1+\mu_1}t^{\la -\mu_1}}.
\end{equation}
Suppose $\beta'<0$. Since $\la\leq\mu_1$, the latter is less than
$$
 \int_{(s_1,s_2)\in[0,1]^2}  \frac{s_1^{\beta}\,  ds_1ds_2}{(\del+s_1+s_2)^{1+\mu_1+(\la-\mu_1-1)/2}}<C\del^{\beta'},
$$
where the last inequality is obtained
by computing  integrals for $\del\geq\max\{s_1,s_2\}$, $ s_1\geq\max\{\del,s_2\}$, and
$s_2\geq\max\{\del,s_1\}$. When $\beta'>0$, it is straightforward that the integral is bounded above by a constant $C$.
\end{proof}

We now treat the concave case. We start with the following.
\le{abcd} Let $D\subset \rr^n$ be a bounded Lipschitz domain.  Let $|\cdot|_a=\|\cdot\|_{\Lambda^a(D)}$ and $\|\cdot\|_a=\|\cdot\|_{C^a(\ov D)}$. Then
\al{}
\label{Hcov}
\|u\|_{a+b}\|v\|_{c+d}&\leq C_{a,b,c,d}(\|u\|_{a+b+d}\|v\|_d+\|u\|_a\|v\|_{c+b+d}),\quad a,b,c,d\geq0;\\\label{Zcov}
|uv|_a&\leq C_a(|u|_a\|v\|_0+\|u\|_0|v|_a),\quad  a\in(0,\infty), \quad a>0;\\
|u|_{a+b}|v|_{c+d}&\leq C_{a,b,c,d}(|u|_{a+b+d}|v|_d+|u|_a|v|_{c+b+d}),\quad a,d>0, b,c\geq0.
\label{eq:abcd}
\end{align}
Further, $C_{a,b,c,d}$ and $ C_a$ depend on the Stein extension $E_D$ and are stable under small perturbations of the Lipschitz domain $D$.
\ele
\begin{proof}
 A proof of \re{Hcov} is  in~\cite{MR2829316}*{Prop.~A.4 (iii)}.
When $D=\rr^n$,  \re{Zcov} is proved in \cite{MR2768550}*{Cor. 2.86, p.~104}. For the $D$, we  use the Stein extension $E_D$ and get
 $|uv|_a\leq|(E_Du)(E_Dv)|_a\leq C_a(|E_Du|_a\|E_Dv\|_0+\|E_Du\|_0|E_Dv|_a)\leq C_a'(|u|_a\|v\|_0+\|u\|_0|v|_a)$.

 We now prove \re{eq:abcd}. Let us first consider the case that $D=\rr^n$ and $u,v$ have compact support in $\rr^n$. Let us use a Littlewood-Paley decomposition
  $$
  u=\sum_{i\geq0}  u_i, \qquad u_0=\Phi\ast u, \quad u_j=\Psi_{2^{-j}}\ast u, \quad j>0
  $$
 where $\Psi_{t^{-1}}(x)=t^{n}\Psi(tx)$; for the definition of   $\Phi$ and $\Psi$, see \cite{MR3243741}*{sect.~1.4.3}. By \cite{MR3243741}*{Thm.~1.4.3, p.~46}, there is a constant $C(n,r)$   such that
 $$
|u|_r/{ C(n,r)}\leq\sup_{j\geq0} 2^{jr}\|f_j\|_0\leq C(n,r)|u|_r.
 $$
 Decompose  $v=\sum_{k\geq0} v_k$ analogously. Then
 \begin{align*}
|u|_{a+b}|v|_{c+d}&\leq C_{a,b,c,d} \sup_{j,k\geq0} 2^{j(a+b)+k(c+d)}\|u_j\|_0\|v_k\|_0\\
 &\leq C_{a,b,c,d} \sup_{j,k\geq0}( 2^{j(a+b+c)+kd}+2^{ja+k(b+c+d)})\|u_j\|_0\|v_k\|_0\\
 &\leq C_{a,b,c,d}' (|u_{a+b+c}|v|_d+|u|_a|v|_{b+c+d}).
  \end{align*}
 Note that we have used $b\geq0$ and $c\geq0$ for the second inequality. By a similar proof of \re{Zcov}, the general case of \re{eq:abcd} is obtained through the  extension $E_D$.
 \end{proof}

 We now derive our main estimates, using Lemmas~\ref{HL}, \ref{inpd} and equivalent norms.

\begin{thm}\label{concave-est} Let $r\in(1,\infty)$ and $1\leq q\leq n-2$.  Let $(D^1,D^2,D^3)$ be a $(q+1)$-concave configuration.
The homotopy operators $H_q$ in Theorems~$\ref{cchf-closed}$ and $\ref{cchf}$ satisfy
\al{}
\label{hqfr12}
\|H_q\var\|_{\Lambda^{r+1/2}(D^{12}_{r_4/2})}&\leq C_r(\nabla\rho^1,\nabla^2\rho^1) (\|\rho^1\|_{\Lambda^{r+5/2}}\|\var\|_{C^1(\ov{D^{12}})}+  \|\var\|_{\Lambda^r(D^{12})}),\\
\label{hqfr120}
\|H_q\var\|_{C^{1+\gamma}(D^{12}_{r_4/2})}&\leq C_\gamma(\nabla\rho^1,\nabla^2\rho^1)(1+\|\rho^1\|_{C^{3+\gamma}} )\|\var\|_{C^1(\ov {D^{12}})}.
\end{align}
Here $ 0\leq\gamma\leq1/2$.  Moreover,  $C_r(\nabla\rho^1,\nabla^2\rho^1), C_\gamma(\nabla\rho^1,\nabla^2\rho^1)$ are stable under small $C^2$ perturbations of $\rho^1$.
\end{thm}
\begin{proof}
To ease notation, write
$$
\|\var\|_{\Lambda^a(D^{12})},\quad \|\var\|_{C^a(D^{12})},\quad\|H_q\var\|_{\Lambda^{a}(D^{12}_{r_4/2})},\quad \|H_q\var\|_{C^{a}(D^{12}_{r_4/2})}
$$
as  $|\var|_{a}, \|\var\|_{a}, |H_q\var|_{ a}, \|H_q\var\|_{ a}$,  respectively. We also write $\|\rho^1\|_{\Lambda^a(U_1)}, \|\rho^1\|_{C^a(U_1)}$ as $|\rho^1|_a,\|\rho^1\|_a$.
    Recall the homotopy operator
$$
H_{q}\var=(H^{(1)}_q+H^{(2)}_q+H^{(3)}_q)\var
$$
where $\var$ has type $(0,q)$ and
$$
H^{(1)}_q\var=R_{U^1\cup D^{12}}^0 E \var+R_{U^1 }^{01}[\db,E] \var.
$$
Near the origin, each term in $H^{(2)}_q\var, H^{(3)}_q\var$ is a linear combination of integrals of $K_0f_0$ with
$$
K_0f_0(z)=A(z,\nabla_z\rho^1,\nabla_z^2\rho^1)\tilde K_0f_0(z)
$$
with $f_0$ being a coefficient of $\var$. Here $\tilde K_0$, which  involves only $\nabla\rho^1$,
is defined by
$$
\tilde K_0f_0(z):=\int_{S^I}\frac{ f_0(\zeta )(\zeta,\bar\zeta)^e\, dV(\zeta)}{
(g^1(z,\zeta)\cdot(\zeta-z))^a(g^2(z,\zeta)\cdot(\zeta-z))^b(g^3(z,\zeta)\cdot(\zeta-z))^c|\zeta-z|^{2d}}
$$
where $a,b,c,d$ are nonnegative integers, $e\in\nn^{2n}$, and $S^I$ is one of $S^1_+, S^{12}, S^2$. Therefore, for the $g^1,g^2,g^3$ that appears in the kernel, we have by \re{gicZ+}
$$
|g^i(z,\zeta)\cdot(\zeta-z)|\geq c_0
$$
when $z\in D^{12}_{r_4/2}$ and $\zeta\in S^I$. By \re{Zcov}, we have
$
|\tilde K_0f_0|_{r+1/2}\leq C_r|\rho^1|_{r+3/2}\|f_0\|_{0}
$
and
$$
|K_0f_0|_{r+1/2}\leq C_r|\rho^1|_{r+5/2}\|\tilde K_0f_0\|_{0}+C_r |\tilde K_0f_0|_{r+1/2}.
$$
Here and in what follows $C_r$ denotes a constant depending on $\rho,\pd\rho,\pd^2\rho$.
Therefore, we have
$$
|K_0f_0|_{r+1/2}\leq C_r|\rho^1|_{r+5/2}\|f\|_0.
$$

We now estimate the main term $H^{(1)}\var$, which is decomposed as
 \eq{dbE}
H^{(1)}\var(z)=\int_{D^{12}\cup U^1}\Om_{(0,q)}^{0}(z,\zeta)\wedge E\var(\zeta)+ \int_{U^1}\Om_{(0,q)}^{01}(z,\zeta)\wedge[\db,E]\var(\zeta).
 \end{equation}
Denote the first integral by $K_1\var$. By an estimate on Newtonian potential~\cite{MR999729}*{p.~316}, we have
$$
|K_1\var|_{r+1/2}\leq C_r|\var|_{r-1/2}.
$$
The last integral in \re{dbE} can be written as a linear combination of
 \gan \label{defnKf}
K_2f(z):=A(z,\nabla_z\rho^1,\nabla_z^2\rho^1 )\tilde K_2f(z)
\end{gather*}
where $f$ is a coefficient of the form $[\db,E]\var$.  Cors.~\ref{all-equiv} and \ref{int-dist} imply
$$
|f|_{r-1}\leq C_r(D^1)|\var|_{D^1,r}, \quad
|f(\zeta)|\leq C_r(D^1)|\var|_r \dist(\zeta,D^1)^{r-1}.
$$
Note that $f$ vanishes on $\ov D$. Further
\gan{}
\tilde K_2f(z):=\int_{U^1} f(\zeta )\f{ N_{1}(\zeta-z )}
 {\Phi^{n-j}(z,\zeta)|\zeta -z |^{2j}}\, dV(\zeta), \quad 1\leq j<n,\\
 \Phi(z,\zeta)=g^1(z,\zeta)\cdot(\zeta-z).
 \end{gather*}

We have
\al{}\label{K2fr}
|K_2f|_{r+1/2}&\leq C_r(\|\rho^1\|_{2})(|\rho^1|_{r+5/2}\|\tilde K_2f\|_0+ |\tilde K_2f|_{r+1/2}),\qquad\ \  r>1;\\
\|K_2f\|_{1+\gamma}&\leq C_\gamma(\|\rho^1\|_{2})(1+\|\rho^1\|_{3+\gamma})\|\tilde K_2f\|_0+ \|\tilde K_2f\|_{1+\gamma}),\quad 0\leq\gamma\leq \f{1}{2}.
\label{K2fr+}\end{align}

The rest of the proof is devoted to the proof of
\begin{gather}\label{tK2f}
|\tilde K_2f|_{r+1/2}\leq C_r( |\rho^1|_{r+5/2}\|\var\|_{1}+ |\var|_{r}), \quad r>1,
\end{gather}
or \re{K2fr} directly for easy cases. When $r=1$, we simply replace $|\var|_1$ by $\|\var\|_1$ in the proof below.
Then combining   above estimates yields the proof for \re{hqfr12}-\re{hqfr120}.

\medskip

A technical difficulty to prove \re{tK2f} is that our approach relies on \re{Zcov} and \re{Hcov} for Zygmund norms, instead of \re{eq:abcd}. See for instance \re{Kmu<} below.  The computation is tedious and our main observation is that the kernel of $\tilde K$ involves only $\nabla^1\rho$ instead of $\nabla^2\rho$.

Let $r=k+\all$ with integer $k\geq1$ and  $0<\all\leq1$.
In the following cases, we will apply \rl{henint} several times. For clarity, we will specify values $(\beta,\beta')$ in \rl{henint} via  $(\beta_i,\beta_i')$ when we use the lemma.

\medskip

$(i)$ $0<\all<1/2$. Recall that  the kernel of $\tilde K_2$  involves only first-order derivatives of $\rho^1$ in $z$-variables and it does not involve $\zeta$-derivatives of $\rho^1$.   Further, $\Phi$ is a linear combination of $\zeta_j-z_j$.    Since $\rho^1\in \Lambda^{k+\all+5/2}\subset C^{k+2}$, we   can express $\nabla^{k+1}\tilde K_2f$ as a sum of
\eq{rhotimesK}
K_{\mu,\nu}^{(k+1)}f:= \nabla^{1+\nu'_1}\rho^1\cdots\nabla^{1+\nu'_{\mu_1}}\rho^1K_\mu f.
 \end{equation}
Here and in what follows $\nabla^{\ell}\rho^1$ stands for a partial derivatives of $\rho^1(z)$ in $z,\ov z$ of order $\ell$ and
 \eq{Kmuf}
 K_\mu f(z):=\int_{U^1}\f{f(\zeta)N_{1-\mu_0+\mu_1-\nu_1''-\cdots-\nu_{\mu_1}''+\mu_2}(\zeta-z)}
{(\Phi(z,\zeta))^{n-j+\mu_1}|\zeta-z|^{2j+2\mu_2}}dV.
\end{equation}
Here $1\leq j<n$,
$\nu_i''=0,1$, and
$$ 
\mu_0+
\mu_2+\sum(\nu'_i+\nu''_i)\leq k+1,\quad \nu_i'+\nu_i''\geq1.
$$ 
To estimate \re{rhotimesK}, we  use \re{LbPhi}-\re{LbPhi+}: For $z\in D^1$, $\zeta\in U\setminus D^1$ and $D^1$ defined by $s_1<0$, we have
\gan  
  |\Phi(z,\zeta)|\geq  c_*(\dist(z,\pd D^1)+ s_1(\zeta)+|s_2(\zeta)|+|t(\zeta)|^2), \\
 C|\zeta-z|\geq |\Phi(z,\zeta)|\geq c_*|\zeta-z|^2,\\ |\zeta-z|\geq c_*(\dist(z,\pd D^1)+ s_1(\zeta)+|s_2(\zeta)|+ |t(\zeta)|.
  \end{gather*}
Consequently, the worst term for $K_\mu f$ occurs when  $j=n-1$. Note that we use $|\nabla\rho^1(z)|\geq c$ in \re{rhotimesK}. Thus
the worst term for $K_{\mu,\nu}^{(k+1)}f$ also occurs when $\mu_0+\mu_2$ is absorbed into $\sum\nu_i''$.
Therefore, it suffices to estimate terms with $\mu_0=\mu_2=0$ and $j=n-1$, which are assumed now.

Throughout the proof, we assume that $\nu_i''=0$ and hence $\nu_i'\geq1$ for $i\leq \mu_1'$ and $\nu_i''=1$ for $i>\mu_1'$. Thus
$$
\mu_1'+\mu_1''=\mu_1, \quad \sum_{i\leq \mu_1'}\nu_i'=\sum_{i\leq\mu_1}\nu_i'.
$$
Thus \re{Kmuf} is simplified in the form
\ga{}\label{Kmuf+}
K_\mu f(z):=\int_{U^1}\f{f(\zeta)N_{1+\mu_1'}(\zeta-z)}
{(\Phi(z,\zeta))^{1+\mu_1}|\zeta-z|^{2(n-1)}}dV,\\
\label{mu1'}
\mu_1\leq\mu_1''+\sum_{i=1}^{\mu_1}\nu_i'\leq \sum(\nu_i'+\nu_i'')\leq k+1.\end{gather}

We will apply  \rl{henint} a few times, using
 $$ 
 \la:=\mu_1-\mu_1',\quad m=2n-3.
 $$ 
We need $\beta_1\geq0$ and $
\beta_1'\leq\beta_1-\f{\mu_1+\la-3}{2}.
$
Thus we take
$$
\beta_1'=\all-1/2<0, \quad
\beta_1=\max\Bigl\{0,\all-1/2+ \f{\mu_1+\la-3}{2}\Bigr\}.
$$
We can verify that $\beta_1\leq r-1$. By  \rl{henint},    we obtain
 for $z\in D^{12}_{r_4/2}$,
\eq{Kmuf0}
|K_\mu f(z)|\leq |\var|_{\beta_1+1}\dist(z,\pd D^1)^{\beta_1'}.
\end{equation}
Thus,
\al \label{Kmu<}
|K_{\mu,\nu}^{(k+1)} f(z)|&\leq C\|\rho^1\|_{1+\nu_1'}\cdots\|\rho^1\|_{1+\nu_{\mu_1}'}|\var|_{\beta_1+1}
 \dist(z,\pd D^1)^{\all-1/2}.
\end{align}
 When $\beta_1=0$, we actually need to replace the above $|\var|_{\beta_1+1}$ by $\|\var\|_1$. Then we  the desired estimate easily since $ \|\rho^1\|_{1+\nu_1'}\cdots\|\rho^1\|_{1+\nu_{\mu_1}'}\leq C\|\rho^1\|_{1+\sum \nu_i}$ and $\sum\nu_i<r+1$. The case  all   $\nu_i'\leq1$  can be estimated directly via \re{hatI} to get   $|K_\mu(z)|\leq C|\var|_{r}\dist(z,D^1)^{(\all+1/2)-1}$
 and hence \re{tK2f}.

 Assume  now  $\beta_1>0$  and  $\nu_1'\geq2$. Thus $\mu_1'\geq1$.
 Let $x_+=\max\{0,x\}$ and
\eq{def-gamma}
\gamma:=\sum(\nu'_i-1)_+=-\mu_1'+\sum\nu_i'.
\end{equation}
  Then
$
\|\rho^1\|_{1+\nu'_1}\cdots\|\rho^1\|_{1+\nu'_{\mu_1}} \leq \|\rho^1\|_{2+\gamma}.
$
Then
$$
\beta_1+\gamma\leq\Bigl[\all-2+ \f{\mu_1+\la}{2}\Bigr]-\mu_1'+{\textstyle\sum}\nu_i'=\all-2- \f{\mu_1'}{2}+\mu_1''+{\textstyle\sum}\nu_i'\leq \all-\f{\mu_1'}{2}+k-1\leq r-\f{3}{2},
$$
where the second inequality is obtained by \re{mu1'} and $\mu'_1\geq1$. Therefore,
\aln{}
&\|\rho^1\|_{1+\nu_1'}\cdots\|\rho^1\|_{1+\nu_{\mu_1}'}|\var|_{\beta_1+1}\leq C_r \|\rho^1\|_{2+\gamma}\|\var\|_{\beta_1+1} \\
&\qquad\leq C_r'\|\rho^1\|_{2}\|\var\|_{1+\beta_1+\gamma}+C_r'\|\rho^1\|_{2+\gamma+\beta_1}\|\var\|_1\leq
C_r''|\var|_{r}+C_r''\|\rho^1\|_{r+1/2}\|\var\|_1.
\end{align*}

By an estimation similar to \re{C-hqfr12-c}, we have $\|\tilde K_2f\|_{1}\leq C(\|\rho^1\|_2)\|f\|_0$.
The above estimate also works for $r=k=1$ (and $\all=0$), when $|\var|_1$ is replaced by $\|\var_1\|_{1}$. Namely, we   have
$ 
\|\tilde K_2f\|_{1+\gamma}\leq C(\|\rho^1\|_{2})\|\rho^1\|_{3+\gamma}\|f\|_0.
$ 
Thus, \re{K2fr+} yields \re{hqfr120}.

\medskip

$(ii)\  1/2<\all\leq1$. In this case we can take an extra derivative since $\rho^1\in\Lambda^{k+\all+5/2}\subset C^{k+3}$. Write
$\nabla^{k+2}\tilde K_2f$ as a sum of
$$
 K_{\mu,\nu}^{(k+2)}f:=\nabla^{1+\nu'_1}\rho^1\cdots\nabla^{1+\nu'_{\mu_1}}\rho^1K_\mu f
 $$
where $K_\mu f$ is defined by \re{Kmuf}.
As before, the worst term occurs for $j=n-1$ and $\mu_0=\mu_2=0$ in $K_\mu$, which are assumed now.
  Then \re{mu1'} in which  $k$ is                                                                                                             replaced by $k+1$  becomes
  $$  
\mu_1\leq\mu_1''+\sum_{i=1}^{\mu_1}\nu_i'\leq k+2.
$$

To apply \rl{henint},
we  need to   $\beta_2\geq0$ and
$
\beta_2'\leq\beta_2-\f{\mu_1+\la-3}{2}.
$  
Thus we take
$$
\beta_2'=\all-3/2<0,
\quad \beta_2=\max\{0,\all-3/2 + (\mu_1+\la-3)/2\}.
$$
 Recall $\la=\mu_1-\mu_1'$. We can verify that $\beta_2\leq r-1$.  We obtain for $z\in D^{12}_{r_4}$
\aln{} 
|K_{\mu,\nu}^{(k+2)} f(z)|&\leq C_r\|\rho^1\|_{1+\nu_1'}\cdots\|\rho^1\|_{1+\nu_{\mu_1}'}|\var|_{\beta_2+1}
 \dist(z,\pd D^1)^{\all-3/2}.
\end{align*}
When $\beta_2=0$, we can replace $|\var|_{\beta_2+1}$ by $\|\var\|_1$. Then, using $$\|\rho^1\|_{1+\nu_1'}\cdots\|\rho^1\|_{1+\nu_{\mu_1}'}\leq C\|\rho^1\|_{k+3}\leq C'|\rho^1|_{r+5/2},
$$
 we get the  estimate immediately. When all  $\nu_i'\leq1$, using \re{tildeI} we get  $|K_\mu(z)|\leq C_r|\var|_{r}\dist(z,D^1)^{\all-3/2}$ and hence \re{tK2f}.

Assume now   $\beta_2>0$ and  $\nu_1'\geq2$. Thus $\mu_1'\geq1$. Let
 $\gamma$ be given by \re{def-gamma}.
  Then
we have
$$
\beta_2+\gamma\leq\Bigl[\all-3 + \f{\mu_1+\la}{2}\Bigr]-\mu_1'+{\textstyle\sum}\nu_i'=\all-3- \f{\mu_1'}{2}+\mu_1''+{\textstyle\sum}\nu_i'\leq \all-1- \f{\mu_1'}{2}+k\leq r-\f{3}{2}.
$$
As in previous case, we get \aln{}
&\|\rho^1\|_{1+\nu_1'}\cdots\|\rho^1\|_{1+\nu_{\mu_1}'}|\var|_{\beta_2+1}\leq
C_r(|\var|_{r}+\|\rho^1\|_{r+1/2}\|\var\|_1).
\end{align*}

\medskip

$(iii)\  \all=1/2$.  We need to estimate $|\tilde K_2f|_{r+1/2}$ with $r+1/2=k+1$.

 Recall that $\rho^1\in\Lambda^{k+3}$.  On  $D^{12}_{r_4/2}$, we write $\nabla^k\tilde K_2f$ as
a sum of $$
 K_{\mu,\nu}^{(k)}f:=\nabla^{1+\nu'_1}\rho^1\cdots\nabla^{1+\nu'_{\mu_1}}\rho^1K_\mu f.
 $$
To estimate $|K^{(k)}_{\mu,\nu}f|_1$, we will use   \rl{HL} to estimate $|K_\mu f|_1$.

As before, the worst term occurs for $j=n-1$ and $\mu_0=\mu_2=0$ in $K_\mu$, which are assumed now.
Then $K_\mu f$ has the form \re{Kmuf+}, while \re{mu1'} has the form
\eq{mu1k}
\mu_1\leq\mu_1''+\sum_{i=1}^{\mu_1}\nu_i'\leq k.
\eeq

 By \re{Zcov}, we   have up to a constant factor
\al\label{halfcase}
 |K_{\mu, \nu}^{(k)}f|_1
 &\leq
  |\nabla^{1+\nu'_1}\rho^1\cdots\nabla^{1+\nu'_{\mu_1}}\rho^1|_1\|K_\mu f\|_0\\
  &\quad
 +\|\nabla^{1+\nu'_1}\rho^1\cdots\nabla^{1+\nu'_{\mu'_1}}\rho^1\|_0|K_\mu f|_1=:\operatorname{I}+\operatorname{II}.
 \nonumber
 \end{align}

 We first estimate $\operatorname{I}$.
Assume first  $\f{\mu_1+\la-3}{2}\geq0$. We set
 $$
 \beta_3'=\e\in(0,1), \quad \beta_3= \beta_3'+\f{\mu_1+\la-3}{2}\geq0.
 $$
Recall $\la=\mu_1-\mu_1'$. Then $\beta_3\leq r-1$.
 Then by \rl{henint} in which $\beta_3'>0$ we obtain
 $$\|K_\mu f\|_0\leq   C\|f\|_{\beta_3}.
 $$
 We also have
 $$ \|\nabla^{1+\nu'_1}\rho^1\cdots\nabla^{1+\nu'_{\mu_1}}\rho^1\|_1\leq C_r\sum \|\nabla^{2+\nu'_i}\rho^1\|_0\prod_{j\neq i}\|\nabla^{1+\nu'_j}\rho^1\|_0\leq C_r'\|\rho^1\|_{2+\tilde\gamma}
 $$
 with $\tilde\gamma=\sum\nu_i'$.   Thus \begin{align*}
\beta_3+\tilde\gamma&=\f{\mu_1+\la-3}{2}+\e +\sum\nu_i'
\\
&=\e-  \f{\mu_1'}{2}-\f{3}{2}+\mu_1''+\sum\nu_i'
\leq\e-  \f{\mu_1'}{2}-\f{3}{2}+ k\leq r-2+\e.
\end{align*}
Thus   we get
$$
\operatorname{I}\leq C_r\|\rho^1\|_{2+\tilde\gamma}\|\var\|_{1+\beta_3} \leq C_r'
\|\rho^1\|_{2}\|\var\|_{r-1+\e}+C_r'\|\rho^1\|_{r+\e}\|\var\|_1.
$$
If $\f{\mu_1+\la-3}{2}<0$, we set $\beta_3=0$ and $\beta_3'=-\f{\mu_1+\la-3}{2}>0$. By \rl{henint} with $\beta_3'>0$, we get $\|K_\mu f\|_0\leq C\|f\|_{\beta_3}=C\|\var\|_1$. Hence
$\operatorname{I}\leq C_r\|\rho^1\|_{2+\tilde\gamma}\|\var\|_1\leq C_r\|\rho^1\|_{r+2}\|\var\|_1$.

Finally, we estimate $\operatorname{II}$ in \re{halfcase}.   \rl{HL} says that we can  estimate $|K_\mu f|_1$  via the pointwise estimate of $\nabla^2K_\mu f$.
Write $\nabla^2K_\mu f$ as a sum of
$$
 K_{\tilde \mu,\tilde \nu}^{(k)}f:=\nabla^{1+\tilde\nu'_1}\rho^1\cdots\nabla^{1+\tilde\nu'_{\tilde\mu_1}}\rho^1K_{\mu+\tilde \mu} f
$$
with
$$
\tilde\mu_1\leq\tilde\mu_1''+\sum\tilde\nu_i'\leq 2.
$$
The worst terms occur in the forms
$$
\operatorname{II}':=\nabla^3\rho^1K_{\mu+\tilde \mu}f,\quad
 \operatorname{II}'':=\nabla\rho^1\nabla^2\rho^1K_{\mu+\tilde{\tilde\mu}}f
$$
for $\tilde\mu=1$ and $\tilde{\tilde\mu}=2$, respectively.

We first estimate $\operatorname{II}'$ where $\tilde\mu=1$. We have
\eq{KmuTm}
K_{\mu+\tilde\mu} f(z):=\int_{U^1}\f{f(\zeta)N_{1+\mu_1'+\tilde\mu_1'}(\zeta-z)}
{(\Phi(z,\zeta))^{1+\mu_1+\tilde\mu_1}|\zeta-z|^{2(n-1)}}dV(\zeta).
\end{equation}
Note that the worst term occurs when $\tilde\mu_1=\tilde\mu=1$.
Let us combine the indices as follows.   Set $\hat \mu=\mu+\tilde\mu=\mu+1$, $\hat\mu_1=\mu_1+\tilde\mu_1=\mu_1+1$.  Set $\hat\nu_i'=\nu_i'$ for $i\leq\mu_1$, $\hat\nu'_{\hat \mu_1}=2$. The latter implies
$$\hat\mu_1'=\mu_1'+1\geq1.
$$
Set  $\hat\lambda=\hat\mu_1''=\hat\mu_1-\hat\mu_1'$. Then $\hat\mu_1''=\mu_1''$. We get
$$
\hat\mu_1\leq\hat\mu_1''+\sum\hat\nu_i'=\mu_1''+2+\sum\nu_i'\leq k+2.
$$
 Thus, we set
 $$
 \beta_4'=-1, \quad \beta_4=\max\Bigl\{0,-1+\f{\hat\mu_1+\hat\la-3}{2}\Bigr\}.
 $$
  Then
 $$
 |\operatorname{II}'(z)|\leq C_r\|\rho^1\|_{1+\tilde\nu_1'}|\var|_{\beta_4+1} \dist(z,\pd D^1)^{\beta_4'}, \quad \tilde\nu_1'=2,
 $$
where $|\var|_{\beta_4+1}$ is replaced by $\|\var\|_1$ when $\beta_4=0$.  For the new indices, let $
\hat\gamma=\sum(\hat\nu'_i-1)_+.
$
 We have $$\| \nabla^{1+\hat\nu'_1}\rho^1\cdots\nabla^{1+\hat\nu'_{\hat\mu_1}}\rho^1\|_0 \leq C_r\|\rho^1\|_{2+\hat\gamma}.$$
Assume first that $\beta_4>0$. Then by $\hat\mu'_1\geq1$, we get from $(\hat \mu_1+\hat\la)/2-\hat\mu_1'=\hat\mu_2''-\hat\mu_1'/2$ and $r=k+1/2$
 $$
\beta_4+\hat\gamma\leq\Bigl[-\f{5}{2}+\f{\hat\mu_1+\hat\la}{2}\Bigr]- \hat\mu_1'+ \sum\hat\nu_i'\leq-\f{\hat \mu_1'}{2}-\f{5}{2}+ (k+2)\leq r-\f{3}{2}.
$$
Thus we obtained the desired estimate for $|\operatorname{II}'(z)|$ from
 $$
\|\rho^1\|_{2+\hat\gamma}|\var|_{1+\beta_4}\leq C'_r \|\rho^1\|_{r+1/2}\|\var\|_1+C_r'\|\rho^1\|_{2}\|\var\|_{r-1/2}.
$$
When $\beta_4=0$,  with $k\geq1$ and \re{mu1k} we  simply use $\hat\gamma=1+\sum(\nu_i'-1)_+\leq k\leq r-1/2$. Thus $ \|\rho^1\|_{2+\hat\gamma}\|\var\|_{1+\beta_4}\leq C_r\|\rho^1\|_{r+3/2}\|\var\|_1$, which gives us the desired estimate.

We now estimate  $\operatorname{II}''$.
Then we have $\tilde{\tilde\mu}=2$, and hence
$K_{\mu+\tilde{\tilde\mu}} f(z)$ has the form \re{KmuTm} in which $\tilde\mu, \tilde\mu_1',\tilde\mu_1$ are replaced by $\tilde{\tilde\mu}, \tilde{\tilde\mu}_1',\tilde{\tilde\mu}_1$ respectively with $$
\tilde{\tilde\mu}_1\leq\tilde{\tilde\mu}_1''+\sum\tilde{\tilde\nu}_i'\leq 2.
$$
The worst term occurs when $\tilde{\tilde\mu}_1=2$.
Let us combine the indices as before.   Set $\hat \mu=\mu+2$, $\hat\mu_1=\mu_1+\tilde{\tilde\mu}_1$, $\hat\mu_1'=\mu_1'+\tilde{\tilde\mu}'$,  and $\hat\lambda=\hat\mu_1''=\mu_1''+\tilde{\tilde\mu}''$. Let $\hat\nu_i'=\nu_i'$ for $i\leq\mu_1$ and  $\hat\nu'_{\hat \mu_1}=2$. The latter implies
$\hat\mu_1'=\mu_1'+1\geq1.
$ Then
\eq{hatmu1}
\hat\mu_1''+\sum\hat\nu_i' \leq k+2.
\end{equation}
 Thus, we set
 $$
 \beta_5'=-1, \quad \beta_5=\max\Bigl\{0,-1+\f{\hat\mu_1+\hat\la-3}{2}\Bigr\}.
 $$
 Then
 $$
 |K_{\mu+\tilde{\tilde\mu}} f(z)|\leq C_r|\var|_{\beta_5+1} \dist(z,\pd D^1)^{\beta_5'}.
 $$
 where $|\var|_{\beta_5+1}$ is replaced by $\|\var\|_1$ when $\beta_5=0$.
The values of $\beta_5,\beta_5'$, $\hat\mu_1'$,  and condition \re{hatmu1} are the same as in the previous case. The same computation yields the desired estimate
for $\operatorname{II}''$.

In summary,  we have obtained desired estimates for $\|\nabla^{k_1} H_q\var\|^{in}_{\Lambda^\beta(D^{12}_{r_4/2})}$ where $r+1/2= k_1+\beta$ with $ k_1\in\nn$ and $0<\beta\leq1$, by using \rl{HL}. Then \rl{inpd} yields the estimate for $\|H_q\var\|^{in}_{\Lambda^{r+1/2}(D^{12}_{r_4/2})}$. The proof is complete by the equivalence of norms in \nrca{all-equiv-excluding-int}.
 \end{proof}

\setcounter{thm}{0}\setcounter{equation}{0}
\section{An estimate of $\db$ solution for $(0,1)$ forms via Hartogs's extension}\label{sec:1form}

 In this section we obtain the regularity of functional $\db$-solutions for $(0,1)$-forms can be achieved via Hartogs's extension for concavity domains that require merely $C^2$ boundary.

We formulate the following result.
\begin{prop}\label{prop:hartogs}Let  $D\subset U\subset\cc^n$ be defined by $\rho<0$ with $\rho\in C^2$. Suppose that $U\cap \pd D$ is strictly  $2$-concave for each $\zeta\in\pd D\cap U$. Assume that $f\in\Lambda^r(D)$ with $r>1$ and $\db u_0=f$ on $D$ with $u_0\in C^0(\ov D)$. Then for open sets $U',U''$ satisfying $U''\Subset U'\Subset U$, we have
\eq{}
\|u_0\|_{\Lambda^{r+1/2}(U''\cap D)}\leq C_r(U',U'',\nabla\rho,\nabla^2\rho)(\|f\|_{\Lambda^{r}(D)}+\|u_0\|_{C^0(\ov D)}).
\end{equation}
Furthermore, $C_r(U',U'',\nabla\rho,\nabla^2\rho)$ is stable under a small $C^2$ perturbation of $\rho$.
\end{prop}
\begin{proof}
 For each $\zeta\in U''\cap\pd D$, by \rp{concave-rho} we find a local biholomorphic map $\psi_\zeta$ such that $\psi_\zeta(0)=\zeta$, and $\rho\circ\psi_\zeta=a_\zeta \rho_\zeta$ has the form
$$
\rho_\zeta(z)=-y_n-3|z_1|^2-3|z_{2}|^2+\sum_{j>2}\la_j|z_j|^2+o(|z|^2), \quad |\la_j|<1/2.
$$
   Then  $D_{\del_0,\e_0}:=(\Del_{\del_0}^{n-1}\times\Del_{\e_0})\cap\psi_\zeta^{-1} (D\cap U)$ contains
   $ \widetilde D_{\del_0,\e_0}=\{z\in \Del_{\del_0}^{n-1}\times\Del_{\e_0}\colon\tilde\rho(z)<0\}$, where
$$
\tilde\rho(z)=-y_n-2|z_1|^2-2|z_{2}|^2+\f{3}{4}\sum_{j=3}^n|z_j|^2.
$$
When $\del,\e$ are sufficiently small and $\del<\e$, both $D_{\del_0,\e_0},\widetilde D_{\del_0,\e_0}$ contain
the Hartogs's   domain
$$H_{\del,\e}:=(0,\dots,0,\sqrt {-1}\del)+
(\Del_{\e+\e^2}\setminus\Del_{\e-\e^2})\times\Del_\del^{n-1}\cup\Del_{\del} ^{n}.
$$
 Here  $H_{\del,\e}$ is the shaded region in \rf{fig:hartogs}.
 Note that $\widetilde D_{\del_0,\e_0}$ has smooth boundary and $\pd D_{\del_0,\e_0}\cap\pd \widetilde D_{\del_0,\e_0}=\{0\}$.

By \rl{concave-rho}, we know that the $C^2$ norm of $\psi_\zeta$ is bounded uniformly  and $\e$ can be chosen uniform in $\zeta\in \pd D$.
Further, $\e$ can be chosen uniformly for small $C^2$ perturbation of $\rho$.

Since $\pd\widetilde D_{\del,\e}$ is smooth, by \rt{concave-est} we have a solution $u$ on $\widetilde D_{\del,\e}$ for a possibly smaller $\del,\e$ such that $\db u= \psi_\zeta^*f$ and $|u|_{\widetilde D_{\del,\e};r+1/2}\leq C_r|f|_{D;r}$.
Then $u_0\circ\psi_\zeta-  u$ admits a holomorphic extension $h$  to
$$
 \widehat D_{\del,\e}:=  \widetilde D_{\del,\e}\cup \Delta_\delta^n
$$
via the Cauchy formula:
$$
h(z)=\f{1}{2\pi}\int_{|\zeta_1|=\e}\frac{h(\zeta_1,z')}{\zeta_1-z_1}\, d\zeta_1.$$
Here $\Del_\e\times\{z'\}$ is the disk indicated by a dotted line in \rf{fig:hartogs}.

\begin{figure}\centering
\begin{tikzpicture}[scale=.65]
 \draw[thin,domain=(-4.25:4.25), name path=Dzeta]
 plot[smooth](\x,-.1*\x*\x);

 \draw[very thick,domain=(-4.25:4.25), name path=De]
 plot[smooth](\x,-.18*\x*\x);

  \fill[gray!40](-3.8,-.9) rectangle (-3.5,1+.8) (-3.8,.75)rectangle (3.8,1+.8) (3.5,-.9) rectangle (3.8,1+.8) ;
\draw[->] (0,0)--(0,2) node[right]{$y_n$};
\draw[->] (0,0)--(5,0) node[right]{$x_1$};
   \draw[line width=.05cm, dashed](-3.8,-.8+.3*4)-- (3.8,-.8+.3*4);
        \draw[line width=.05cm, dashed](-(3.8,-.8+.3*2)-- ((3.8,-.8+.3*2);
        \draw[line width=.05cm, dashed](-(3.8,-.8)-- ((3.8,-.8);

 \node at (4.65,-2.6) {$D_{\del_0,\e_0}$};
   \node at (4.65,-1.3) {$\widetilde D_{\del_0,\e_0}$};
   \node at (.75,1.25) {$H_{\del,\e}$};
\end{tikzpicture}
   \vspace{3ex}
    \caption{$D_{\del,\e}\subset\widetilde D_{\del,\e}$ with $\pd D_{\del,\e}\in C^2$ and $\pd\widetilde D_{\del,\e}\in C^\infty$}
      \vspace{4ex}
       \label{fig:hartogs}
\end{figure}

Let $m$ be the largest integer less than $r+1/2$. By the Cauchy inequalities, we obtain $\|h\|_{\Del_{\del/2}^n;m+1}\leq C\|h\|_{L^\infty (H_{\del,\e})}\leq C_m(\|u_0\|_{ D;0}+|f|_{D;r})$.
In particular,
\eq{u0tD}
\|u_0\|_{\Lambda^{r+1/2}(D^*_\zeta)}\leq C_r(\|u_0\|_{C^0( D)}+\|f\|_{\Lambda^r(D)}), \quad D^*_\zeta:=\psi_\zeta(\widetilde D_{\del',\e'}).
\end{equation}
Here $\del',\e'$ are chosen so that $\widetilde D_{\del',\e'}$ are contained in $\Del_{\del/2}^n$. Note that $\del',\e'$ can be chosen uniformly in $\zeta\in U''\cap \pd D$. The constant $C_r$ is independent of $\zeta\in U''\cap \pd D$ and stable under a small $C^2$ perturbation of $\rho$.

\vspace{1ex}
\begin{figure}
\centering
\begin{tikzpicture}[scale=.65]
 \draw[very thick,domain=(-.5*pi:2.75*pi), name path=D1]
 plot[smooth]({(\x )},{.25*sin(\x r)});

  \draw[dashed,domain=-.25:1.75,rotate around={-3:(pi-3.2,-.15)}
  ]
  plot[smooth](pi-3.2+2*\x-1.35,-.15+.35*\x*\x);
  \node at (1-.75,-.5+1.5) {$D^*_{(z-w)^*}$};

\draw[thin,domain=-2:2.25, rotate around={-15:(pi+.1,.5)}
  ]
  plot[smooth](pi+.1+2*\x-.1,.05+.5*\x*\x-.05);

 \draw[dashed,domain=-1.5:1.25,rotate around={7:(pi+3,.1)}
  ]
  plot[smooth](pi+2.8+2*\x-.2,.5*\x*\x-.025);

   \node at (pi+3.1,-.5+1.5) {$D^*_{(z+w)^*}$};
   \node at (pi+5,.5) {$D$};

    \filldraw[black] (pi-3,.25-.1) circle (2pt) (pi,.25-.1) circle (2pt) (pi+3,.25-.1) circle (2pt);

    \filldraw[black] (pi-3+.8,.25+.75*4) circle (2pt) (pi+.8,.25+.75*4) circle (2pt) (pi+3+.8,.25+.75*4) circle (2pt);


  \node at (pi+.4,-.5+1.5) {$D^*_{z^*}$};
\end{tikzpicture}
  \vspace{3ex}
    \caption{Lift   $z\pm w\in D^*_{(z\pm w)^*}$ into $D^*_{z^*}$}
     \vspace{3ex}
      \label{fig:non-concave}
\end{figure}

Next we need to estimate the $\Lambda^\all$ norms of $m$-th derivatives of $u_0$.
Set $v$ be a derivative of
$u_0$ of order $m$.  It remains to estimate
$
v(z+w)+v(z-w)-2v(z)$ when $ z,z\pm w\in D.
$

We may assume that $z$ is sufficiently close to the origin and $|w|$ is small. Let $z^*\in\pd D$ be the closed point to $z$.
 Suppose $\del=|w|$ is small. Let $\tilde w=\delta(z-z^*)/{|z-z^*|}$. Since $\pd D^*_{z^*}$ is tangent to $\pd D$ at $z^*$, then $z+t\tilde w$ and $z\pm w+t'\tilde w\in D^*_{z^*}$ for $t\in(0,2)$ and $t'\in[1,2]$.
Here we cannot claim any concavity of the intersection of $D$ with the complex line through $z,z-w$, as shown in \rf{fig:non-concave} for $\pd D$. Nevertheless, for $\pd D^*_{z^*}$ is tangent to $\pd D$ at $z^*$. Consequently, the smooth domain $D^*_{z^*}$ contains a cone $V_{z^*}$ with vertex $z^*$.
We now use decomposition \re{thedecom} for $u_0$.
We can estimate each row in \re{thedecom} via \re{u0tD} because  the triple points in each row are in same smooth domain $D^*_\zeta$ for $\zeta=(z-w)^*,z^*$ or $(z+w)^*$.
With the estimate for the Zygmund ratio, we can conclude that $v\in \Lambda^\all$ when $0<\all<1$.

We now consider $\all=1$. Since our solution operator $H_1$ in \rt{concave-est} is linear, by interpolation we conclude that when $f\in\Lambda^r$, $u_0=H_1f$ is in $\Lambda^{r+1/2}(D^{12}_{r_4/2})$. Therefore, the estimate for Zygmund ratio is valid   too.
\end{proof}

 To conclude this section, we note that \rt{conv-est}, \re{hqfr120} in \rt{concave-est} (for $C^2$ domains $D$) and \rp{prop:hartogs} yield the local version of \rt{regsol} $(a)$.
 We also remark that the $C^{1/2}$ estimate in~\cite{MR986248}*{Thm.~14.1} seems to require $\pd D\in C^{5/2}$ to repeat the proof of ~\cite{MR986248}*{Thm. 9.1}.

 \setcounter{thm}{0}\setcounter{equation}{0}
\section{Proof of \rt{regsol} via canonical solutions}\label{sec8}
The proof of regularity of the solutions from local to global uses some standard approaches. See Kerzman~\cite{MR0281944} for the case when $D$ is a domain in $\cc^n$. We will also derive a global estimate reflecting the norm convexity and this estimate will be used in the next section to prove \rt{regsol+}.

Let $D\subset X$ be a relatively compact $C^1$ domain with a $C^1$ defining function on $X$. Fix a relatively compact neighborhood $\cL U$ of $\ov D$.
Using a partition of unity on $\ov D$ and \re{Zcov}, we can define a Stein extension $E\colon C^0(D)\to C^0_0(\cL U)\cap C^0(X)$ such that
$$
\|Ef\|_{C^a(X)}\leq C_a\|f\|_{C^a(\ov D)}, \quad \|Ef\|_{\Lambda^a(X)}\leq C_a\|f\|_{\Lambda^a( D)}
$$
where $C_a$ depends on $\rho$ and but is independent of small $C^1$ perturbations of $\rho$. Using local coordinates on $\cL U$, for a fixed $L$ we can define a Moser smoothing operator
$
S_t\colon C^0_0(\mathcal U)\cap C^0(X)\to C_0^\infty(X)
$
such that
\begin{align*}
\|S_tf-f\|_{C^a(\cL U)}&\leq C_{a,b}t^{b-a}\|f\|_{C^b(\cL U)}, \quad0\leq a<b\leq a+L;\\
\|S_tf\|_{C^b(\cL U)}&\leq C_{a,b}t^{a-b}\|f\|_{C^a(\cL U)}, \quad \forall b\geq a\geq0;\\
\|S_tf-f\|_{\Lambda^a(\cL U)}&\leq C_{a,b}t^{b-a}\|f\|_{\Lambda^b(\cL U)}, \quad0<a\leq b\leq L+a;\\
\|S_tf\|_{\Lambda^b(\cL U)}&\leq C_{a,b}t^{a-b}\|f\|_{\Lambda^a(\cL U)}, \quad \forall b\geq a>0.
\end{align*}
See~\cite{GG}*{(3.19)-(3.22)} for details when $X=\cc^n$, and the general case follows from a partition of unity and  \re{Zcov}.

We start with the following.
\le{8.1}Let $ D\Subset X$ be a domain defined by a $C^2$ function $\rho<0$ and let $D_a$ be defined by $\rho<a$. Let $\rho_t=S_t\rho$, where $S_t$ is the Moser smoothing operator. Suppose that $\pd D$ is an $a_q$ domain. Let $ D^t_a$ be defined by $\rho_t<a$.  There exists $t_0=t_0(\nabla\rho,\nabla^2\rho)>0$ and $C>1>c>0$ such that if $0\leq t<t_0$, then
\begin{gather}
\label{14t}
\|\rho_t-\rho\|_{C^0(X)}\leq t/4,\\
\pd D^t_{-t}\subset D_{-ct}\setminus D_{-Ct} \quad \pd D^t_{t}\subset \ov{D_{Ct}}\setminus D_{ct}
\label{14t+}
\end{gather}
while $ D_{b}$ and $ D^t_{b}$ still satisfy the condition $a_q$ for $b\in(-t_0,t_0)$.
Here $t_0,c,C$ are independent of small $C^2$ perturbations of $\rho$.
\ele
\begin{proof}Let $\rho_t=S_t\rho$. We have
$
\|\rho_t-\rho\|_{C^0(X)}\leq C_2t^2\|\rho\|_{C^2(X)}.
$
This shows that
\eq{C2t2}
\dist(\pd D^t_s,\pd D_s)\leq C_2\|\rho\|_{C^2(X)}t^2,
\end{equation}
Using a local $C^1$ diffeomorphism  whose first component is $\rho$, we can verify that
  $\max\{|s|,|s'|\}<\|\rho\|_{C^1(X)}/C_1$  and $s'>s$,  we also have
$$
c_1(s'-s)\leq \dist( D_s,\pd D_{s'})\leq C_1(s'-s),\quad  c_1(s'-s)\leq \dist( D^t_s,\pd D^t_{s'})\leq C_1(s'-s).
$$
Thus \re{C2t2} implies that
\begin{align*}
\dist(\pd D^t_{t}, D)&\geq \dist(\pd D_{t}, D)-\dist(\pd D^t_{t},\pd D_{t})\\
&\geq c_1t-C_2\|\rho\|_{C^2(X)}t^2>c_1t/2.
\nonumber
\end{align*}

 Suppose $L_\zeta\rho$ has  $(q+1)$ negative Levi eigenvalues bounded above by $-\la$  or $(n-q)$ positive Levi eigenvalues bounded below by $\la$ for $\zeta\in U$, where $U$ is neighborhood of $\pd D$ and $\la$ is positive number. We find a subspace $W$ of $T_\zeta^{1,0}\pd D$ dimension $q+1$ such $L\rho(\zeta,v)\leq-\la$ for $v$ in the unit sphere of $W$. Projecting $W$ onto $\widetilde W\subset T_{\tilde \zeta}^{(1,0)}\pd D^t_a$ when $\tilde\zeta\in \pd D^t_a$ is sufficiently close to $\zeta$ and $t$ is close to zero. Then $\dim W\geq q+1$ and $L\rho(\tilde\zeta,v)\leq-\la/2$ for $\tilde v$ in the unit sphere of $\widetilde W$.
One can also verity that if $L_\zeta\rho$ has at least $(n-q)$ positive eigenvalues, so is $L_{\tilde\zeta}\rho^t$ when $\tilde\zeta$ is sufficiently close to $\zeta$. Therefore, $D^t_a$ still satisfies the condition $a_q$.\end{proof}

We need to deal with domains of which different boundary components may have different smoothness. Let $D$ be a relatively compact $a_q$ domain in a complex manifold $X$. Assume that the $(n-q)$ strictly convex components $b^+_{n-q}D$ are of class $C^2$ and $(q+1)$ strictly concave components $b^-_{q+1}D$ are of class $\Lambda^s$ with $s>2$. Using a partition of unity, we can find a defining function $\rho$ on $X$ such that $D$ is defined by $\rho<0$ and $\nabla\rho(\zeta)\neq0$ for $\zeta\in\pd D$, $\rho\in C^2(\ov{U'})$ and $\rho\in \Lambda^{s}(U'')$, where $U'$
(resp. $U''$) is a neighborhood of $b^+_{n-q} D$ (resp. $b^-_{q+1} D$) in $X$.
\begin{defn}\label{def|-}Set $\|\rho\|_{C^{a}}^-:=\max\{\|\rho\|_{C^2(\ov{U'})},\|\rho\|_{C^a(\ov{U''})}\}$
and
$$\|\rho\|_{\Lambda^{s}}^-:=\max\{\|\rho\|_{C^2(\ov{U'})},
\|\rho\|_{\Lambda^{s}(U'')}\}.$$
\end{defn}

 We know formulate the main result of this paper in details.

\th{regsol-full} Let $r\in(1,\infty)$ and $q\geq 1$.
Let $D:=\{\rho<0\}$ be a relatively compact   domain with $C^2$ boundary in a complex manifold $X$ satisfying the condition $a_q$. Let $V$ be a holomorphic vector bundle of finite rank over $X$.
 Then there exists a
linear $\db$ solution operator $H_q\colon   \Lambda_{(0,q)}^r({D}, V)\cap\db L^2_{loc}(D,V)\to \Lambda_{(0,q-1)}^{r+1/2}({D}, V)$  satisfying the following
\bpp
\item When
$q=1$ or $\pd D$ is strictly $(n-q)$ convex, we have $\|H_qf\|_{\Lambda^{r+1/2}(D)}\leq C_r(\nabla^1\rho,\nabla^2\rho)\|f\|_{\Lambda^{r}(D)}$.
\item When $q>1$ and the components of $\pd D$ that are $(q+1)$--concave are in $\Lambda^{r+\f{5}{2}}$, we have
$$\|H_qf\|_{\Lambda^{r+1/2}(D)}\leq C_{r}(\nabla\rho,\nabla^2\rho) ( (\|\rho\|^-_{C^{3 }})^m \|\rho\|^-_{\Lambda^{r+5/2 }}\|f\|_{C^{1}(\ov D)}+\|f\|_{\Lambda^{r}(D)}).$$
    \item In both cases, $H_q$ is independent of $r$ and hence $H_qf\in C^\infty(\ov D)$ when $f\in C^\infty(\ov D)$.
\epp
Further, $m$ is independent of $r,s$ and remain the same under small $C^2$ perturbation of $\rho$.  Furthermore, the constant $C_r(\nabla\rho,\nabla^2\rho)$ is  stable under small $C^2$ perturbations of $\rho$; more precisely there exists $\del(\nabla\rho,\nabla^2\rho)>0$ such that if  $\|\tilde\rho-\rho\|_{C^2(X)}<\del(\nabla\rho,\nabla^2\rho)$, then
   \eq{defupst}
   C_r(\nabla\tilde\rho,\nabla^2\tilde\rho)<C'_rC_r(\nabla\rho,\nabla^2\rho).
   \end{equation}
We emphasize that $C_r(\nabla\rho,\nabla^2\rho)$ involves an unknown constant that is $C_*$ from \rta{3.4.6}.
\eth
\begin{rem}
The stability of estimates on $\db$ solutions has been discussed extensively in literature; see Greene--Krantz~\cite{MR644667} for strictly pseudoconvex domains in $\cc^n$, Lieb-Michel~\cite{MR1900133} strictly pseudoconvex domains with smooth boundary in a complex manifold. The stability in terms \re{defupst} is called upper-stability in  Gan--Gong~\cite{GG} where the reader can find a version of lower stability and its use.
\end{rem}

\begin{proof}The proof is a combination of the following: the local regularity results obtained, Grauert's bumping method,    the stability of solvability of the $\db$-equation after the bumping is applied \cite{MR0179443}*{Thm.~3.4.1} (see \rt{3.4.6}  for the vector bundle version), and the interior estimates of $\db$-solutions on Kohn's canonical solutions.

We will complete the proof in three steps.

\medskip
\noindent{\it Step 1. Reduction to interior regularity.}
Let $D$ be a relatively compact subset of $ \mathcal U$,  defined by $\rho<0$ in $\mathcal U$ in $X$ with $\|\rho\|_{C^{5/2}(\mathcal U)}<\infty$   and $\nabla \rho\neq0$ when $\rho=0$. For each $p\in\pd D$, we have a configuration $(U_p,D_{p}^1, D_{r_2}^2,\psi_p,\rho^1_p, \rho^2)$ with $\rho^1_p=\rho\circ\psi_p^{-1}$ as in \rd{def:3.4}. Recall that $
D_{p}^1:=\psi_p(U_p\cap D)$. Let
$$
\om_p':=\psi_p^{-1}(D^1_p\cap D^2_{r_2}).
$$
Let $\chi\geq0$ be a smooth cut-off function with compact support in $B_{r_4/2}$  such that $\chi$ equals $1$ on $B_{r_4/3}$.  
Applying Theorems~\ref{conv-est} and~\ref{concave-est}.  On $D^{1}_p\cap D^2_{r_4/2}$, we solve the $\db$ equation $\db u=(\psi_p^{-1})^*f$ with $u\in \Lambda^{r}(D^{1}_p\cap D^2_{r_4/2})$.  Then $f_1=f-\psi_p^*\db(\chi u_p)$ is still $\db$-closed and  in $\Lambda^r$ as
$$
\hat f_1:=(\psi_p^{-1})^*f-\db(\chi u_p)=(1-\chi)(\psi_p^{-1})^*f-\db\chi\wedge u_p.
$$
Recall that
\gan 
\|u_p\|_{\Lambda^{r+1/2}(D_p^1\cap D^2_{r_4/2})}\leq C_r(|f|_{D,r}+|\rho|^-_{r+5/2} \|f\|_{C^1(D)}), \quad r>1,\\
\|\hat f_1\|_{C^{1}(D_p^1\cap D^2_{r_4/2})}\leq C_\e(1+\|\rho\|^-_{3}) \|f\|_{C^1(D)}.  
\end{gather*}

 In fact, setting $f_1=0$ on $X\setminus D$, we have $f_1\in \Lambda^r( D\cup\ov{ \tilde B_{r_4/4}(p)})$ for $\tilde B_{a}(p)=\psi_p^{-1}(B_{a})$.
Let $D_p$ be defined by $\rho_p:=\rho-\e\chi\circ\psi_p<0$. We have
$$
D\cup\om_p\subset D_p\subset D\cup \tilde B_{r_4/2}(p),
$$
where $\om_p$ is an open set containing $p$; for instance we can take $\om_p=\{\rho<\e/C\}\cap B_{r_4/3}$, where $C$ depends only on $\|\nabla\rho\|_0$. Then $\om_p\cap \pd D$ contains $\pd D\cap \tilde B_{r_4/3}(p)$ no matter how small $\e$ is.
As in~\cites{MR1900133,MR3848426},   we find finitely many $p_1,\dots, p_m\in\pd D$ independent of $\e$ so that $\{\om_{p_1},\dots, \om_{p_m}\}$ covers   $\pd D$. Also $\sum\chi\circ\psi_{p_j}>0$ on $\pd D$.

\begin{figure}
  \begin{tikzpicture}[scale=1]

   \draw [very thick] (-4,0) 
      to [out=275,in=180]
      (0,-3)
      to  [out=0,in=270]
      (4,0)
      to  [out=90,in=0]
       (0,3)
       to [out=180,in=90]
        (-4,0)
        ;

 \draw [very thick, dashed] (-4,0+.4) 
      to [out=275,in=180]
      (0,-3+.2)
      to  [out=0,in=270]
      (4.1,0+.2)
      to  [out=90,in=0]
       (0,3+.2)
       to [out=180,in=90]
        (-4,0+.4)
        ;

  \draw
   (3.75,0) ellipse[x radius=.5cm, y radius=1cm];

   \draw[rotate around={25:(3.7,1.2)}]
   (3.7,1.2) ellipse[x radius=.5cm, y radius=1cm];
    \node[black] at (3.45,1.2) {$p_k$};
     \filldraw[black] (3.7,1.2) circle (2pt);

       \draw[rotate around={325:(3.4,-1.7)}]
   (3.4,-1.7) ellipse[x radius=.5cm, y radius=1cm];
    \node[black] at (3.75,-1.7) {$p_\ell$};
     \filldraw[black] (3.34,-1.7) circle (2pt);

    \draw [thick] (3.52,-1.5) 
      to [out=60,in=220]
      (4.25,-.95)
      to  [out=40,in=280]
       (4.75,0)
        to [out=280,in=330]
      (4.25,1)
       to [out=150,in=300]
       (3.4,1.65);

       \node at (3.8,0) {$p_1$};
          \node at (5.25,0) {$D_1$};
                \node at (-2,.5) {$\pd  D_0\subset\bigcup\omega_j\subset D_m$};
                \node at (-2+1,-.5) {a small perturbation $\pd\widetilde  D_0$ of $\pd D_0$};
                \node at (3,1) {$\omega_k$};
                   \node at (3,0) {$\omega_0$};
                      \node at (3,-1) {$\omega_\ell$};
\draw[->](-3,.75)[out=90,in=340] to (-3.5,1.4);
\draw[->](-3,-.75)[out=280,in=5]
 to (-3.3,-1.3);

      \filldraw[black] (4,0) circle (2pt);
 \end{tikzpicture}
   \vspace{3ex}
   \caption{Stability of bumping:\emph{} $\pd\widetilde D_0\subset \cup\om_j\subset D_m$} 
     \vspace{3ex} \label{fig:bumps}
    \end{figure}

 With $\rho_0=\rho$, $D_0=D$ and $\e>0$, set
\eq{rhoj=}
\rho_j=\rho_{j-1}-\e\chi\circ\psi_{p_j}=\rho_0-\e\sum_{\ell\leq j}\chi\circ\psi_{p_\ell}
\end{equation}
 and $D_j:=(D_{j-1})_{p_j}\colon\rho_j<0$ for $j\geq1$. We have $D_j\setminus D_{j-1}\subset \tilde B_{r_4/2}(p_j)$. Also, $D_{j}$   contains $  D\cup \om_{p_j}$
 and $D_j\subset D_{j+1}$. Thus $D_{m}$ contains $\ov{D_0}$ and
 $
 \dist(\ov {D_0},\pd D_m)>1/C_{\e}.
 $

 Let $\widetilde D_0$ be a perturbed  domain of $D_0$, which is defined by $\tilde\rho<0$, and let $\widetilde D_m$ be defined by $\tilde\rho_m<0$, where $\tilde\phi_m$ is defined by \re{rhoj=} with $\rho_0$ being
replaced by $\tilde\rho_0=\tilde\rho$.

Fix $\e>0$ and then fix $\del_1$ such that
\gan\label{delta_1}
0<\max\{\|\rho_m-\rho_0\|_{C^2(X)},\del_1\}<\min\left\{c_*,\f{1}{2}\del(\nabla\rho,\nabla^2\rho),\f{1}{2}\del_*, \f{1}{2}t_0(\nabla\rho,\nabla^2\rho)\right\},\\
\{\rho<2\del_1\}\subset \{\rho_m<-\del_1\},
\end{gather*}
where $t_0(\nabla\rho,\nabla^2\rho)$ is given by \rl{8.1},   $\del(\nabla\rho,\nabla^2\rho)$
is given by Theorems~\ref{regsol-full} and \ref{3.4.6}, and $c_*,\del_*$ are given by Corollary~\ref{corb11}.

Suppose $\|\tilde\rho-\rho\|_2\leq \del_1$. Thus
$$
\{\tilde\rho<a-\del_1\}\subset\{\rho<a\}\subset\{\tilde\rho<a+\del_1\}.
$$
Since $\tilde\rho_m-\rho_m=\tilde\rho-\rho$, then  $\|\tilde\rho_m-\rho_m\|_2\leq \del_1$ and hence
$$
\{\tilde\rho_m<a-\del_1\}\subset\{\rho_m<a\}\subset\{\tilde\rho_m<a+\del_1\}.
$$
Therefore,
$$
\{\rho<-c_*\}\subset\{\tilde\rho<0\}\subset  \{\rho<\del_1\}\subset \{\rho<2\del_1\}\subset \{\rho_m<-\del_1\}\subset \{\tilde\rho_m<0\}.
$$

Set $\rho'=S_{\del_1}\rho$.  Then by \re{14t}
\eq{rhop}
\|\rho'-\rho\|_{C^0(X)} \leq\f{\del_1}{4}.
\end{equation}
Define
$$
\Omega'=\{\rho'<\frac{5\del_1}{4}\}, \quad \Om=\{\rho'<\f{7\del_1}{4}\}.
$$
Then \re{rhop} implies that $\{\rho<\del_1\}\subset \Om'$ and $\Om\subset \{\rho<2\del_1\}$.  Therefore,
\eq{Om*}
\{\rho<-c_*\}\subset\{\tilde\rho<0\}\subset  \{\rho<\del_1\}\subset  \Om'\subset\Om\subset    \{\tilde\rho_m<0\}.
\end{equation}

For the configuration  $(D_{p_j}^1, D_{r_2}^2,U_{p_j},\psi_{p_j})$ or $(D_{p_j}^1, D_{r_2}^2,D_{p_j}^3,U_{p_j}, \psi_{p_j})$ as in
Definitions~\ref{def:3.4}  or \ref{ccav}, we have the $\dbar$ solution operator $T_j$ satisfying $\dbar T_j(\psi_{p_{j+1}}^{-1})^*f_{j}=f_{j}$, where $f_0:=f$. Define
$$
f_{j+1}=f_j-\db\{\psi_{p_{j+1}}^*(\chi\cdot ( T_j(\psi_{p_{j+1}}^{-1})^*f_{j}))\}.
$$
Then $f_{m}\in \Lambda^r( D_m)$ is $\dbar$ closed on $D_{m}$.
We remark that $f\mapsto f_{m}$ is a linear operator $\mathcal G_D\colon  \Lambda^r( D)\cap \ker\db\to \Lambda^r( D_m)\cap\ker\dbar$, and $\cL G_D$ is independent of $r$.

We write
$$
f=\db u_m+f_{m}, \quad u_m=\sum_{j=1}^{m} \psi_{p_j}^*(\chi\cdot ( T_{j-1}\psi_{p_j}^{-1})^*f_{j-1}.
$$
Therefore, we  can focus   on the $\db$ equation for a {\it  fixed} $a_q$ domain $\Om$ with smooth boundary. It is important that $\Om\subset\widetilde D_m$ as long as $\|\tilde\rho-\rho\|_{C^2(X)}<\del_1$.
 We will apply \rt{3.4.6} to the domain $\widetilde D$ with $\del$ being as in Cor.~\ref{corb11}. In what follows, $\widetilde D, \widetilde D_m$ is denoted by $D, D_m$ respectively.

To ease notation, we write
$$
\|\cdot\|_{D;a}=\|\cdot\|_{C^a(\ov D)}, \quad|\cdot|_{D;a}=\|\cdot\|_{\Lambda^a(\ov D)},\quad \|\rho_j\|^-_{\Lambda^a}=|\rho_j|^-_a. $$
Recall that $\|\rho_j\|^-_{\Lambda^a}$ is defined in Definition~\ref{def|-}. We also need to estimates the norms for $f_m$.  We have for $r>1$
\al
|f_{j+1}|_{D_{j+1};r}&\leq C_r(D_j)(|f_j|_{D_{j};r}+ |\rho_j|^-_{r+2} \|f_j\|_{D_{j};1}),\\
\|f_{j+1}\|_{D_{j+1};1}&\leq C_r(D_j)(1+\|\rho_j\|^-_{3})\|f_j\|_{D_{j};1},\\
|u_{j+1}|_{D_{j+1};r+1/2}&\leq C_r(D_j)(|f_j|_{D_{j};r}+|\rho_j|^-_{r+5/2}\|f_j\|_{D_{j};1}).
\end{align}
By \re{rhoj=}, we have $|\rho_j|^-_{r+5/2}\leq C_{j,r}(1+|\rho|^-_{r+5/2})\leq 2C_{j,r}|\rho|^-_{r+5/2}$. Thus,
\ga\label{f_test}
|f_{m}|_{D_{m};r}\leq C_r(D)(|f|_{D_{0};r}+(\|\rho_j\|^-_{3})^m|\rho|^-_{r+2}\|f\|_{D_{0};1}),\\
 |u_{m}|_{D_m;r+1/2}\leq C_r(D)(|f|_{D_0;r}+(\|\rho_j\|^-_{3})^{m} |\rho|^-_{r+5/2}\|f\|_{D_0;1}).
\end{gather}

\medskip
\noindent{\it Step 2. Smoothing for interior regularity.}
To obtain the interior regularity, we will use regularity in Sobolev spaces. We need to avoid the loss in H\"older exponent from the Sobolev embedding. To this end, we will again use a partition of unity  to overcome the loss.
We can make $f_{m}$ to be $C^\infty$  on any relatively compact subdomain $U'$ of $\Om$  via local solutions as follows. Fix $x_0\in \ov{\Om'}$ where $\Om'$ is given in Step 1 with $D\Subset\Om'\Subset D_m$. We solve $\db u=f_{m}$ on an open set $ \om$ in $D_m$ that contains $x_0$. Let $\chi$ be a smooth function with compact support in $\om$ such that $\chi=1$ on a neighborhood $\om'$ of $x_0$. Then $\tilde f=f_{m}-\db(\chi u)=(1-\chi)f_{m}+\db\chi\wedge u$ is still in $\Lambda^r$, while $\tilde f=0$ on $\om'$. In particular, $\tilde f\in C^\infty(\om')$. Repeating this finitely many times, we can find $\tilde u\in\Lambda^{r+1}(\Om)$ with compact support in $D_m$ such that
\gan{}\label{tffm}
\tilde f=f_{m}-\db \tilde u\in C^\infty( \ov{\Om'}),\\
|\tilde f|_{\Om';r'}\leq C_{r'}|f_{m}|_{D_m;r}, \quad |\tilde u|_{\Om';r+1}\leq C_{r}|f_{m}|_{D_m;r}
\nonumber
\end{gather*}
for \emph{any}   $r'>r$. Recall that $\Om'$ still satisfies the condition $a_q$ and $\pd \Om'\in C^\infty$.

By \re{Om*}, we have $D_{-c_*}\subset \tilde D\subset\Om'\subset\Om\subset D_{2\del_1} \subset\widetilde D_{m}$.
Since $2\del_1<\del_*$,  Cor.~\ref{corb11} implies that $\db u=\tilde f_m$ admits a $L^2$ solution $u$ on $\Om'$ and hence it admits an $L^2$ solution $u$ on $\Om'$.  Furthermore, the defining function $\rho':=\rho'-\f{5\del_1}{4}$ of $\Om'$  satisfies $\|\rho''\|_a\leq C_a\|\rho\|_2$.

\medskip
\noindent{\it Step 3. An application of  canonical solutions.} We fix a smooth hermitian metric on $T^{(1,0)}X$ and a smooth hermitian metric on $V$. Let $\var=e^{\tau\rho}$,
 where $\rho(z)=-\dist(z,{\pd\Om'})$.
Let $L_{p,q}( \Om',V,\var)$ be the space of $V$-valued $(p,q)$ forms $f$ on $ \Om'$ such  that
$$
\|f\|_{\var}^2:=\int_ \Om |f(x)|^2e^{-\var(x)}\, \dv(x)<\infty,
$$
where $\dv$ is the volume form on $X$ with respect to the hermitian metric on $X$.  Write $\|f\|_\var$ as $\|f\|$ when $\var=0$.

We want to apply \cite{MR0179443}*{Thm.~3.4.6}, i.e. \rt{3.4.6} for the vector bundle version to
$
\Om'.
$
 Then $\var$ satisfies the condition $a_q$ on $\Om\setminus D_{-c*}$.
Let $\var_k=\chi_k(\var)$.
 Let $k_*$ to be the integer in \rt{3.4.6}   and let $\tilde\var=\var_{k_*}$.

 We now consider $T_{q}^{\tilde\var}=\db$ as densely defined from  $ L_{(0,q-1)}^2(\Om',\tilde\var)$ into $ L_{(0,q)}^2( \Om',\tilde\var)$ and $T_{q}^{\tilde\var *}$ its adjoint. Let $f_{m}$ be the $(0,q)$ form derived in Step 1.  By Corollary~\ref{corb11}, we find a solution $u_0$ satisfying
 $
 \db u_0=f_{m}
 $ on $\Om'$ and
 $$
 \|u_0\|_{\tilde\var}\leq C_*\|f_{m}\|_{\tilde\var}.
 $$

For the estimate, we need Kohn's canonical solutions.  By  \cite{MR0179443}*{Thm.~1.1.1}, $R_{(T^{\tilde\var}_q)^*}$ is also closed  and $R_{(T^{\tilde\var}_q)^*}=N_{(T^{\tilde\var}_q)^*}^\perp$.
We now apply the decomposition
\ga
u_0=u+h, \quad u\in N_{T^{\tilde\var}_q}^\perp, \quad h\in N_{T^{\tilde\var}_q}.
\end{gather}
  Thus, $u=(T^{\tilde\var}_q)^*v$ and $u$ satisfies
  \ga\label{du=f}
  \db u=\tilde f\quad \text{on $\Om'$},\\
  \label{uvark}
  \|u\|_{\tilde\var}\leq\|u_0\|_{\tilde\var}\leq C_*\|f_m\|_{\tilde\var}.
  \end{gather}
   In particular, in the sense of distributions, we have
  $
  \vartheta_q^{\tilde\var} v=u$ on $\Om'$. Here $\vartheta_q^{\tilde\var}$ is the formal adjoint (acting on test forms) of $\db_q$ in the $L^2$ spaces with weight $\tilde\var$.
Since $\vartheta^{\tilde\var}_{q-1}\vartheta^{\tilde\var}_q=0$, we get in the sense of distributions
 \eq{keyid}
  \theta^{\tilde\var}_{q-1}u=0  \quad \text{on $\Om'$.}
\end{equation}

Our last step is to  use the system of elliptic equations for $u$ to derive the {\it interior} estimates, using \re{du=f}-\re{keyid}. For its independent interest,   the last step follows from \rp{prop:int} below.
\end{proof}

\begin{prop}\label{prop:int} Let $\Om'$ be an open set in $X$.  Let $f\in  L^2_{(0,q-1)}(\Om',V,\tilde\var)$ and $u\in L^2_{(0,q-1)}(\Om',V,\tilde\var)$ with $\tilde\var\in C^\infty(\ov{\Om'})$. Suppose
\eq{keyid+}
\db u=f, \quad \theta^{\tilde\var} u=h  \qquad \text{on $\Om'$}.
\end{equation}
If $\Om''\Subset\Om'$, then for $2\leq p<\infty,0<\all<1$ and $k\in\nn$, we have
\al\label{int-1}
\|u\|_{W^{k+2,p}(\Om'')}&\leq C_{k,p}(\Om',\Om'')(\|u\|_{W^{0,p}(\Om')}+\|(f,h)\|_{W^{k+1,p}(\Om')});\\
\label{int-2}\|u\|_{C^{k+\all+1}(\Om'')}&\leq C_{k,\all}(\Om',\Om'')(\|u\|_{L^2(\Om')}+\|(f,h)\|_{C^{k+\all}(\Om')}),\quad k>1.
\end{align}
Here $C_{k,\all}(\Om',\Om''),C_{k,p}(\Om',\Om'')$  also depend on  $\|\var\|_{C^{k+3}(\Om')}$.
\end{prop}
\begin{proof}   Let $\|\cdot\|_{W^{k,p}}:=\|\cdot \|_{W^{k,p}(\Om')}$ denote the norm for space $W^{k,p}(\Om')$. Recall that
$T^{(1,0)}X$ and $V$ are endowed with smooth hermitian metrics.
Let $\om_1,\dots, \om_n$ be a local unit frame for $(1,0)$ forms and let  $e_1,\dots, e_m$ be the local unit frame of $V$.
Following notation in Appendix~\ref{sec:vb},
we have $u={\sum}^{\prime}u^\nu_{J}d\ov \om^J\otimes e_\nu$ and
\begin{align*}
\db u=\sum_\nu\sum_k{\sum_J}^{\prime}\f{\pd u^\nu_J}{\pd \ov\om_k}\ov \om_k\wedge \bar\om^J\otimes e_\nu+Ru,\\
\vartheta^{\tilde\var}_{q-1} u=-\sum_\nu\sum_j{\sum_K}^{\prime} \f{\pd u^\nu_{jK}}{\pd\om_j}  \ov \om^K\otimes e_\nu+Bu,
\end{align*}
where $Ru$ and $Bu$ are of order zero in $u$.
Write $\vartheta^\var$ as $\vartheta$ when $\var=0$.
Let $\chi$ be a  smooth function with support in a ball $B_R\Subset\Om'$ of radius $R$ centered at a point in $x_0\in\Om'$. Let $\tilde u=\chi u$. We abbreviate $(f,h)$ by $\tilde f$. We have
$$
\|\db\tilde u\|_{W^{0,2}}+\|\vartheta \tilde u\|_{W^{0,2}}\leq C(\|\tilde f\|_{W^{0,2}}+\|u\|_{W^{0,2}}).
$$
By~\cite{MR1045639}*{Lem.~4.2.3, p.~86},   for any relatively compact subset $\Om''$ of $\Om'$, we have
\eq{Uom}
\|u\|_{W^{1,2}(\Om'')}\leq C(\Om'',\Om')C_*(\|\tilde f\|_{W^{0,2}}+\|u\|_{W^{0,2}}).
\end{equation}

Recall the Sobolev compact embedding of $W^{j,q}$ in $W^{j+1,p}$ for $q=\f{p}{1-\f{p}{2n}}$ when $1\leq p<2n$.
For the following, we take $p_0=2$ and fix any $2<q<\infty$. We then fix $2<p_1<\cdots< p_{n^*}$ such that $p_{n^*-1}<2n$, $p_{j+1}\leq \f{p_j}{1-\f{p_j}{2n}}$ and  $p_{n^*}>q$.

Let $\Box_{\tilde\var}= \dbar_{q-1}\vartheta^{\tilde\var}_{q-1}+\vartheta_q^{\tilde\var}\dbar_{q}$. Then the principal part $\Delta_g=-\sum g^{jk}\f{\pd^2}{\pd z^j\pd \ov z^k}$ of $\Box_{\tilde\var}$ has smooth coefficients. Further, $\Delta_g$
is independent of $\tilde\var$,  diagonal and elliptic, where $g$ is the smooth hermitian metric $X$ (see~\cite{MR2109686}*{p.~154, p.~160} or Appendix~\ref{sec:vb}). In the sense of distribution, we have
$$
\Delta_gu= b(\nabla\tilde\var)\nabla \tilde f+c_1(\nabla\tilde\var)\nabla u+ c_0(\nabla^2\tilde\var)u,
$$
 where $b,c_1,c_0$ and $c_1',c_2$ below are matrices of polynomials whose coefficients depending on $g$. Then as a weak solution, $\tilde u$ satisfies $\Delta_g\tilde u= v$ for
\ga\label{Dgtu}
  v:=\chi [ b(\nabla\tilde\var)\nabla \tilde f+c_1(\nabla\tilde\var)\nabla u+ c_0(\nabla^2\tilde\var)u]
+c_1'(\nabla\chi))\nabla u+  c_2(\nabla^2\chi)u.
\end{gather}

Next,  we recall two interior estimates on systems of elliptic equations from Morrey~\cite{MR0202511}*{Thm 6.4.4., p.~246}:
\ga{}\label{sobe}
\|\tilde u\|_{W^{k+2,p} }\leq C_{k,p} \|v\|_{W^{k,p}}+C_R\|\tilde u\|_{W^{0,1}},\quad 1<p<\infty;\\
\|\tilde u\|_{k+2+\all}\leq C_{k,\all}\|v\|_{k+\all}+C_R\|\tilde u\|_{Lip}, \quad  \supp\tilde u\subset B_R
\label{holde}
\end{gather}
provided  the right-hand sides are finite. We may assume that $c_0<R<C_0$.

By the Sobolev inequality, $\|\tilde u\|_{W^{j,p_1}}\leq  C_j\|\tilde u\|_{W^{j+1,2}}$. Estimating the latter via \re{sobe} with $k=0$ and $p=2$ and \re{Uom}, we get  $\|\tilde u\|_{W^{1,p_1}}\leq C_2(\|v\|_{W^{0,2}}+\|u\|_{W^{0,2}})\leq C_2'(\|u\|_{W^{0,2}}+\|\tilde f\|_{W^{1,2}})$. Note that all constants depend on $\dist(\Om'',\pd\Om')$.  Repeating this bootstrapping argument, we can show that for any $\ell<\infty$, we have $\|\tilde u\|_{W^{1,\ell}}\leq C_\ell  (\|\tilde f\|_{W^{1,\ell}}+\|u\|_{W^{0,2}})$.
By \re{sobe} and the latter, we get $\|\tilde u\|_{W^{2,\ell}(\Om'')}\leq C_\ell (\|v\|_{W^{0,\ell}(\Om''')}+\|u\|_{W^{0,2}}) \leq C''_\ell (\|u\|_{W^{1,\ell}(\Om''')}+\|\tilde f\|_{
W^{1,\ell}})\leq C'''_\ell  (\|\tilde f\|_{W^{1,\ell}}+\|u\|_{W^{0,2}})$
for $\Om''\Subset\Om'''\Subset\Om'$. Repeating this for higher order derivatives, we obtain \re{int-1}.

Recall Sobolev inequality $C^{k,\all}\subset W^{k+1, \ell}$ for $\all:=1-\f{2n}{\ell}>0$.  We obtain
\eq{est1}
\|u\|_{\Om'',Lip}\leq C\| u\|_{\Om'';1+\all}\leq C_\all(\|\tilde f\|_{W^{1,\ell}}+\|u\|_{W^{0,2}})\leq C'_\all(\|\tilde f\|_{Lip}+\|u\|_{W^{0,2}}).
\end{equation}

Next, we prove by induction that
\eq{uk1}
\|u\|_{\Om'';k+1+\all}\leq C_{k,\all}'C_*(\|\var\|_{\Om';k+1+\all}\| \tilde f\|_{\Om';Lip}+\|\tilde f\|_{\Om';k+\all}+\|u\|_{W^{0,2}}).
\end{equation}
By \re{est1}, the above holds for $k=1$. Suppose the above hold and we want to verify it when $k$ is replaced by $k+1$. We have for $v$ in \re{Dgtu}
\aln{}
\|v\|_{\Om'';k+\all}&\leq C\{\|\tilde f\|_{\Om'';k+1+\all}
+C_*\|\var\|_{\Om'';k+1+\all}(\|\tilde f\|_{\Om'';0}+\|u\|_{\Om'';0} ) +\|u\|_{\Om'';k+1+\all}\}.
\end{align*}
By \re{uk1} and \re{est1}, we get $\|v\|_{\Om'';k+\all}\leq C\|\tilde f\|_{\Om';k+1+\all}
+2CC_*\|\var\|_{\Om'';k+1+\all}\|\tilde f\|_{\Om';Lip}+C\|u\|_{W^{0,2}}$. Then \re{holde}  yields
$$\| u\|_{\Om'';k+2+\all}\leq (2CC_*\|\var\|_{\Om';k+2+\all}+
C_\all')\|\tilde f\|_{\Om';Lip}+C\|\tilde f\|_{\Om';k+1+\all}+C\|u\|_{W^{0,2}}.
$$
This proves \re{uk1} and hence \re{int-2}.
  \end{proof}

  We conclude this section with an isomorphism theorem on cohomology groups with bounds.
  Define
\aln{}
H^{r,r+1/2}_{(0,q)}(\Om, V)&=\frac{\Lambda_{(0,q)}^r(\Om,V)\cap\ker\db}{\Lambda_{(0,q)}^r(\Om,V)\cap\db\Lambda_{(0,q-1)}^{r+1/2}(\Om,V)
}, \\ H^{r,loc}_{(0,q)}(\Om, V)&=\frac{\Lambda_{(0,q)}^r(\Om,V)\cap\ker\db}{\Lambda_{(0,q)}^r(\Om,V)\cap\db L^{2,local}_{(0,q-1)} (\Om,V)
}.
\end{align*}
\begin{thm}Let $r\in(1,\infty]$. Let $\Om$ be relatively compact  $a_q$ domain in $X$ and let $V$ be a holomorphic vector bundle on $X$. There exists $c>0$ such that the restriction $ H_{(0,q)}^{r,r+1/2}(\Om, V) \to  H^{r,loc}_{(0,q)}({\Om_{-c}}
, V)$ is an isomorphism for the following cases.
\bpp
\item     $q=1$ or $\pd \Om$ is strictly $(n-q)$ convex.
\item  $\pd \Om\in \Lambda^{r+5/2}$ and $r\in(1,\infty]$.
\epp
\eth
\begin{proof}The injectivity follows from the stability result proved in Appendix B. The surjectivity is obtained by the Grauert bumping method.
\end{proof}
We remark that using the $C^{1/2}$ local solution operators in~\cite{MR986248}*{Thms.~9.1 and 14.1}, we can also show that  the restriction
$$
\frac{C_{(0,q)}^0(\ov\Om,V)\cap\ker\db}{C_{(0,q)}^{0}(\ov\Om,V)\cap\db C_{(0,q-1)}^{1/2}(\ov \Om,V)
}
\to \frac{C_{(0,q)}^0(\ov\Om,V)\cap\ker\db}{C_{(0,q)}^{0}(\ov\Om,V)\cap\db L^{2,loc}_{(0,q-1)}( \Om,V)
}
$$
is an isomorphism,  provided $(a)$ $\Om$ is $(n-q)$ strictly convex   or $q=1$,   or $(b)$ $\Omega$  is $(q+1)$ concave with $\pd\Om\in C^{5/2}$.

\setcounter{thm}{0}\setcounter{equation}{0}

\section{Proof of \rt{regsol+} via a Nash-Moser method}\label{sec:NM}
In this section, we will prove \rt{regsol+} for $q\geq2$ when the domains have negative Levi eigenvalues by using the Nash-Moser smoothing operators. Our approach was inspired by a method of
 Dufresnoy~\cite{MR526786} for the $\db$-equation on a compact set that can be approximated from outside by strictly pseudoconvex domains of which the Levi eigenvalues are well controlled.
It is interested that V. Michel~\cite{MR1198845} showed that if the number of {\it non-negative} Levi eigenvalues of $\pd \Om$ is {\it exact} $n-q'$ near $z_0\in\pd\Om$, then   $\db u=f_{(0,q)}$ has a solution in $C^\infty(U\cap\Om)$ when $\pd\Om\in C^2$ for all $q\geq q'$. When $\pd\Om\in C^4$, there is (a possibly different) solution $u$  in $C^\infty(\ov\Om\cap U)$. For pseudoconvex domains with $C^2$ boundary in $\cc^n$, the $C^\infty$ regularity of $\db$ solutions under suitable assumptions on the Levi-form has been proved by Zampieri~\cites{MR1757879, MR1749685} and Baracco--Zampieri~\cites{MR2145559, MR2178735}. The reader is referred to the thesis of Yie~\cite{MR2693230} for the
global regularity of $\db$ solutions with $\pd D\in C^4$.  When $\pd D\in C^\infty$ additionally, the existence of $u$ was proved by Kohn~\cite{MR344703}.  Michel--Shaw~\cite{MR1675218}
 obtained smooth regularity of the $\db$ solutions on annulus domain $D_1\setminus \ov{D_2}$ where $D_1$ is a pseudoconvex domain with piecewise smooth boundary and $D_2$ is the intersection of bounded pseudoconvex domains.
 The $H^{s}$ solutions was proved by Harrington~\cite{MR2491606} for pseudoconvex domains $D$ with
 $\pd D\in C^{k-1,1}$,  $k>s+1/2, k\geq2$, and $s\geq0$.

 As observed in~\cite{MR3961327}, to the author's best knowledge it remains an open problem if $C^\infty(\ov D)$ solutions $u$ to $\db u=f$ exist on a bounded weakly pseudoconvex domain $D$ in $\cc^n$ with $C^2$ boundary.

We now state a detailed version of \rt{regsol+}.
\begin{thm}\label{nash-moser-c12} Let $q>1$. Let $D$ be a relatively compact domain with $C^s$ boundary in a complex manifold $X$ satisfying the condition $a_q$. Let $V$ be a holomorphic vector bundle on $X$.    Fix $\hat r$ as follows.
Let $m$ be as in \rt{regsol-full}, which is the number of times the Grauert bumping is used in Step 1 of the proof of \rt{regsol-full}.
\bppp
\item When $s=2$ and $r>r_0:=\max\{2m-1,s+5/2\}$,  fix $\hat r<r-5/2+\e$ for some sufficiently small $\e>0$.
\item When $s\geq 7/2$ and $r>s+5/2$, fix
$$\hat r<r+\f{1}{2}-\f{1}{2}\f{s}{s-1}\f{r-1}{r}.
$$
\eppp
There exists a linear $\db$-solution operator
$$  H_{q}^{r}\colon \Lambda_{(p,q)}^r( D, V)\cap\db L^2_{loc}(D,V)\to \Lambda_{(p,q-1)}^{r'}( D, V)
$$
satisfying
$
|H_q^{r}f|_{r'}\leq C_{r,r'}(D)|f|_r$ for all $r'\leq \hat r.
$
 Furthermore,
 \bpp \item  $C_{r,r'}(D)$ is stable under small $C^2$ perturbations of $D$, and $H_{q}^{r}f\in C^\infty(\ov D)$ if $f\in
  C_{(p,q)}^\infty(\ov D, V)\cap\db L^2_{loc}(D,V)$ additionally.
   \item The $r_0$ as the $m$ is stable under small $C^2$ perturbations of $\pd D$.
   \epp
\end{thm}
\begin{proof}
It suffices to prove the theorem when $r,r+1/2,r',r'+1/2$ are not integers, by replacing $r$ by a smaller number and $r'$ by a larger number. This allows us to identity the Zygmund spaces with the H\"older spaces for these orders. This also allows us to use the Taylor theorem and the global estimates for smooth domains.

Let $D$ be a $(q+1)$ concave domain with $C^{2}$ boundary defined by $\rho<0$. Denote by $D^a$ the domain defined by $\rho<-a$. Suppose that $f\in C^r$.
Using the Moser smoothing operator $S_t$, we
define
$$\tilde S_tu=S_tE_{D^\e}u.
$$
We remark that we need to use values of $u$ on $D^\e\Subset D$ and the value of  $t$ in $S_t$ is much larger than $\e$. Thus we need to use  $S_tE_{D^\e}$.
We have
\gan \label{tildeS}
\|\tilde S_tu-u\|_{D^\e,a}=\|S_tE_{D^\e}u-E_{D^\e}u\|_{D^\e,a}\leq C_{a,b}t^{b-a}\|u\|_{D^\e,b},\quad
a< b<L+a;\\
\label{tildeS+}
\|\tilde S_tu\|_{D^\e,b}\leq C_{a,b}t^{a-b}\|u\|_{D^\e,a}.
\end{gather*}
The $L$ is fixed for the rest of the proof.

Assume that $\pd D\in C^{s}$ with
$$
s\geq 2, \quad r\in(s+5/2,\infty), \quad t_{i+1}=t_i^d, \quad i\geq1
$$
with $d>1$ and $t_1\in(0,1)$ to be determined.
As in Yie~\cite{MR2693230}, we apply the above to the defining function $\rho$ of $D$ by setting
\eq{cstar}
\rho_{t_i}=\rho \ast\chi_{t_i}, \quad t_i=:(c_*\e_i)^{1/s}, \quad  D^{\e_i}_i=\{\rho _{t_i}<-\e_i\}, \quad D^{\e_i}=\{\rho<-\e_i\}.
\end{equation}
Here $c_*\in(0,1)$ is to be determined.
Also we need $\e_1=t_1^{s}/{c_*}<t_0$ for the $t_0$ in \rl{8.1}.  We now have
\al{}
\|\rho_{t_i}-\rho\|_{0}&\leq C_s t^s_i\|\rho\|_{s};\label{rho1-1-25} \\
\|\rho_{t_i}-\rho\|_{2}&\leq \e_*,\label{rho1-2c-25} \quad t_i<t^*(\nabla^2\rho,\e_*); \\
  \|\rho_{t_i}\|_b&\leq C_{b,s}t_i^{s-b}\|\rho\|_{s}, \quad b\geq s.
\label{rho1-3-25}
\end{align}
We emphasize that $t^*(\nabla^2\rho,\e_*)$ in \re{rho1-2c-25}
depends on the {\it modulus of continuity} of $\nabla^2\rho$.
Note that $t^{s}_1=c_*\e_1$. When $c_*$ is sufficiently small, we have by \re{rho1-1-25} and \re{14t+}
$$
\pd D^{\e_i}\subset D_i^{\e_i-c\e_i}\setminus D_i^{\e_i+c\e_i}, \quad  \pd D_i^{\e_i}\subset D^{\e_i-c\e_i}\setminus D^{\e_i+c\e_i}
$$
where $c\in(0,1/2)$ is independent of  $c_*$.
This however does not cause any difficulty for   domains satisfying the condition $a_q$, because the Levi eigenvalues do not decay towards the boundary. This is  decisive for the Nash--Moser iteration to succeed in our proof for $C^2$ domains.

We now determine $\e_*$ and $t_1$. Choose $\e_*$ so that $\tilde\rho:=\rho_{t_i}$ meet the requirements to apply the stability of constants in  \rt{3.4.6} and \rt{regsol-full}. We then fix $t_1<t_0$ so that $\e_1<t_0$ for $t_0$ in \rl{8.1}
and \re{rho1-2c-25} hold.

We will apply Taylor's theorem for H\"older spaces. Our estimates of gaining $1/2$ derivatives for the homotopy formula are for Zygmund spaces.  Therefore, in the following argument, we will work on H\"older spaces $C^a$ with $a\not\in\f{1}{2}\nn$.
By \re{rho1-3-25} and \rt{regsol-full} applied to the domain $D_1^{\e_1}$,  we have a $\db$ solution operator $ H_{q,D_1^{\e_1}}$ on $D_1^{\e_1}$ satisfying
$$
|H_{q,D_1^{\e_1}}\var|_{D_1^{\e_1}, r+1/2}\leq C_{r,\e^*}(\nabla\rho_{t_1},\nabla^2\rho_{t_1}) (|\var|_{D_1^{\e_1},r}+|\rho_{t_1}|_{r+5/2}|\rho_{t_1}|^m_{3+\e^*}
|\var|_{D_1^{\e_1},1+\e^*})
$$
where $\e^*>0$ is a small positive number.
Set $f_1=f$ and $v_1=H_{q,D_1^{\e_1}}f_1$ and $\tilde v_1=E_{D_1^{\e_1}}v_1$. Therefore, we obtain
\ga\label{v1t}
| v_1|_{D_1^{\e_1},r+1/2}\leq C_{r,\e^*}( |f_1|_{D_1^{\e_1},r}+t_1^{s-s_*-r}|f_1|_{D_1^{\e_1},1+\e^*}), \quad r\in(1,\infty).
\end{gather}
with
\eq{defs*}
s_*=\max\{0,5/2-s\}+m\max\{0,3+\e^*-s\}.
\eeq

Smoothing with a different parameter $\e_1$, we  define $w_1=S_{\e_1}\tilde v_1$ on $X$.  On $D$  define
$$
f_2=f_1-\db w_1.
$$
We iterate this.
  We also find a solution $v_i=H_{q,D_i^{\e_i}}f_i$ 
so that
$$
f_i=\db v_i,\quad \text{on $D_i^{\e_i}$}.
$$
Define $\tilde v_i=E_{D_i^{\e_i}}v_i$, $w_i=S_{\e_i}\tilde v_i$ and $f_{i+1}=f_i-\db w_i$. Then we want to show that  $u:=\sum w_i$ is the desired solution to $\db u=f_1$ on $D$, where $f_j$ tends to zero on $D$ as $j\to\infty$.

We have
$$
\db  v_2= f_1-\db w_1, \quad \text{on $D_2^{\e_2}$}.
$$
On $D_1^{\e_1}$, $\tilde v_1=E_{D_1^{\e_1}}v_1=v_1$ and hence $\db v_1=\db\tilde v_1$.  Therefore
$$
 f_2= f_1-\db w_1=\db(\tilde v_1-w_1)=\db(\tilde  v_1-S_{\e_1}\tilde v_1),\quad \text{on $D_1^{\e_1}$}.
$$
We have
$$
\|w_1-\tilde v_1\|_{D_1^{\e_1},a}=\|S_{\e_1}\tilde v_1-\tilde v_1\|_{D_1^{\e_1},a}\leq C_{b,a}\e_1^{b-a}\|\tilde v_1\|_{b}, \quad a< b<L+a.
$$
Hence by Taylor's theorem, we get for $b\geq1+\e^*$
\AL{\label{b2t}
\| f_2\|_{D,1+\e^*}&\leq C_b\sum_{k=1}^{[b]-1}\e_1^{k-1}\| f_2\|_{D_1^{\e_1},k}+C_b\e_1^{b-1-\e^*}\| f_2\|_{D,b}\\
&= C_b\sum_{k=1}^{[b]-1}\e_1^{k-1}\|\db\tilde v_1-\db w_1\|_{D_1^{\e_1},k}+C_b\e_1^{b-1-\e^*}\| f_2\|_{D,b}
\nonumber\\
&\leq C_b\sum_{k=1}^{[b]-1}\e_1^{k-1}C_b''\e_1^{b-k-\f{1}{2}}\|v_1\|_{D_1^{\e_1},b+
\f{1}{2}}+C_b\e_1^{b-1-\e^*}\| f_2\|_{D,b}\nonumber\\
&= C'_b \e_1^{b-\f{3}{2}}\|v_1\|_{D_1^{\e_1},b+\f{1}{2}}+C_b\e_1^{b-1-\e^*}\| f_2\|_{D,b}.\nonumber
}

Recall that $\pd D\in C^s$. By analogy of \re{v1t} and \re{b2t}, we have
\AL{\label{vii-25}
|v_{i}|_{D_i^{\e_{i}},r+1/2}&\leq C_{r,\e^*}(  | f_i|_{D_i^{\e_i},r}+t_i^{s-s_*-r}| f_i|_{D_i^{\e_i},1+\e^*}), \quad r\in(1,\infty)\\
\label{bi+-25}
\| f_{i+1}\|_{D,1+\e^*}
&\leq C_{r,\e^*} \e_i^{r-\f{3}{2}}\|v_i\|_{D_i^{\e_i},r+\f{1}{2}}+C_r\e_i^{r-1-\e^*}\| f_{i+1}\|_{D,r}.
}

To simplify notation, let $\|f_i\|_r=\|f_i\|_{D,r}$ and $\|v_i\|_r=\|v_i\|_{D_i^{\e_i},r}$. Thus the $r$-norm is estimated by
\gan\label{m-norm--25}
\| f_{i+1}\|_r\leq \| f_i\|_r+C_r\|w_i\|_{r+1}\leq \| f_i\|_r+C_r\e_i^{-1/2}\|v_i\|_{r+1/2}.
\end{gather*}
Therefore, by \re{vii-25}, we obtain
\AL{\label{m-norm-25}
\| f_{i+1}\|_{r}\leq  2 C_{r,\e^*} \e_i^{-1/2}(   \| f_i\|_{r}+ t_i^{s-s_*-r}\| f_i\|_{1+\e^*}), \quad r\in(1,\infty)^*.
}
Here $(1,\infty)^*=(1,\infty)\setminus \f{1}{2}\nn$. We exclude this discrete set of values since we need Taylor theorem for H\"older spaces while our estimates for the homotopy formula need Zygmund spaces.

We now define
\eq{Bi1-25}
B_{i+1}=\hat C^2_rt_{i}^{-s/2}B_{i},
\end{equation}
with $B_1\geq1$ being fixed and $\hat C_r\geq 1$  to be determined. Fix $r\in(s+\f{5}{2},\infty)^*$. Since our solution $u$ to $\db u=f_1$ is linear in $f$, by rescaling we may assume that
\eq{f1-B1+} \max\left\{\| f_1\|_r, t_1^{s-s_*-r}\| f_1\|_{1+\e^*},  \|v_1\|_{r+1/2}\right\}\leq 1.
\end{equation}

By induction, let us show that
\AL{\label{bi-25}
\| f_i\|_r&\leq   B_{i},\\
\label{bi0-25}
t_i^{s-s_*-r}\| f_i\|_{1+\e^*}&\leq\hat C_r  B_i,\\
\label{vi-25}
\|v_i\|_{r+\frac{1}{2}}&\leq \hat C^2_rB_i.
}
Suppose that the three inequalities hold. We want to verify them when $i$ is replaced by $i+1$.
Clearly, $\re{bi-25}_{i+1}$ follows from \re{m-norm-25}, \re{Bi1-25},  \re{bi-25}, and $\hat C_r\geq 4C_{r,\e^*}$.  By \re{cstar}, \re{bi+-25}, \re{vi-25} and $\re{bi-25}_{i+1}$, we obtain
\aln
t_{i+1}^{s-s_*-r}\| f_{i+1}\|_{1+\e^*}&\leq C_{r,\e^*} t_{i+1}^{s-s_*- r}\e_i^{r-\f{3}{2}}\|v_i\|_{r+\f{1}{2}}+C_{r,\e^*}
 t_{i+1}^{s-s_*- r}\e_i^{r-1-\e^*}\| f_{i+1}\|_r\\
&\leq C'_{r,\e^*} c_*^{\f{3}{2}-r}t_{i+1}^{ s-s_*-r}t_i^{sr-\f{3s}{2}}\hat C_r^2B_i+C_{r,\e^*}t_{i+1}^{ s-s_* -r}\e_i^{r-1-\e^*}B_{i+1}
\nonumber\\
\nonumber&= C_{r,\e^*}' c_*^{{\f{3}{2}-r}}t_{i+1}^{ s-s_*- r}t_i^{sr-s} B_{i+1}+C_{r,\e^*}c_*^{{\f{3}{2}-r}}t_{i+1}^{ s-s_*- r}t_i^{s(r-1-\e^*)}B_{i+1},
\end{align*}
which gives us $\re{bi0-25}_{i+1}$, provided  $\hat C_r\geq(C_{r,\e^*}'+C_{r,\e^*})c_*^{\f{3}{2}-r}$ and
\eq{condd-25}
-d(r+s_*-s)+s(r-1-\e^*)>0.
\end{equation}
The latter is assumed now.
Then $\re{vi-25}_{i+1}$ follows from $\re{vii-25}_{i+1}$,  $\re{bi-25}_{i+1}$ and $\re{bi0-25}_{i+1}$ and $\hat C_r\geq 2C_{r,\e^*}$. Therefore, it suffices to take
\eq{hCr}
\hat C_r:=\max\{ 4C_{r,\e^*},(C_{r,\e^*}'+C_{r,\e^*})c_*^{\f{3}{2}-r}\}.
\eeq

By interpolation, we get
$$
| f_i|_{(1-\theta)(1+\e^*)+\theta r}\leq C_{r,\theta}| f_i|_{1+\e^*}^{1-\theta}| f_i|_r^\theta\leq C_{r,\theta}\hat C_r t_i^{(1-\theta)(r+s_*-s)}B_i.
$$
By definitions, we have
$$
t_{i}=t_{1}^{d^{i-1}}, \quad
B_i=\hat C_r^{i-1}(t_1\cdots t_{i-1})^{-\f{s}{2}}B_1=\hat C_r^{i-1}t_1^{-\f{s}{2}\f{d^{i-1}-1}{d-1}}B_1\leq \hat C_r^{i-1}t_i^{-\f{s}{2(d-1)}}B_1.
$$
Therefore, $| f_i|_{1-\theta+\theta r}\leq C_{r,\theta}\hat C_r^{i}B_1t_i^{\la}$ converges rapidly if  \re{condd-25} holds and
$$
\la:=(1-\theta)(r+s_*-s)-\f{s}{2(d-1)}>0.
$$
By \re{vii-25}, we have for $\theta r-\theta\in(0,\infty)^*$
\begin{align*}\label{vii-+-25}
|v_{i}|_{D^{\e_{i}},3/2-\theta+\theta r}&\leq C_{r,\e^*}(   | f_i|_{D^{\e_i},1-\theta+\theta r}+ t_i^{s-s_*-1+\theta-\theta r}| f_i|_{D^{\e_i},1+\e^*})\\
&\leq   B_1C_{r,\e^*}C_{r,\theta}\hat C_r^{i} t_i^{\la}+C_{r,\e^*}\hat C_rt_i^{r-1+\theta-\theta r-\f{s}{2(d-1)}}B_1.
\nonumber\end{align*}
 This shows that
 \eq{vj}
 |v_i|_{\f{3}{2}-\theta+\theta r}\leq 2 B_1C_{r,\e^*}C_{r,\theta}\hat C_r^{i}t_i^{\la_*}
 \eeq
  converges rapidly if
\eq{}\label{mula}
\la_*:= (1-\theta)r_*  -\f{s}{2(d-1)}>0, \quad r_*:=r-1+\min\{0, s_*-s+1\}.
\end{equation}

We want to maximize $\theta r$ or $\theta$.
Now, $\la_*>0$ implies that
 $$
 \theta <1 -\f{s}{2(d-1)r_*}.
 $$
The latter is an increasing function of $d$. Assume
 \eq{m1d1-25}
1 -\f{s}{2(d-1)r_*}>0.
 \end{equation}
We now specify the parameters. Note that \re{condd-25} and \re{m1d1-25} are equivalent to
\eq{existd*}
1+\f{s}{2r_*}<d<d_*,\quad d_*:=\f{s(r-1-\e^*)}{r+s_*-s}.
\end{equation}
Under the above restriction on $d$,  $\f{3}{2}-\theta+\theta r$ has maximum value
$$
 \hat r=r+\f{1}{2}-\f{s(r-1)}{2(d_*-1)r_* }.
$$
Then $|v_i|_{r'}<C^i_{r'}t_i^{\la'}$ with $\la'>0$ where $\la'$ depends on $r'$ and $r'<\hat r$. Consequently the solution $u=\sum S_{t_i}E_{D^{\e_i}}v_i$ to $\db u=f$ is also in $\Lambda^{r'}( D)$.


We now compute the value of $\hat r$.
When $s=2$, we have $s_*=m(1+\e^*)$ and $r_*=r-1$. Thus for
\eq{r0value}
r>\max\{2m-2,s+5/2\},
\eeq we have $d_*>4/3$. Then we choose $\e^*$ sufficiently small so that \re{existd*} is satisfied and
$$
\hat r> r-\frac{5}{2}.
$$

Suppose that $s> 7/2$, $r>s+5/2$ and $\e^*$ is sufficiently small. We have $s_*=0$ and $r_*=r-s$.   One can check that there exists $d$ satisfying \re{existd*}.
Then
$$
\hat r= r+\f{1}{2}-\f{1}{2}\f{s}{s-1}\f{r-1}{r-\f{s}{s-1}\e^*}.
$$
Note that from $s\geq7/2$ and $r>s+5/2$, we see that
$$
\hat r>r-1/{12}
$$
Further, we still have $\hat r>r-1/{12}$ for $s=7/2$.

As stated in the theorem, the dependence of $H^{r,r'}$ on $r,r'$ for $r'<\hat r$ stated in  arises from the restriction for $r\in(1,\infty)^*$. When $r$ is already in $(1,\infty)^*$, we obtain a solution operator $H^{r,r'}$ depends only on $r$ for any $r'<\hat r$.

Finally, we show that if $f\in C^\infty(\ov D)$, then the constructed $u$ is smooth on $\ov D$. Fix $s\geq2$. For notations, we fix $r,r'$ and rename them as $r_0,r_0'$. Fix any $r>r_0$ satisfying $r,r+\f{1}{2}\not\in\nn$. We may assume \re{f1-B1+}, by the linearity of the solution $u$ in $f_1$. Thus we have found solutions $\sum w_i$ with
$$
\|w_i\|_0\leq C\hat C_r^{i}t_i^{\la_*}, \quad i=1,2,\dots.
$$
by \re{vj}.
Here for $\la_*$ we have fixed $d\in(1,2),\theta,t_0$ so that \re{condd-25} and \re{mula} in which $r$ is replaced by $r_0$ are satisfied and hence by the larger $r$.
  With \re{f1-B1+}, we have \re{bi-25}-\re{vi-25} for $i=1$.
  For the $\hat C_r$ defined by \re{hCr}, the  same proof shows that \re{bi-25}-\re{vi-25} hold for all $i$ since \re{condd-25} and \re{mula}-\re{m1d1-25} hold for the fixed $\theta,d$ (depending on $r_0,r_0'$) and the $r$.   This shows that $|v_i|_{1-\theta+\theta r+1/2}$ converges rapidly. Therefore, $u:=\sum w_i\in\Lambda^{1-\theta+\theta r+1/2}$. We conclude that $u\in C^\infty(\ov D)$.
\end{proof}

\appendix

\setcounter{thm}{0}\setcounter{equation}{0}
\section{Distance function to $C^2$ boundary}\label{sec:dist}

  The following is  proved in Li-Nirenberg~\cite{MR2305073} for $\pd\Om\in C^{k,\all}$ for $k\geq 2$ and $0<\all\leq1$ and stated in Spruck~\cite{MR2351645} for $C^2$ case and proved in Crasta-Malusa~\cite{MR2336304} for $C^2$ boundary.
We provide a proof here, including a stability property for our purpose. See Gilbarg--Trudinger~\cite{MR1814364} and Krantz--Park~\cite{MR614221} when domains are in $\rr^N$.
\begin{prop}[\cite{MR2305073}*{Thm.~1}, \cite{MR2351645}*{Prop.~4.1}]\label{LNS}
Let $s\in[2,\infty]$. Let $M$ be
a smooth Riemannian manifold. Let $\Om\subset M$ be a bounded domain   with $C^s$ (resp. $\Lambda^s$ with $s>2$) boundary. Let $\rho$ be the signed distance function of $\pd\Om$.  There is $\del=\del(\Om)>0$ so that $\rho\in C^s$ (resp. $\Lambda^s$ with $s>2$) in $B_\delta(\pd\Om)=\{x\in M\colon\dist(x,\pd\Om)<\del\}$. Furthermore, $\del(\Om)$ is upper-stable under small $C^2$ perturbations of $\pd\Om$.
\end{prop}
\begin{proof}
We are given a smooth Riemannian metric on $M$. Let $N$ be a subset of $M$.  Let $\gamma_p\colon [0,d]\to M$ be a geodesic connecting  $\gamma(0)=p\in M\setminus\ov N$ and $p^*=\gamma(d)\in\ov N$. Suppose that $\gamma$ is normal, i.e. $|\gaa_p'|=1$ and length $|\gamma_p|$ of $\gamma_p$ equals $\dist(p,N)$.  Then $$|\gamma_p([t,t'])|=t'-t, \quad \dist(\gamma_p(t),N)=d-t.
$$
Thus if $N$ is a $C^1$ hypersurface near $p^*\in M$, then $\gamma_p$ is orthogonal to $N$ at $p^*$.

\medskip

We recall some facts about   geodesic balls; see~\cite{MR1138207}*{Chap.~3, Sect.~4}:

$(i)$ For $p\in M$,  there is $0<r(p)\leq\infty$ so that the geodesic ball $B:=B_r(p)$, centered at $p$ with radius $r$, is strictly geodesic convex for $0<r<r(p)$. Specifically, any two points $p_0,p_1\in B$ are connected by a unique shortest normal geodesic in $M$ and  the geodesic is contained in $B$. Here the uniqueness is up to a reparameterization $t\to \pm t+c$.

\medskip
$(ii)$ $\pd B$ is a smooth compact hypersurface.

Let $K$ be a compact set in $M$. By a), we have $r(K):=\inf_{p\in K}r(p)>0$. We   cover $K$ by finitely many open sets $U_i$ and choose coordinate chart $x_i$ on $U_i$    such that  if a normal geodesic $\gamma$  is contained in   $K$ then $|(x_i\circ\gamma)^{(j)}(t)|\leq C_j(K)$ wherever $x_i\circ\gamma(t)$ is defined.
This implies that if  $p\in N$,  and $N\cap B$ is closed in $B$, then any point $q\in B_{r/2}(p)$ is connected to a point $q^*\in N\cap B$ via a geodesic $\gamma_q$ in $B$  with length  $\dist(q,N)$. However, $q^*$ may not be unique.

Let us use the above facts for $N=\pd\Om$. Set $r_0=r(\pd\Om)$.   Let $\nu$ be the unit inward normal of $\pd\Om$ with respect to the Riemann metric and choose sign so that $\rho<0$ in $\Om$. Then $\nu$ is $C^{s-1}$ on $\pd\Om$. Let $\gamma(t,p)$ be the  geodesic through $p\in\pd\Om$ with $\pd_t\gamma(0,p)=\nu(p)$ and $|t|<r_0$. Then $\gamma$ is $C^{s-1}$ on $(-r_0,r_0)\times\pd \Om$. Since $s\geq2$, $\pd_t\gamma(0,p)$ is non zero, normal to $\pd\Om$, and $\gamma$ fixes $\{0\}\times\pd\Om$ pointwise, then the Jacobian of $\gamma$ is non-singular. By the inverse mapping theorem,  there  exists a unique  solution $( \tilde\rho,P)\in\rr\times\pd\Om$ satisfying
\eq{gPx}
\gamma( \tilde\rho(x),P(x))=x,  \quad |\tilde\rho(x)|<r_1
\end{equation}
for $x\in B_{r_2}(\pd\Om)$. Here $0<r_2<r_1<r_0$ and $r_1,r_2$ are sufficiently small. Furthermore, $\tilde\rho,P$ are in $C^{s-1}(B_{r_2}(\pd\Om))$.

We want to show that $\tilde\rho=\rho$ on $B_{r_2}(\pd\Om)$. Fix $x\in B_{r_2}(\pd\Om)$ and take $x^*\in\pd\Om$ with $\dist(x,x^*)=\rho(x)(=\dist(x,\pd\Om))$. Since $r_2<r(\pd\Om)$, then $x$ is connected to the center $x^*$ by a geodesic $\tilde\gaa$  in the geodesic ball $B_{r_2}(x^*)$. Since $\dist(x,x^*)=\rho(x)$, then $\tilde\gamma$ is orthogonal to $\pd\Om$ at $x^*$.   Then $\tilde\gamma$ must be contained in the normal geodesic $\gamma(\cdot,x^*)$ with tangent vector $\nu(x^*)$. Next we choose the parametrization of $\tilde \gaa$ so that $\tilde\gamma'(0)=\nu(x^*)$. We get
$$
\gamma(\rho( x),x^*)=\tilde \gamma(\rho(x))=x, \quad\tilde\gamma(0)=x^*=\gamma(0,x^*),\quad
|\rho(x)|=\dist(x,x^*)<r_2.
$$
 By the uniqueness of solution to \re{gPx}, we conclude that $\tilde\rho(x)=\rho(x)$ (and $x^*=P(x)$, $|\tilde\rho(x)|<r_2$).

Next we verify $\rho\in C^s$. Recall that the vector field $X_1(x):=\pd_t\gamma(\rho(x),P(x))$ is $C^{s-1}$ in $x=\gamma(\rho(x),P(x))$.  Fix $x_0\in  B_{r_2}(\pd\Om)$.   In a small neighborhood $U$ of $x_0$ in $B_{r_2}(\pd\Om)$, we use the Gram--Schmidt orthogonalization to find pointwise linearly independent vector fields $X_2,\dots, X_n$ of class $C^{s-1}$ that are   orthogonal to $X_1$. We already know that $\rho=\tilde\rho\in C^{s-1}\subset C^{1}$. This allows us to compute a directional derivative of $\rho$ via any $C^1$ curve that is tangent to the direction.   Since $\gamma$ is  normal, then $X_1\rho=1$. Let $j>1$ and we want to show that $X_j\rho=0$. At $x\in U\setminus\pd\Om$, Gauss lemma says that $X_2,\dots, X_n$ are tangent to the smooth geodesic sphere  $\pd B_{|\rho(x)|}(P(x))$. We have $|\rho(x)|=\dist(x,\pd\Om)>0$ and
$$
|\rho(q)|=\dist(q,\pd\Om)\leq \dist(q,P(x))=|\rho(x)|, \quad \forall q\in\pd B_{|\rho(x)|}(P(x)).
$$
Hence on this geodesic sphere, the $C^{s-1}$ function  $\rho$ attains a local extreme  at $x$. This shows that $X_j\rho=0$ on $U\setminus\pd\Om$.   By continuity,  $X_j\rho(x)=0$ on $U$.
Therefore all $X_i\rho$ are $C^\infty$ functions on $U$. As observed by Spruck~\cite{MR2305073}, since all $X_i$ are $C^{s-1}$,  then  $\rho\in  C^{s}(U)$. Therefore, $\rho\in  C^{s}(B_{r_2}(\pd\Om))$.

Finally, the stability of $\dist(\cdot, \pd\Om)$ near $\pd\Om$ is a consequence of the geodesic equations of the second-order ODE system. We leave the details to the reader.
\end{proof}

\setcounter{thm}{0}\setcounter{equation}{0}
\section{Stability of $L^2$ solutions on pseudoconvex manifolds with $C^2$ boundary satisfying   condition $a_q$}\label{sec:vb}

Let $\Om$ be a relatively compact domain in a complex manifold $X$. Let $V$ be a holomorphic vector bundle on $X$. Let $f$ be a $V$-valued $(0,q)$ form on $\Om$.  Suppose that $\db u=f$ can be solved on $\Om$ for some $u\in L^2_{loc}(\Om)$ and $f=\tilde f+\db v$ for $v\in L_{loc}^2(\Om)$ and $\tilde f$ is a $\db$ closed form on a larger domain $\Om'$ containing $\ov\Om$. We want to know if there exists a neighborhood $\tilde\Om$ of $\ov\Om$, that is independent of $f$
 such that  $\tilde f=\db \tilde u$ for some $\tilde u\in L^2_{loc}(\tilde\Om)$. If such a domain $\tilde \Om$ exists, we say the solvability of the $\db$-equation on $\Om$ is {\it stable}.
This stability was proved by H\"ormander~\cite{MR0179443} when $\Om$ is an $a_q$ domain and $V$ is the trivial bundle and by Andreotti--Vesentini~\cite{MR0175148}*{Lem.~29, p.~122} for vector bundles on domains that are strictly $(n-q)$ convex with smooth boundary. For  completeness, we sketch a proof for the case of vector bundle.   We also take this opportunity to  relax the boundary condition $\pd\Om\in C^3$ to the minimum $C^2$ smoothness and we also formulate a stability for the $L^2$ bounds of the $\db$ equation on small $C^2$ perturbations of $a_q$ domains.
 For the reader's convenience, we will give our statements for the vector bundle case the references in~\cite{MR0179443}. 

We fix smooth hermitian metrics on $X$ and  $V$.    Cover $\ov \Om $ by finitely open sets $U_1,\dots, U_{m_0}$ of $X$. We assume that  each $U_j$ is biholomorphic to the unit ball in $\cc^n$
  by a coordinate map $z_j$ which is biholomorphic on $\ov{U_j}$.
 In what follows, we will denote by $U$ one of $U_1,\dots, U_{m_0}$ or their subdomains.

Let $\{e_1,\dots, e_m\}$ be a smooth unitary basis of $V$ on $U$. Let $u=\sum u^\mu  e_\mu\in C^1_{(p,q-1)}(\Om,V,loc)$. We have
$$
\db u=Au +Ru, \quad A(u^\nu e_\nu)=(\db u^\nu) e_\nu,
$$
where $Au^\nu=\db u^\nu$ is as in~\cite{MR0179443} for the scalar case, and $Ru$ involves no derivatives of $u$, i.e. $Ru$ is of  {\it order zero} in $u$ with smooth coefficients.
Therefore, the principal part $A$ of $\db$
 is locally {\it diagonal}. 
 This is an important property allowing the proofs   for the scalar case  can be adapted  to the vector bundles case without difficulty.

Let $\om_1,\dots, \om_n$ be unitary smooth $(1,0)$ forms on $X$. For a $V$-valued $(p,q)$-forms  $f=\sum f^\nu e_\nu$, define
$$
\left|\sum f^\nu e_\nu(x)\right|=\sum_\nu{\sum_{I,J}}^{\prime}|f^\nu_{I,J}(x)|^2, \quad f^\nu= {\sum_{I,J}}^{\prime}f^\nu_{I,J}\om^I\wedge\ov\om^J.
$$
The volume form on $X$  will be
$$
d\upsilon=
 (\sqrt{-1})^n\om^1\wedge\cdots\wedge \om^n\wedge\ov\om^1\wedge\cdots\wedge\ov\om^n.
$$
  For $q\geq0$, let $L_{(p,q)}^2(\Om ,V,loc)$ and $\Lambda^r_{(p,q)}(\Om ,V)$ be the spaces of $V$-valued $(p,q)$ forms of which the coefficients on $U$ are in $L^2_{loc}(\Om\cap U)$, $\Lambda^r({ \Om\cap  U})$, respectively. Let $\cL D^{(p,q)}(\Om ,V)$ be the space of  smooth $V$-valued $(p,q)$ forms of which the coefficients are in $\cL D(\Om )$, i.e. smooth functions with compact support. Let $\cL D'_{(p,q)}(\Om ,V)$ be the space of $V$-valued $(p,q)$ forms of which the coefficients are distributions in $\Om $.

If $\var$ is a real $L^\infty$ function in $\Om$, let $L_{(p,q)}^2(\Om ,V,\var)$ be the space  of sections of $V$-valued $(p,q)$ forms satisfying
$$
\|f\|_{\var}^2=\int_\Om |f(x)|^2e^{-\var(x)}\, \dv(x)<\infty.
$$
We will write $\jq{\cdot,\cdot}_\var$ for the induced hermitian product and $\|\cdot\|_\var$ for the norm on $L^2_{(p,q)}(\Om ,V, \var)$. The  norms are equivalent for all weights $\var\in L^\infty$.

Throughout the appendix, we assume  $q\geq1$. The operator $\db$ defines a  linear, closed, densely defined operator
$$
T\colon L^2_{(p,q-1)}(\Om, V ,\var)\to L_{(p,q)}^2(\Om ,V,\var)
$$
while $Tu=f$ holds if $\db u=f$ in the sense of distributions.
We abbreviate $
T=T_{q}, S=T_{q+1}.
$
We will write $T_q$ for $T$ if needed.  The domain $D_T$ and range $R_T$ are independent of  $\var\in L^\infty$.
For  $f\in L_{(p,q)}^2(\Om,V,\var)$, write  $v=T^* f$
if
$
\jq{u,v}_\var=\jq{\db u, f}_\var
$
for all $u\in D_{T}$.

 Throughout the section, we assume that $\pd\Om\in C^2$.
By \rp{LNS},
$\Om$ has a $C^2$ defining function $\rho$ in $X$ satisfying
$$
2|\db\rho|=1 \quad \text{on $\pd\Om$}.
$$
We also assume that   $\var\in Lip(\Om)$. Then  $D_{T^*}$ is independent of $\var$, while $R_{T^*}$  depends on $\var$.

\begin{rem}
As in~\cite{MR0179443},  $|T^*f|_\psi$ is  {\it always} referred to as the dual with respect to $\psi$, where $\psi$ will be chosen appropriately. For clarity, we write $T_\psi^*$ for $T^*$ when $\psi$ needs to be specified.
\end{rem}

Using integration by parts, we can verify that if $f\in C^1_{(p,q)}(\ov \Om,V )$ has compact support in $U\cap\ov \Om $, then  $f\in D_{T^*}$ if and only if
$$
\sum_{j=1}^nf^\nu_{I,jK}\DD{\rho}{\om^j}=0, \quad\text{on $U\cap\pd \Om $}, \quad \nu=1,\dots, m.
$$
Define
$
\cL D^1_{T^*}(\ov\Om):=C^1(\ov\Om,V)\cap D_{T^*}.$
We have from~\cite{MR0179443}*{p. 148}
$$
T^* f=B f +R^*f, \quad B f:=(B f^\nu)e_\nu,
\quad\text{on $U$}
$$
with
$
B f^\nu=-\sum_j{\sum_{I,K}}'\del_j f^\nu_{I,jK}\om^I\wedge\ov\om^K.
$
Thus $R^*$, {\it independent} of $\var$,  is an operator of the zero-th order with smooth coefficients, and the $B$ is {\it diagonal} and its principle part  is also independent of $\var$. Thus the boundary condition is principle and of order zero.

In summary, we have
\pr{density}Let $\Om$ be a relatively compact $C^2$ domain in $X$ and let $\var\in Lip(\Om)$.
%
Then   $D_{T^*}$ is independent of $\var$. Let  $\psi\in L^\infty(\Om)$ be  a  real function.
\bpp
\item For all  $f=\sum f^\nu e_\nu\in C_{(p,q)}^1(U\cap\ov\Om,V )$,
\ga\label{deco}
\left|\|Sf\|_\psi^2-\|Af\|_\psi^2\right|\leq C(\Om)\|f\|_\psi^2.
\end{gather}
\item
For all $f=\sum f^\nu e_\nu\in C_{(p,q)}^1(\ov \Om,V )\cap D_{T^*}$ with compact support in $U\cap\ov \Om $,
\ga\label{deco+} \left|\|T^*f\|_\psi^2- \|B f\|_\psi^2\right|\leq C(\Om)\|f\|_\psi^2.
\end{gather}
\item   $\cL D^1_{T^*}(\ov\Om)$ is  dense in $D_{T^*}\cap D_{S}$ w.r.t. the graph norm $ |f|_\psi+|Sf|_\psi+|T^*f|_\psi$.
 \epp
Further,   $C(\Om)$, independent of $\var,\psi$,    depends only on the diameter of $\Om$.
\epr
Here the last assertion follows from  \cite{MR0179443}*{p.~121}.
   We also have
\begin{prop}[\cite{MR0179443}*{eq.~(3.1.9)}]\label{MKH}  Let $\Om$ be a relatively compact $C^2$ domain in $X$. Let $\rho$ be the signed distance function of $\pd\Om$.  Let $\var\in C^{1,1}(\Om)$.
For all $f\in  C^1_{(p,q)}(\ov \Om,V )$ with compact support in $U\cap\ov\Om$, we have
\aln\label{t4}
\|Af\|_\var^2+\|B f\|_\var^2&=\sum_{\nu=1}^{m}\|Af^\nu \|_\var^2+\|B f^\nu \|_\var^2\\
&=\sum_{\nu=1}^{m}(Q_1+Q_2+t_1+t_2+t_3+t_4)
(f^\nu,f^\nu),\nonumber\end{align*}
with
\aln{}
\nonumber
Q_1(f^\nu ,f^\nu )&:=\sum_{I,J}\sum_j\int_{U\cap \Om }\left|\DD{f^\nu _{I,J}}{\ov\om^j}\right|^2e^{-\var}\, \dv,\\
\nonumber
Q_2(f^\nu ,f^\nu )&:= \sum_{I,K}\sum_{k,j}\int_{U\cap \Om }\var_{j\ov k}f^\nu _{I,jK}\ov{ f^\nu _{I,kK}}e^{-\var}\, \dv,
\\
t_1(f^\nu ,f^\nu )&:=\sum_{I,K}\sum_{k,j}\int_{U\cap\pd \Om }\left(f^\nu _{I,jK}\DD{\rho}{\om^j}\ov{\del _kf ^\nu_{I,kK}}-
f^\nu _{I,jK}\DD{\rho}{\ov\om^k}\ov{\DD{f^\nu _{I,kK}}{\ov\om^j}}\right)e^{-\var}\, \ds,\\
t_2(f^\nu ,f^\nu )&:=\sum_{I,K}\sum_{k,j}\int_{U\cap \Om }\left(f^\nu _{I,jK}\ov\sigma_j\ov{\del _kf ^\nu_{I,kK}}-
f^\nu _{I,jK}\sigma_k\ov{\DD{f^\nu _{I,kK}}{\ov\om^j}}\right)e^{-\var}\, \dv,\\
t_3(f^\nu ,f^\nu )&:=\sum_{I,K}\sum_{i,j,k}\int_{U\cap \Om } f^\nu _{I,jK}\ov c_{jk}^i\ov{\del_if^\nu_{I,kK}} \, \dv,\\
t_4(f^\nu ,f^\nu )&:=-\sum_{I,K}\sum_{i,j,k}\int_{U\cap\pd \Om } f^\nu _{I,jK} c_{jk}^i\ov{\del_if^\nu_{I,kK}} \, \dv.
\end{align*}
\epr
\begin{proof} The proof   for $f\in C^2$ and $\var\in C^2$ is in~\cite{MR0179443}. The case for $f\in C^1$ can be obtained by $C^1$ approximation of $C^2$ forms as in~\cite{MR0179443}*{p.~101}.
\end{proof}

Let $\var$ be a $C^2$ real function defined in a neighborhood of $z_0\in X$. Let $1\leq q<n$. Suppose $\var$ satisfies  the condition $a_q$ at $z_0$.
Let $\mu_1(z)\leq\mu_2(z)\leq\cdots\leq\mu_{n-1}(z)$ be the eigenvalues of the Levi form $L_{z_0}\var$ and let $\la_1(z)\leq\la_2(z)\leq\cdots\leq\la_n(z)$ be the eigenvalue of the hermitian form $H_\zeta\var(t):=\sum\DD{^2\var}{\om_j\ov\om_k}t_j\ov t_k$. The minimum-maximum principle for the eigenvalues says that
$$
\la_j(z)=\min_{\dim W=j}\left\{\max_{v\in W, |v|=1}\{H_z\var (v)\}\right\}.
$$
Thus $\la_1(z)\leq\mu_1(z)\leq\cdots\leq\mu_{n-1}(z)\leq \la_n(z)$.
 Let $r^-=\max(-r,0)$ for a real $r$. Then at $z_0$, the condition $a_q$ is valid if and only if
 $$
 \mu_1+\cdots+\mu_q+\sum_{j=1}^{n-1}\mu_j^->0.
 $$
  If $\psi<\psi(z_0)$ is strictly pseudoconvex at $z_0$ then $\psi$ satisfies the $a_q$ condition for $q=0,\dots, n-1$.
 As in~\cite{MR0179443}*{Def. ~3.3.2}, we say $\psi$ satisfies {\it the condition $A_q$} at $z_0$, if
$\nabla\var(z_0)\neq0$
 and
 $$
 \la_1(z_0)+\cdots+\la_q(z_0)+\sum_{j=1}^{n-1}\mu_j^-(z_0)>0.
 $$
When needed, we denote the above eigenvalues $\la_j,\mu_k$ by $\la_j(z_0,\var),\mu_k(z_0,\var)$.
Let us prove the following estimate for {\it weighted} eigenvalues.
\begin{lemma}[\cite{MR0179443}*{Lem.~3.3.3}]\label{3.3.3}
Suppose that $\var$ satisfies the condition $a_q$ at $\zeta$. Then $e^{\tau\var}$ satisfies the condition $A_q$. More specifically, there exist $c(\var)>0$ and $\tau_0(\var)$ such that for $
\tau>\tau_0$
\eq{newAq}e^{-\tau\var(\zeta)} \Bigl\{ \la_1(\zeta,e^{\tau\var})+\cdots+\la_q(\zeta,e^{\tau\var})+\sum_{j=1}^{n-1}\mu_j^-(\zeta,e^{\tau\var})\Bigr\}
>c(\var)\tau.
\end{equation}
Furthermore, $c(\var), \tau_0(\var)$ are stable under small $C^2$ perturbation of $\var$.
\ele
\begin{proof}The proof in~\cite{MR0179443} uses a proof-by-contradiction argument. For  stability, we need a direct proof. We have
$$\la_1(z_0)+\cdots+\la_q(z_0)+\sum_{j=1}^{n-1}\mu_j^-(z_0)\geq
\la_1(z_0)+\cdots+\la_q(z_0)+\sum_{j=1}^{n-1}\mu_j^-(z_0).
$$  For $t$, decompose $t\cdot\f{\pd}{\pd \zeta}=t'+t''$ where $t'$ is in   the complex tangent space $T'_\zeta \var$  and $t''$ is in its orthogonal complement.   We have
$$
\tilde H^\tau_\zeta(t):=\tau^{-1} e^{-\tau\var}He^{\tau\var}_\zeta(t)= H_\zeta\var(t)+
\tau|\pd\var(\zeta)|^2|t''|^2.
$$
Restricted on $T_\zeta'$, the above is still the Levi form $L_\zeta\var$ of which the eigenvalues are $\mu_1\leq\cdots\leq\mu_{n-1}$.
Let $\la_1(\tau),\dots,\la_n(\tau)$ be the eigenvalues of the above quadratic form. We still have
$\la(\tau)\leq\mu_1\leq\cdots\leq\mu_n\leq\la_n(\tau)$.  For any $\del>0$, we choose $\tau_0$ so that
$$
\tau|\pd\var(\zeta)|^2\del^2>\la_{n}(0)+1, \quad \forall\tau>\tau_0.
$$
Then $\la_1(\tau)\geq H^\tau_\zeta\var(t)\geq \mu_1-\e$ when $\del$ is sufficiently small. Analogously, we get $\la_j(\tau)\geq \mu_j-\e$ for $j=1,\dots, n-1$ when $\del$ is sufficiently small. We can choose $\e$ depending on
$
\mu_1+\cdots+\mu_q+\sum_{j=1}^{n-1}\mu_j^-$
 and modulus of continuity of $\partial^2\var$ to obtain \re{newAq}.
\end{proof}

\begin{thm}[\cite{MR0179443}*{Thm.~3.3.1}]\label{3.3.1}
Let $\Om$ be a relatively compact $C^2$ domain in $X$. Let $z_0\in\Om$. Suppose that $\var\in C^{2}(\ov\Om)$.  Then
\eq{eq3.3.1}
\tau \|f\|^2_{\tau\var}\leq C_\var^*\left\{   \|T_{\tau\var}^*f\|_{\tau\var}^2+\|S f\|^2_{\tau\var} + |f|_{\tau\var}^2\right\}
\end{equation}
holds for some neighborhood $U\subset\Om$ of $z_0$, some $C_\var,\tau_\var$, and all $\tau>\tau_\var$ and all $f\in C_{(p,q)}^1(\ov \Om,V ) $ with compact support in $U\cap\ov\Om$, if and only if the hermitian form $\sum\var_{jk}(z_0)t_j\ov t_k$ on $\cc^n$ has either at least $q+1$ negative or at least $n-q+1$ positive eigenvalues.
Furthermore, we can take
$$
C_\var^*=\f{C(\Om)}{\min_{z_0\in \Om\setminus\Om_c}(\sum_1^{n-1}\mu^-_j(z_0,\var)+\sum_{j=1}^q\mu_j(z_0,\var))}
$$
where $\mu_1(z_0,\var)\leq \cdots\leq \mu_{n-1}(z_0,\var)$ are eigenvalues of
$L_{z_0}\var$
with respect the hermitian metric on $X$, while $U$ depends on the modulus of continuity of $\partial^2 \var$. The constants $C_\var^*, \tau_\var$ are stable under $C^2$ perturbation of $\pd\Om$.
\eth
\begin{proof}Take any $g\in C^1_{(p,q)}(\ov\Om)\cap D_{T^*}$ with compact support in $U\cap\ov\Om$.
Apply \re{eq3.3.1}   to $f=ge_1$, which is actually proved in~\cite{MR0179443}  for the $g$; see (3.3.4)-(3.3.6) in~\cite{MR0179443}. By \re{deco} we get
\aln
\tau \|g\|^2_{\tau\var}&\leq C_\var(\|T^*(ge_1)\|_{\tau\var}^2+\|\db (ge_1)\|^2_{\tau\var})\\
&\leq C_\var(\|T^*g\|_{\tau\var}^2+\|\db g\|^2_{\tau\var}+C_1\|g\|^2_{\tau\var}),
\end{align*}
where $C_\var$ depends on the eigenvalues of $\var$ and $C_1$ is independent of $\tau$ and $\var$. Assume further that $\tau>2C_\var C_1$. Then we get \re{eq3.3.1} in which $f,C_\var$ are replaced by $g, 2C_0$.
Note that the constant $C$ in \re{deco} is independent of $\tau$. By \cite{MR0179443}*{Thm 3.3.1}, we get the eigenvalue condition. Assume that the eigenvalue condition holds. Then \re{eq3.3.1} holds when  $f$ is
replaced by $f^\nu$ for each $\nu$. By \re{deco} again, we get \re{eq3.3.1} by adjusting $\tau_0$ and $C_0$.
\end{proof}

\begin{thm}[\cite{MR0179443}, Thm.~$3.3.5$]\label{3.3.5}Let $\Om$ be a relatively compact $C^2$ domain in $X$. Let $\var$ satisfy the condition $A_q$ at $z_0\in\ov\Om$. If $z_0\in\pd\Om$ assume further that $\var<\var|_{\pd\Om}=\var(z_0)$ in $\Om$. Then there are a neighborhood $U$ of $z_0$
and a constant $C_\var^*$ such that for all convex increasing function $C^2$ function $\chi$ in $\rr$ we have
$$
\int\chi'(\var)|f|^2e^{-\chi(\var)}\, \dv\leq C^*_\var(\|T^*f\|_{\chi(\var)}^2+
\|\db f\|_{\chi(\var)}^2+\|f\|_{\chi(\var)}^2)
$$
for all $f\in C_{(p,q+1)}^1(\ov \Om,V )\cap D_{T^*}$ with compact support in $U\cap\ov\Om$.
\eth
\begin{proof}We apply the scalar version of the result as in the proof of \rt{3.3.1}. The proof in \cite{MR0179443} is valid via $C^1$ approximations for $f$.
\end{proof}

By partition of unity, the above yields the following.

\begin{prop}[\cite{MR0179443}*{Prop.~3.4.4}]\label{3.4.4} Let $\Om$ be a relatively compact $C^2$ domain in $X$.
 Let $\var<0$ in $\Om$ and vanish in $\pd\Om$ with $\var\in C^2(\ov\Om)$. Let $\Om_a=\{z\in\Om\colon\var(z)<a\}$. Suppose that $\var$ satisfies the condition $A_q$ in $\ov\Om\setminus\Om_{-c}$ for some $c>0$. Then there are a compact subset $K$ of $\Om_{-c}$ and a constant $C_\var^*$ such that for all convex increasing function $\chi\in C^2(\rr)$
\eq{eq3.4.2}
\int_{\Om\setminus K}\chi'(\var)|f|^2e^{-\chi(\var)}\, \dv\leq C^*_\var(\|T^*f\|_{\chi(\var)}^2+
\|S f\|_{\chi(\var)}^2+\|f\|_{\chi(\var)}^2)
\end{equation}
 holds for   all $f\in C^1_{(p,q)}(\ov\Om, V)\cap D_{T^*}$.
\epr

\begin{thm}[\cite{MR0179443}*{Thm.~3.4.1}]\label{3.4.1} Let $\Om$ be a relatively compact $C^2$ domain in $X$.  Suppose that $\pd\Om$ satisfies the condition $a_q$. Fix  $C^2$ defining function $\rho$ of
$\Om$ such that $\rho$ is the signed distance function to $\pd\Om$ and fix
 $\var=e^{\la\rho}-1$ with $\la$ sufficiently large. Then there exist compact subset $K$ of $\Om$ and  constant $\tau_\var$
  such that  if $\tau>\tau_\var$
  and $ f\in D_{S}\cap D_{T^*}\cap L^2_{p,q}(\Om,V)$  we have
\eq{eq3.4.1}
\int_{\Om\setminus K}|f|^2e^{-\tau\var}\leq\|T^*f\|_{\tau\var}^2+\|S f\|_{\tau\var}^2+\int_K|f|^2e^{-\tau\var}\, \dv.
\end{equation}
The latter implies that
 $R_{T}$ is closed and   finite codimension in $N_{S}$.
\eth
\begin{proof} 
 Here we need to go through the proof of  \cite{MR0179443}*{Thm.~3.4.1}.
There is a compact set $K$ in $\Om$ such that
$$
 \tau\int_{\Om\setminus K}|f|^2e^{-\tau\var}\, \dv\leq C^*_\var(\|T^*f\|^2_{\tau\var}+\|Sf\|^2_{\tau\var}+\|f\|^2_{\tau\var}),
$$
where $ C^*_\var$ is independent of $\tau$. The above is proved in \cite{MR0179443}*{Thm.~3.4.1} when $V$ is trivial. Thus it also holds for any $V$ by \re{deco} and \re{deco+}.
 We get
\re{eq3.4.1} for $f\in C^1(\ov\Om, V)\cap D_{T^*}$ when $\tau>2 C^*_\var$. By the density theorem, it holds for $f\in D_{T^*}\cap D_S$. The proof for the other direction in \cite{MR0179443}*{Thm.~3.4.1} is valid without any change.
 \end{proof}

So far, all the constants in the estimates are stable under $C^2$ perturbations of the domain $\Om$ and these constants are explicit to some extent. The next constant is however not explicit since it comes from a proof by contradiction. Nevertheless, it leads no essential difficulty in our applications.

Fix $\gaa>2C_\var^*$, where $C_\var^*$ is in \re{eq3.4.1}.
Let $\chi_k\in C^2$ be an increasing sequence of convex increasing functions such that
\eq{chik}
\chi_k(\tau)=\gamma\tau,\quad \text{when $\tau<-
 c
$}; \qquad \chi'_k(\tau)\to\infty, \quad \text{as $k\to\infty, \  \tau>-c$}.
\end{equation}
Set $\var_k=\chi_k(\var)$. Note that $\var_k\in C^2(\ov\Om)$.
Define
$$
N_{(p,q)}(\Om_{-c}, V,\gaa\var):=N_{S_c}\cap N_{T_c^*},
$$
where $T_c^*$ is the adjoint of $T_c=\db\colon L^2_{(p,q-1)}(\Om_{-c},V,\gaa\var)\to L^2_{(p,q)}(\Om_{-c},V,\gaa\var)$, while $S_c$ is the operator $\db\colon L^2_{(p,q)}(\Om_{-c},V,\gaa\var) \to L^2_{(p,q+1)}(\Om_{-c},V,\gaa\var)$.
We have the following.
\begin{prop}[\cite{MR0179443}*{Prop.~3.4.5}]\label{3.4.5} Fix $q>0$.
Let $\Om, \Om_{-c}, \var$ satisfy the hypotheses  in \rpa{3.4.4}.In particular, $\var$ satisfies the condition $A_q$ in $\ov\Om\setminus\Om_{-c}$.
 There exist constants $C_*$ and $k_*$, depending on $\var,c,\gaa$, and the sequence $\chi_k$ such that for $k>k_*$
\eq{eq3.4.4a}   \|f\|^2_{\var_k}\leq C_*(\|T^*f\|_{\var_k}^2+\|S f\|^2_{\var_k})
\end{equation}
provided $f\in D_{T^*}\cap D_{S}\cap L_{(p,q)}^2(\Om,V)$ with $q\geq1$ and
\eq{eq3.4.4b}
\int_{\Om_{-c}}\ip{f,h}e^{-\gaa\var}\, \dv=0, \quad \forall h\in N_{(p,q)}(\Om_{-c}, V,\gaa\var).
\end{equation}
Further, $k_*, c, C_*$ are stable under $C^2$ perturbation of $\pd\Om$   as follows: There is $\delta_*=\del_*(\nabla\rho,\nabla^2\rho)>0$, depending on $\rho,c$ such that if a real function $\tilde\varphi^0$ satisfies  $\|\tilde\var^0-\var\|_{C^2}<\delta_*$ then there exists a real function $\tilde\var$ satisfying the following
\bpp
\item $\|\tilde\var-\var\|_{2}<C(\|\var\|_2)\delta_*$
and
$$
\tilde\var=\var \quad \text{on $\Omega_{-c}$},\quad \tilde\var=\tilde\varphi^0 \quad\text {on $  \{z\in X\colon \tilde\varphi^0(z)>-c/2\}$}
$$
\item For $\tilde\Omega=\{z\in X\colon\tilde\rho^0<0\}$,  $\tilde\var_k=\chi_k(\tilde\var)$ and $k\geq k_*$, we have
\eq{eq3.4.4a-tilde}   \|f\|^2_{\tilde\var_k}\leq C_*(\|T^*f\|_{\tilde\var_k}^2+\|S f\|^2_{\tilde\var_k})
\end{equation}
provided $f\in D_{T_{\tilde\Omega}^*}\cap D_{S_{\tilde\Omega}}\cap L_{(p,q)}^2(\tilde\Om,V)$ with $q\geq1$ and $f$ satisfies \rea{eq3.4.4b}.
\epp
\epr
\begin{proof}
Fix a smooth function $\chi$ on $X$ such that $0\leq\chi\leq1$, $\chi=1$ on $\Omega_{-8c/9}$ and $\supp\chi\subset\Omega_{-7c/9}$. Define
$$
\tilde\var=\chi\var+(1-\chi)\tilde\var^0.
$$
Then $\tilde\var$ satisfies $(a)$.

 Fix $\gaa>2C^*_\var$ for the constant $C^*_\var$ in \rea{eq3.4.2}.

Assume that estimate \re{eq3.4.4a} is false.  Then we can find a sequence $\var_j$ of $C^2$ defining functions satisfying the following. For each $j$, we have the following.
\bppp
\item $\varphi_j\in C^2(X)$ and $\|\varphi_j-\varphi\|_2<\frac{1}{j}$,  $\Omega^j=\{\varphi_j<0\}$,   the operators $T_j,S_j$ for $\Omega^j$.
\item  $\Omega^j_{-c}=\{\var_j<-c\}=\Om_{-c}$.
\item   $f_j\in D_{T_j^*}\cap D_{S_j}\cap L_{(p,q)}^2(\Omega_j,V)$ with $q\geq1$ and for some $k_j>j$
\ga\label{eq3.4.5a}
\|f_j\|_{(\var_{j})_{k_j}}=1, \quad \|T_j^*f_j\|_{(\var_j)_{k_j}}+\| S_jf_j\|_{(\var_{j})_{k_j}}<\frac{1}{j},\\
\label{eq4.4.5aa}
\int_{\Om_{-c}}\ip{f_{j},h}e^{-\gaa \var}\, \dv=0, \quad \forall h\in N_{(p,q)}(\Om_{-c}, V,\gaa\var).
\end{gather}
\eppp
By the density theorem, we may assume that $f_j\in C_{(p,q)}^1(\ov {\Om^j},V )\cap D_{T_j^*}$, while \re{eq3.4.5a} still holds and \re{eq4.4.5aa} is, however, weakened to
\eq{eq3.4.5b}
\Bigl|\int_{\Om_{-c}}\ip{f_{j},h}e^{-\gaa\var}\, \dv\Bigr|<\f{1}{j}, \quad \forall h\in N_{(p,q)}(\Om_{-c}, V,\gaa\var), \quad \|h\|_{\Om_{-c},\gamma\var}\leq1.
\end{equation}
(Compare \re{eq3.4.5a} and \re{eq3.4.5b} with \cite{MR0179443}*{eq.~(3.4.5)}.)  Here we have used  $(\var_j)_{j}=\gamma\var$ on $\Om_{-c}$ and Cauchy-Schwarz inequality
\begin{align*}
\Bigl|\int_{\Om_{-c}}\ip{\tilde f_{j},h}e^{-\gaa\var}\, \dv\Bigr|&=\Bigl|\int_{\Om_{-c}}\ip{\tilde f_j-f_j,h}e^{-(\var_j)_{k_j}}\, \dv\Bigr|\\
& \leq C\left\{\int_{\Om^j}|\tilde f_j-f_j|^2e^{-(\var_j)_{k_j}}\, \dv\right\}^{1/2}
\end{align*}
for $\tilde f_j\in C_{(p,q)}^1(\ov {\Om^j},V )\cap D_{T_j^*}$. We still call the approximation  $\tilde f_j$ satisfying \re{eq3.4.5a} and \re{eq3.4.5b}  by $f_j$.

The rest of proof is based on some minor changes of the proof in~\cite{MR0179443}. We give details for completeness.
 The increasing function $\chi'$ is bounded below by $\gaa$. Let $K$ be any compact subset of $\Omega_{-c}$.  By \re{eq3.4.2} we have
\gan\label{eq3.4.6-0}
\int_{\Om^j\setminus K}|f_j|^2e^{-(\var_j)_{k_j}}\, \dv\leq C^*\gaa^{-1}(1+j^{-1}),
\\
1\geq\int_{  K}|f_j|^2e^{-\gamma\var}\, \dv\geq 1-C_\var^*\gaa^{-1}(1+j^{-1}).
\label{eq3.4.60}
\end{gather*}
Thus
\ga\label{eq3.4.6-}
\int_{\Om_{-c}\setminus K}|f_j|^2e^{-\gaa\var}\, \dv\leq C_\var^*\gaa^{-1}(1+j^{-1}),
\\
1\geq \int_{  K}|f_j|^2e^{-\gaa\var}\, \dv\geq 1-C_\var^*\gaa^{-1}(1+j^{-1}).
\label{eq3.4.6}
\end{gather}
By \re{eq3.4.6-}-\re{eq3.4.6} and the Banach-Alaoglu theorem,  we can find a subsequence of $f_j$, still called it $f_j$,  that converges weakly  to $f$ in $L^2_{(p,q)}(\Om_{-c}, V)$.
Hence, by \re{eq3.4.5b} and the weak convergence,  we get \re{eq3.4.4b}.
By \re{eq3.4.5a}, we have  $\db f_j$ tends to $0$ in $L^2$ norm on each compact subset of $\Om_c$. Therefore, $\db f=0$ on $\Om_{-c}$.

Take any function $\psi\in C_0^\infty(\Om_{-c})$ with $\psi=1$ on $K$.
Then by \re{eq3.4.5a},
$$
\|T_c^*(\psi f_j)\|_{\gaa\var}+\|S_c(\psi f_j)\|_{\gaa\var}
+\|\psi f_j\|_{\gaa\var}\leq C_\var.
$$
Further, by~\re{eq3.4.5a} and ~\cite{MR1045639}*{Lem.~4.2.3, p.~86}, we have
$\|Df_j\|_{K}\leq C'$ for all $f_j$. 
By Rellich--Kondrachov compactness theorem~\cite{MR2597943}*{p.~272}, there is a subsequence $f_j$ converges strongly to $f$ on $K$. Still denote the subsequence by $f_j$.
 By \re{eq3.4.6}, we get from $\gamma>2C_\var^*$
\eq{nozero}
\int_{  K}|f|^2e^{-\gaa\var}\, \dv=1-C_\var^*\gamma^{-1}>\f{1}{2}.
\end{equation}


Next we want to show that   $f\in N_{T_{c}^*}$.
Set $g_j=f_je^{-(\var_j)_{k_j}}$. By \re{chik} and the convexity of $\chi_{k_j}$, we have
$(\var_j)_{k_j}\geq\gamma\var_j\geq\gamma(\var-\f{1}{j})$. Thus, for any relatively compact subdomain $\Omega'$ of $\Omega$, $\Om'$ is contained in $\Omega^j$ for large $j$ and hence
\eq{norm-bound}
\int_{\Om'}|g_j|^2e^{\gaa\var}\, \dv\leq e^\gamma
\int_{\Om^j}|g_j|^2e^{(\var_j)_{k_j}}\, \dv\leq e^\gamma\|f_j\|^2_{(\var_j)_{k_j}}=e^{\gamma}.
\end{equation}
Fix $\Omega_i\Subset\Omega$ with $\Omega_{i}\subset\Omega_{i+1}$ and $\Omega=\cup\Omega_i$.
Take a subsequence $g_{k_{1,j}}$ of $g_j$ converging weakly to $g_{k_{1,\infty}}$ on $\Omega_1$. Inductively, take a subsequence $\{g_{k_{i+1,j}}\}_{j=1}^\infty$ of $\{g_{k_{i,j}}\}_{j=1}^\infty$ converging weakly to $g_{k_{j+1,\infty}}$ on $\Omega_{j+1}$.  Then $g_{k_{i+1,\infty}}=g_{k_{i,\infty}}$ on $\Omega_i$, which is denoted by $g$ now. Then $g_{k_{j,j}}$, still denoted by $g_j$,
 converges to $g$ weakly on $L^2(\Omega_i, V,\gamma\varphi)$ for each $i$. By the  weak convergence of $g_j$, we have  $\int_{\Om_i}|g|^2e^{\gaa\var}\, \dv
=\lim_{j\to\infty}\int_{\Om_i}\jq{g,g_j} e^{\gaa\var}\, \dv$. By \re{norm-bound} and the Cauchy-Schwarz inequality, we obtain $(\int_{\Om_i}|g|^2 e^{\gaa\var}\, \dv)^{1/2} \leq e^{\gaa/2}.
$
Hence $g\in L^2(\Omega, V,\gamma\var)$.

On $\Om_{-c}$, we have $g=fe^{-\gaa\var}$.  As $j\to\infty$,
\aln
\Bigl|\int_{\Om_i } \ip{g_j,\chi_{\Om_i\setminus\Om_{-c}}g}e^{\gaa\var}\, \dv\Bigr|&\leq\int_{\Om_i}|g_j|^2e^{(\var_j)_{k_j}}\, dv\Bigl [\int_{\Om_i\setminus\Om_{-c}}
|g|^2e^{2\gaa\var-(\var_j)_{k_j}}\, \dv\Bigr]^{1/2}\\
&\leq\Bigl [\int_{\Om_i\setminus\Om_{-c}}
|g|^2e^{2\gaa\var-(\var_j)_{k_j}}\, \dv\Bigr]^{1/2} \to0.
\end{align*}
Hence $
\int_{\Om_i\setminus\Om_{-c}}
|g|^2 e^{\gaa\var}\, \dv=\lim_{j\to\infty}\int_{\Om_i} \ip{g_j,\chi_{\Om\setminus\Om_{-c}}g}e^{\gaa\var}\, \dv=0.
$
Therefore,  $g=0$ on $\Om\setminus\Om_{-c}$.

In the sense of distribution, we have $\vartheta g_j=e^{-(\var_j)_{k_j}}T_j^*f_j$ where $\vartheta$ is defined by
$$
\int_{\Om^j}\ip{\vartheta g_j,u}\, \dv=\int_{\Om^j}\ip{g_j,\db u}\, \dv
$$
for all $u\in C^1_{(p,q-1)}(\Om^j, V)$ with compact support in $\Om^j$. Now
\eq{thgk}
\int_{\Om^j}|\vartheta g_j|^2e^{\gaa\var_j}\, \dv\leq
\int_{\Om^j}|\vartheta g_j|^2e^{(\var_j)_{k_j}}\, \dv=\|T_j^*f_j\|_{(\var_j)_{k_j}}^2\to0, \quad j\to\infty.
\end{equation}
For any $u\in C^1_{(p,q-1)}(\ov{\Om_{-c}}, V)$, we extend it to $\tilde u\in C^1_{p,q-1}(\Om\cup\Om^j, V)$ with compact support in $\Omega^j\cap\Om$ when $j$ is sufficiently large.
Recall that on $\Om_{-c}$, $g=fe^{-\gamma\var}$, $g=0$ on $\Om\setminus\Om_{-c}$.
Thus
\aln
\int_{\Om_{-c}}\ip{f,\db u}e^{-\gaa\var}\, \dv&=\int_{\Om}\ip{g,\db\tilde u}\, \dv=\lim_{j\to\infty}\int_{\Om}\ip{g_j,\db\tilde u}\, \dv.
\end{align*}
We have for large $j$
\aln
\int_{\Om}\ip{g_j,\db\tilde u}\, \dv
&= \int_{\Om^j}\ip{g_j,\db\tilde u}\, \dv=
\int_{\Om^j}\ip{\vartheta g_j,\tilde u}\, \dv.
\end{align*}
By \re{thgk} the last integral tends to $0$ as $j\to\infty$. This shows $\int_{\Om_{-c}}\ip{f,\db u}e^{-\gaa\var}\, \dv=0$.
Since $C^1_{(p,q-1)}(\ov{\Om_{-c}}, V)$ is dense in $D_{T_{c}^*}$, we conclude that $f\in N_{T_{c}^*}$.
By \re{eq3.4.5a}, $\db f_j\to 0$ in $L^2$ norm on any compact subset of $\Om_{-c}$.  Thus $\db f=0$ on $\Om_{-c}$.
Therefore,
$f\in N_{S_c}\cap N_{T_{c}^*}$.

Using the weak convergence of $f_j$ to $f$ on  $\Om_{-c}$   and letting $j\to\infty$  in \re{eq3.4.5b}, we obtain $f\perp N_{(p,q)}(\Omega_{-c}, V, \gamma\varphi)$.
Then $f=0$ on $\Om_{-c}$, which contradicts \re{nozero}.
\end{proof}

\begin{thm}[\cite{MR0179443}*{Thm.~3.4.6}]\label{3.4.6}Let $\Om$ be a relatively compact $C^2$ domain in $X$. Let $\Om_{-c}, \var, k,\var_k$ be as in \rpa{3.4.4}.  Let $V$ be a holomorphic vector bundle in $X$. Assume that $\var$ satisfies the condition $a_q$ in $\ov \Om\setminus\Om_{-c}$. There exist $k_*$ and $C_*$ satisfying the following.
\bpp\item If $f\in L^2_{(p,q)}(\Om,V)$
 and the equation
 $\db u_0=f$ has a solution $u_0$ in $L^2(\Om_{-c},V)$, then
 it has a solution $u$ in $L^2_{(p,q-1)}(\Om,V)$. In other words, the restriction $\ov H_{(p,q)}(\Om, V)\to\ov H_{(p,q)}(\Om_{-c}, V)$ is injective. Moreover,
$$
 \|u\|_{\var_k}\leq C_*\|f\|_{\var_k}, \quad k\geq k_*.
$$
 where $C_*, c, k_*$ are the constants in \rea{eq3.4.4a}.
 \item Furthermore,  $C_*,c,  k_*$ are  stable under $C^2$ perturbations of $\pd\Om$ in the following sense: Let $\Om,\tilde\Om$ be defined by $\rho<0, \tilde\rho<0$ respectively.  There exists $\del_*=\del_*(\nabla\rho,\nabla^2\rho)>0$, depending on $c,\rho$ such that  with $\var=e^{\la\rho}-1, \tilde\var=e^{\lambda\tilde\rho}-1$ for some $\lambda>0$ depending only on $\rho$ and $\|\tilde\rho-\rho\|_2<\del_*$,
 if $f\in L^2(\tilde\Om)$ is $\db$ closed and $f=\db u_0$ for some   $u_0\in L^2_{loc}(\Om_{-c})$, then there is a solution $u\in L^2(\tilde\Om)$ such that $\db u=f$ on $\tilde\Om$ and
 $$
  |u|_{\tilde\Om,\tilde\var_k}\leq C_*|f|_{\tilde\Om,\tilde\var_k}, \quad k\geq k_*.
  $$
  \epp
Here $C_*, k_*,\del_*$ are independent of  $\tilde\rho$ and  unknown.
\eth
\begin{proof}Assertion $(a)$ is proved in \cite{MR0179443}*{Thm.~3.4.6}, using \re{eq3.4.4a} and
\cite{MR0179443}*{Lem.~3.3.3 and Thm.~4.1.4}. The same proof is valid for $(b)$ via \re{eq3.4.4a-tilde}.
\end{proof}

Denote the $c$ in \rt{3.4.6} by $c_*$. Taking $\tilde\rho=\rho-a$, we obtain
\begin{cor}\label{corb11}Let $\rho,c_*,\del_*, \Omega$ be as above. Then the restriction
$\ov H_{(p,q)}(\Omega_a, V)\to\ov H_{(p,q)}(\Omega_{-c_*}, V)$ is injective for any $|a|<\del_*$.
\end{cor}

There is a detailed study in Lieb--Michel~\cite{MR1900133}*{Chapt.~VIII, Sect.~8} on the stability of estimates for the $\db$-Neumann operator on $\Om_c$, when $\Om$ is a strictly pseudoconvex manifold with smooth boundary. In our case, we must treat a slightly more general situation where $\tilde \Om$ can be any small $C^2$ perturbations of $\Om$.
We do not know if $C_*, k_*$ are stable under $C^2$ perturbations in the sense of \re{defupst}; nevertheless, using Grauert's bumping method for $a_q$ domains, Corollary~\ref{corb11} suffices our purposes.


\newcommand{\doi}[1]{\href{http://dx.doi.org/#1}{doi:#1}}
\newcommand{\arxiv}[1]{\href{https://arxiv.org/pdf/#1}{arXiv:#1}}

\bibliographystyle{alpha}



\end{document}